\documentclass[10pt]{amsart}
\usepackage{graphicx}
\usepackage{verbatim}
\usepackage[square]{natbib}

\font\elevenss=cmss11

\font\eightss=cmss8

\font\sixss=cmss8 at 6pt

\newfam\ssfam
\textfont\ssfam=\elevenss \scriptfont\ssfam=\eightss
 \scriptscriptfont\ssfam=\sixss
\def\ss{\fam\ssfam \elevenss}%

\newtheorem{thm}{Theorem}[section]

\newtheorem{pr}[thm]{Proposition}
\newtheorem{cor}[thm]{Corollary}
\newtheorem{defn}[thm]{Definition}
\newtheorem{ass}[thm]{Standing assumption}
\newtheorem{example}[thm]{Example}

\newtheorem{alg}[thm]{Algorithm}
\newtheorem*{unremark}{Remark}

\newcommand{\Em}[1]{\textbf{#1}}

\newcommand{\noproof}{\hfill \qedsymbol}

\def\disp{\displaystyle}
\def\grad{\nabla}
\def\ee{\epsilon}

\def\vv{{\bf v}}
\def\xx{{\bf x}}
\def\yy{{\bf y}}
\def\zz{{\bf z}}
\def\pp{{\bf p}}
\def\mm{{\bf m}}

\def\ww{{\bf w}}
\def\rr{{\bf r}}
\def\sss{{\bf s}}
\def\nn{{\bf n}}
\def\vb{{\bf b}}
\def\P{{\mathbb P}}

\def\I{{\bf I}}
\def\Y{{\bf Y}}
\def\V{{\bf V}}
\def\L{{\bf L}}
\def\J{{\bf J}}
\def\simplex{{\Delta^{d-1}}}
\def\torus{{\bf T}}
\def\disk{{\bf D}}
\def\poly{{\bf D}}
\def\patch{\mathcal{P}}
\def\cone{{\bf K}}
\def\linear{{\bf L}}
\def\CC{\mathcal{C}}

\def\sing{\mathcal{V}}

\def\denom{H}
\def\numer{G}
\def\dom{\mathcal{D}}
\def\domain{\mathcal{D}}
\def\logdom{{\rm log} \dom}

\def\dirset{{\bf \Xi}}
\def\one{{\bf 1}}
\def\zero{{\bf 0}}
\def\rhat{\hat{\rr}}
\def\rbar{\overline{\rr}}
\def\shat{\hat{\sss}}

\def\stratum{\mathcal{S}}
\def\nbd{\mathcal{N}}

\def\dir{{\ss dir}}
\def\Cox{\hfill \Box}
\def\C{{\mathbb C}}
\def\Z{{\mathbb Z}}
\def\N{{\mathbb N}}
\def\R{{\mathbb R}}
\def\Q{{\mathbb Q}}
\def\CP{{\mathbb C}{\mathbb P}}
\def\RP{{\mathbb R}{\mathbb P}}
\def\orthant{\overline{\mathcal{O}}}
\def\Real{{\rm Re\,}}
\def\Imag{{\rm Im\,}}

\def\M{\mathcal{M}}

\def\contrib{{\mbox{\ss contrib}}}
\def\connect{{\mbox{\ss connect}}}
\def\crit{{\mbox{\ss crit}}}
\def\alg{{\mbox{\ss alg-log}}}

\def\loggrad{{\grad_{\log}}}
\def\logsing{\log \sing}
\def\ft{{\tilde{f}}}
\def\hess{\mathcal{H}}
\def\homot{{\bf H}}
\def\eul{{\hat{A}}}
\def\diag{{\rm diag \,}}
\def\trace{{\rm trace \,}}
\def\map{\Psi}
\def\Eta{\Theta}
\def\btheta{{\mathbb \theta}}
\def\normalB{\mathfrak{n}_B}

\newcommand{\var}{\mathcal{W}}

\makeatletter \def\romenumi{ \def\theenumi{\roman{enumi}}
\def\p@enumi{\theenumi} \def\labelenumi{(\@roman\c@enumi)}} 
\makeatother

\begin{document}
\renewcommand{\thepage}{\roman{page}}

\begin{titlepage}
\begin{center}
{\large \bf Twenty combinatorial examples of asymptotics 
derived from multivariate generating functions} \\
\end{center}
\vspace{5ex}
\begin{flushright}
Robin Pemantle \footnote{Research supported in part by
National Science Foundation grant \# DMS 0103635}$^,$\footnote{University
of Pennsylvania, Department of Mathematics, 209 S. 33rd Street, Philadelphia,
PA 19104 USA, pemantle@math.upenn.edu} ~\\
Mark Wilson \footnote{Department of Computer Science, University of Auckland, 
Private Bag 92019, Auckland, NEW ZEALAND, mcw@cs.auckland.ac.nz} 
\end{flushright}

\vfill

\noindent{\bf ABSTRACT:} \hfill \\[1ex]

Let $\{ a_\rr : \rr \in \N^d \}$ be a $d$-dimensional array
of numbers, for which the generating function $F(\zz) := 
\sum_\rr a_\rr \zz^\rr$ is meromorphic in a neighborhood of 
the origin.  For example, $F$ may be a rational multivariate 
generating function.  We discuss recent results that allow the 
effective computation of asymptotic expansions for the 
coefficients of $F$.

Our purpose is to illustrate the use of these techniques 
on a variety of problems of combinatorial interest.  The survey
begins by summarizing previous work on the asymptotics of 
univariate and multivariate generating functions.  Next we
describe the Morse-theoretic underpinnings of some new asymptotic
techniques.  We then quote and summarize these results in such 
a way that only elementary analyses are needed to check hypotheses 
and carry out computations.  

The remainder of the survey  focuses on combinatorial applications, 
such as enumeration of words with forbidden substrings, edges and 
cycles in graphs, polyominoes, and descents in permutations.  After the 
individual examples, we discuss three broad classes of examples, 
namely functions derived via the transfer matrix method, those 
derived via the kernel method, and those derived via the method of 
Lagrange inversion.  These methods have the property that 
generating functions derived from them are amenable to our 
asymptotic analyses, and we describe further machinery that 
facilitates computations for these classes of examples.
\vfill



\end{titlepage}

\tableofcontents
\pagebreak
\begin{center} {\bf List of generating functions by subsection} \end{center}
{\large
\begin{tabular}{lcl}
\ref{ss:binomial} &~~~~~~~~~~~~~~& $\frac{1}{1-x-y}$ \\[2ex]
\ref{ss:delannoy} && $\frac{1}{1-x-y-xy}$ \\[2ex]
\ref{ss:distinct subseq} && $\frac{1}{1-x-xy(1-x)^d}$ \\[2ex]
\ref{ss:fine} && $ \frac{2}{ 1 + 2x + \sqrt{1-4x} - 2xy } $
   \\[2ex]
\ref{ss:polyomino} && $\frac{xy(1-x)^3}{(1-x)^4 - xy(1-x-x^2+x^3+x^2y)}$ \\[2ex]
\ref{ss:euler} && $\frac{e^x - e^y}{xe^y - ye^x}$ \\[2ex]
\ref{ss:smirnov} && $\frac{1}{1-\sum_{j=1}^d \frac{z_j}{1+z_j}}$ \\[2ex]
\ref{ss:alignments} && $\frac{1}{2 - \prod_{j=1}^k (1+z_j)}$ \\[2ex]
\ref{ss:edge} && $\frac{1}{1-z(1-x^2 y^2)} \frac{1}{1-x(1+y)}$ \\[2ex]
\ref{ss:queueing} && $\frac{\exp (x+y)} { (1 - \rho_{11} x - \rho_{21} y)
   (1 - \rho_{12} x - \rho_{22} y) }$ \\[2ex]
\ref{ss:switching} && $\frac{P(u,v,z)}{1 + z^3 v + u z^3 + u z^2 + u z^2v + z^2 v - uz - 2z - zv -  u z^4 v}$ \\[2ex]
\ref{ss:connector} && $\scriptstyle{\left [ 1 - (x_1 + \cdots + x_d)
   - \trace \left ( (\I + \V)^{-1} \L \J \right ) \right ]^{-1}}$ \\[2ex]
\ref{ss:forbidden} && $\frac{1 + x^2 y^3 + x^2 y^4 + x^3 y^4 - x^3 y^6}
   {1 - x - y + x^2 y^3 - x^3 y^3 - x^4 y^4 - x^3 y^6 + x^4 y^6}$ \\[2ex]
\ref{ss:core} && $\frac{3xz(1-z)(3-z)}{(1-3y(1+z)^2) (27 -xz (3-z)^2)}$ \\[2ex]
\ref{ss:random walk} && $\frac{1}{1 + \sqrt{1-x} - y}$ \\[2ex]
\ref{ss:paths} && $\frac{2}{1 + \sqrt{1 - 4 x^2} - 2xy}$ \\[2ex]
\ref{ss:paths} && $\frac{2} {1 + \sqrt{1 - 2x - 3x^2} - x - 2xy}$ \\[2ex]
\ref{ss:paths} && $\frac{2} {1 + \sqrt{1 - 6 x^2 + x^4} - x^2 - 2xy}$ \\[2ex]
\ref{ss:pebble} && $\frac{\eta (x,y)}{x - (x+y)^2}$ \\[2ex]
\ref{ss:diagonal} && $\frac{1 + xy + x^2 y^2}{1-x-y+xy-x^2y^2}$
\end{tabular}
}
\pagebreak

\renewcommand{\thepage}{\arabic{page}}
\setcounter{page}{1}

\setcounter{equation}{0}
\section{Introduction} \label{ss:intro}

The purpose of this paper is to review recent developments
in the asymptotics of multivariate generating functions, and
to give an exposition of these that is accessible and centered
around applications.  The introductory section lays out notions 
of generating functions and their asymptotics, and delimits the 
scope of this survey.

We employ the standard asymptotic notation: ``$f = O(g)$''
means $\limsup_{z \to z_0} |f(z) / g(z)| < \infty$, where the
limit $z_0$ is specified but does not appear in the notation;
``$f = o(g)$'' means $f/g \to 0$ and ``$f \sim g$'' means 
$f/g \to 1$, again accompanied by a specification of which
variable is bound and the limit to which it is taken.
When we say ``$a_n = O(g(n))$'' we always mean as $n \to \infty$;
for multivariate arrays and functions, statements such as
``$a_\rr = O(g(\rr))$'' must be accompanied by a specification of
how $\rr$ is taken to infinity.  We typically use $\rr$ for
$(r_1 , \ldots , r_d) \in \Z^d$ and $\zz$ for $(z_1 , \ldots , z_d) 
\in \C^d$, but sometimes it is clearer to use $(r,s,t)$ for
$(r_1, r_2, r_3)$ or $(x,y,z)$ for $(z_1, z_2, z_3)$.  The norm
$|\rr|$ of a vector $\rr$ denotes the $L^1$ norm $\sum_j |r_j|$.

\subsection{Background: the univariate case}

Let $\{ a_n : n \geq 0 \}$ be a sequence of complex numbers and let
$f(z) := \sum_n a_n z^n$ be the associated generating function.
For example, if $a_n$ is the $n^{th}$ Fibonacci number, then
$f(z) = 1 / (1-z-z^2)$.
Generating functions are among the most powerful tools 
in combinatorial enumeration.  In the introductory section
of his graduate text~\cite{Stan1997}, Richard Stanley
deems a generating function $f(z)$ to be ``the most useful
but most difficult to understand method for evaluating'' $a_n$,
compared to a recurrence, an asymptotic formula and a complicated
explicit formula.  

There are two steps involved in using generating functions to
evaluate a sequence: first, one must determine $f$ from the
combinatorial description of $\{ a_n \}$, and secondly one
must be able to extract information about $a_n$ from $f$.  
The first step is partly science and partly an art form.  Certain
recurrences for $\{ a_n \}$ translate neatly into functional equations
for $f$ but there are numerous twists and variations.  Linear
recurrences with constant coefficients, such as in the Fibonacci
example, always lead to rational generating functions, but often
there is no way to tell in advance whether $f$ will have a sufficiently
nice form to be useful.  A good portion of many texts on 
enumeration is devoted to the battery of techniques available
for producing the generating function $f$.

The second step, namely estimation of $a_n$ once $f$ is known,
is reasonably well understood and somewhat mechanized.  
Starting with Cauchy's integral formula
\begin{equation} \label{eq:cauchy0}
a_n = \frac{1}{2 \pi i} \int z^{-n-1} f(z) \, dz \, ,
\end{equation}
one may apply complex analytic methods to obtain good estimates
for $a_n$.  As a prelude to the multivariate results which are
the main subject of this survey, we will give a quick primer 
on the univariate case (Section~\ref{ss:univariate} below), 
covering several well known ways to turn~(\ref{eq:cauchy0})
into an estimate for $a_n$ when $n$ is large.

\subsection{Multivariate asymptotics} \label{ss:GF-seq}

In 1974, when Bender published the review 
article~\cite{Bend1974}, the extraction of 
asymptotics from multivariate generating functions was 
largely absent from the literature.  Bender's concluding section 
urges research in this area:
\begin{quote}
Practically nothing is known about asymptotics for recursions in two 
variables even when a generating function is available.  Techniques for
obtaining asymptotics from bivariate generating functions would be quite
useful.
\end{quote}

By the time of Odlyzko's 1995 survey~\cite{Odly1995}, a
single vein of research had appeared, initiated by 
Bender and carried further by Gao, Richmond and others.
Let the number of variables be denoted by $d$ so that $\zz$ denotes the
$d$-tuple $(z_1 , \ldots , z_d)$. We use the multi-index notation
for products: $\zz^\rr := \prod_{j=1}^d z_j^{r_j}$.  Suppose
we have a multivariate array $\{ a_\rr : \rr \in \N^d \}$
for which the generating function
$$F(\zz) = \sum_\rr a_\rr \zz^\rr$$
is assumed to be known in some reasonable form.  We are interested
in the asymptotic behavior of $a_\rr$.  For concreteness, one may
keep in mind the following examples, which are discussed in detail
in Sections~\ref{ss:binomial} and~\ref{ss:delannoy}.
\begin{example}[binomial coefficients] \label{eg:pascal}
Let $a_{rs} = \binom{r+s}{r,s}$.  The array $\{ a_{rs} \}$
is Pascal's triangle, oriented so that the rays of ones emanate
from the origin along the positive $r$ and $s$ axes.  From the recursion
$a_{rs} = a_{r-1,s} + a_{r,s-1}$, holding whenever $(r,s) \neq (0,0)$,
one easily obtains 
$$F(x,y) := \sum_{r,s \geq 0} a_{rs} x^r y^s = \frac{1}{1-x-y} \, .$$ 
\end{example}
\begin{example}[Delannoy numbers] \label{eg:delannoy}
Let $a_{rs}$ count the number of paths from the origin to the lattice
point $(r,s)$ made of three kinds of steps: 1 unit east, 1 unit north,
and $\sqrt{2}$ units northeast.  These are called \Em{Delannoy
numbers}. The recursion $a_{r,s} = a_{r-1,s} + a_{r,s-1} +
a_{r-1,s-1}$ leads immediately to 
$$
F(x,y) := \sum_{r,s \geq 0} a_{rs} x^r y^s = \frac{1}{1 - x - y - xy}
\, .
$$ 
\end{example}

In both of these examples, $F$ is a rational function.  In one
variable, obtaining asymptotics for the coefficients of a rational
function is a quick exercise (see~(\ref{eq:partial fractions}) in
Section~\ref{ss:univariate} below) but in more than one variable this
is far from true.  For example, whereas the binomial coefficients are
exact expressions easily approximated via Stirling's formula, the
Delannoy numbers do not have such a simple expansion (the expression
$a_{rs} = \sum_{i} \binom{r}{i}\binom{s}{i}2^i$ is readily derived
from the generating function, but no simpler one is forthcoming).

In~\cite[Example~6.3.8]{Stan1999}, it is shown how to get more
information about the \Em{central} Delannoy numbers $a_{nn}$ by
finding the generating function $\sum_n a_{nn} z^n$.  Uniform
estimation of the general term $a_{rs}$ is not possible by this
diagonal method (see Section~\ref{ss:diagonal}), though the problem
succumbs easily to multivariate methods (see~(\ref{eq:delannoy})
below). The class of multivariate generating functions addressed in
this survey is larger than the rational functions; the exact
hypothesis is that the function be meromorphic in a certain domain
which is spelled out in the remark following
Theorem~\ref{th:universal}.

The first paper to concentrate on extracting asymptotics from
multivariate generating functions was~\cite{Bend1973},
already published at the time of Bender's survey, but the seminal
paper is~\cite{BeRi1983}.  In this paper, Bender and Richmond assume that
$F$ has a singularity of the form $A / (z_d - g (\xx))^q$ near the
graph of a smooth function $g$, for some real exponent $q$, where
$\xx$ denotes $(z_1 , \ldots , z_{d-1})$. They show, under appropriate
further hypotheses on $F$, that the probability measure $\mu_n$ one
obtains by renormalizing $\{ a_\rr : r_d = n \}$ to sum to~1 converges
to a multivariate normal when appropriately rescaled.  Their method,
which we call the \Em{GF-sequence method}, is to break the
$d$-dimensional array $\{ a_\rr \}$ into a sequence of
$(d-1)$-dimensional slices and consider the sequence of
$(d-1)$-variate generating functions 
$$
f_n (\xx) = \sum_{\rr : r_d = n} a_\rr \zz^\rr \, .
$$
They show that, asymptotically as $n \to \infty$,
\begin{equation} \label{eq:quasi-power}
f_n (\xx) \sim C_n g(\xx) h(\xx)^n
\end{equation}
uniformly over $\xx$ in a certain ball, 
and that sequences of generating functions obeying~(\ref{eq:quasi-power}) 
satisfy a central limit theorem and a local central limit theorem.

We will review these results in more detail in
Section~\ref{ss:compare}, but one crucial feature is that they always
produce Gaussian (central limit) behavior.  The applicability of the
entire method is therefore limited to the single, though important,
case where the coefficients $a_\rr$ are nonnegative and possess a
Gaussian limit. The work of~\cite{BeRi1983} has been greatly expanded upon,
but always in a similar framework.  For example, it has been extended
to matrix recursions~\cite{BRW1983} and the applicability has been
extended from algebraic to algebraico-logarithmic singularities of the
form $F \sim (z_d - g(\xx))^q \log^\alpha (1/(z_d -
g(\xx)))$~\cite{GaRi1992}. The difficult step is always deducing
asymptotics from the hypotheses~\eqref{eq:quasi-power}.  Thus
some papers in this stream refer to such an assumption in their
titles, and the term ``quasi-power''
has been coined for such a sequence $\{ f_n \}$.

The theory has also been pushed forward via its use in applications.
The forthcoming textbook by Flajolet and
Sedgewick~\cite{FlSe} devotes a
chapter of nearly 100 pages to multivariate asymptotics, in which many
of the basic results on quasi-powers are reviewed and extended.  The
reader is referred to this substantial body of work for an extensive
collection of applications that can be handled by reducing
multivariate problems to univariate contour integrals. The limit
theorems in~\cite[Chapter~IX]{FlSe} (outside of
the graph enumeration example in Section~11) are all Gaussian.  The
large deviation results of Hwang~\cite{Hwan1996, Hwan1998a} are not restricted to the Gaussian case,
but give asymptotics on a cruder scale.

\subsection{New multivariate methods}

Odlyzko's survey of asymptotic enumeration
methods~\cite{Odly1995}, which is meant to be somewhat
encyclopedic, devotes fewer than~6 of its 160 pages to multivariate
asymptotics.  Odlyzko describes why he believes multivariate
coefficient estimation to be difficult. First, the singularities are
no longer isolated, but form $(d-1)$-dimensional hypersurfaces.  Thus,
he points out, ``Even rational multivariate functions are not easy to
deal with.'' Secondly, the multivariate analogue of the
one-dimensional residue theorem is the considerably more difficult
theory of Leray~\cite{Lera1959}. This theory was later fleshed out by
Aizenberg and Yuzhakov, who spend a few
pages~\cite[Section~23]{AiYu1983} on generating functions
and combinatorial sums.  Further progress in using multivariate
residues to evaluate coefficients of generating functions was made by
Bertozzi and McKenna~\cite{BeMc1993}, though at the
time of Odlyzko's survey none of the papers based on multivariate
residues such as~\cite{Lich1991,BeMc1993} had resulted
in any kind of systematic application of these methods to enumeration.

The topic of the present review article is a recent vein of research
begun in~\cite{PeWi2002} and continued in~\cite{PeWi2004,BaPe,Llad2003} and
in several manuscripts in progress.  The idea, seen already to some
degree in~\cite{BeMc1993}, is to use complex methods
that are genuinely multivariate to evaluate coefficients of
multivariate generating functions via the multivariate Cauchy formula.
By avoiding symmetry-breaking decompositions such as $F(\zz) = \sum f_n
(z_1 , \ldots , z_{d-1}) z_d^n$, one hopes that the methods will be
more universally applicable and the formulae more canonical.  In
particular, the results of Bender {\it et al.} and the results of
Bertozzi and McKenna are seen to be two instances of a more general
result that estimates the Cauchy integral via topological reductions of
the cycle of integration.  These topological reductions, while not
fully automatic, are algorithmically decidable in large classes of
cases and are the subject of Section~\ref{ss:results}. An ultimate
goal, stated in~\cite{PeWi2002,PeWi2004}, is to develop computer software to
automate all of the computation.

Aside from providing a summary and explication of this line of
research, the present survey is meant to serve several other purposes.
First, results from~\cite{PeWi2002,PeWi2004,BaPe} are presented in streamlined
forms, stated so as to avoid the scaffolding one needs to prove them.
This is to make the results more comprehensible.  Secondly, by
focusing on combinatorial applications, we hope to create a sort of
user's manual: one that contains worked examples akin to those a
potential user will have in mind.  Many of the applications are to
abstract combinatorial structures but we also include direct
applications to computational biology and to formal languages and
automata theory.  Finally, we present a number of results that
``pre-process'' the basic, general theorems, providing useful
computational reductions of the hypotheses or conclusions in specific
cases of interest.  We now give the notation for generating functions
and their asymptotics that will be used throughout, and then briefly
describe a prototypical asymptotic theorem.

Throughout this survey, we let
\begin{equation} \label{eq:F}
F(\zz) = \frac{\numer(\zz)}{\denom(\zz)} 
   = \sum_{\rr} a_\rr \zz^\rr
\end{equation}
be a generating function in $d$ variables, where $\numer$ and $\denom$ 
are analytic and $\denom(\zero) \neq 0$.  Recall that in the 
bivariate case ($d=2$) we write
$$F(x,y) = \frac{\numer(x,y)}{\denom(x,y)} 
   = \sum_{r,s = 0}^\infty a_{rs} x^r y^s \, .$$
The representation of $F$ as a quotient of analytic functions is
required only to hold in a certain domain, described in the remark
following Theorem~\ref{th:universal}, though in the majority of the
examples $F$ is meromorphic on all of $\C^d$. We will assume
throughout that $\denom$ vanishes somewhere, since the methods in this
paper do not give nontrivial results for entire functions.

We are concerned with asymptotics when $|\rr| \to \infty$ with $\rhat
:= \rr / |\rr|$ remaining in some specified set, bounded away from the
coordinate planes.  Thus for example, when $d=2$, the ratio $s/r$ will
remain in a compact subset of $(0,\infty)$. It is possible via our
methods to address the other case, where $r = o(s)$ or $s = o(r)$,
(see, for example,~\cite{Llad2003}), but our main purpose in
this paper is to give examples that require, among all the methods and
results cited above, only those from~\cite{PeWi2002} together with the
simplest methods from~\cite{PeWi2004,BaPe}.

To illustrate what sort of basic results we quote from~\cite{PeWi2002,PeWi2004}, 
we quote the following combination of Corollary~\ref{cor:procedure} 
and Theorem~\ref{th:smooth nice} from Section~\ref{ss:results}.
\begin{thm} \label{th:prototype}
Let $F$ be as in~(\ref{eq:F}) and suppose $a_\rr \geq 0$.
\begin{enumerate}
\romenumi
\item For each $\rr$ in the positive orthant, there is a unique
$\zz (\rr)$ in the positive orthant satisfying the equations 
$r_d z_j \partial \denom / \partial z_j = r_j z_d \partial \denom / 
\partial z_d \medspace (1 \leq j \leq d-1)$ from~(\ref{eq:zero-dim}) below, 
and lying on the boundary of the domain of convergence of $F$; 
the quantity $\zz(\rr)$ depends on $\rr$ only through the direction 
$\rr / |\rr|$.
\item With $\zz (\rr)$ defined in this way, if $\numer (\zz(\rr)) \neq 0$,
$$a_\rr \sim (2 \pi )^{-(d-1)/2} \hess^{-1/2} \frac{\numer (\zz (\rr))}
   {-z_d \partial \denom / \partial z_d (\zz (\rr))}
   (r_d)^{-(d-1)/2} \zz (\rr)^{-\rr}$$
uniformly over compact cones of $\rr$ for which $\zz (\rr)$ is a
smooth point of $\{ \denom = 0 \}$ uniquely solving~(\ref{eq:zero-dim}) 
on the boundary of the domain of convergence of $F$, and for which 
$\hess$ is nonzero, where $\hess$ is the determinant of the
Hessian matrix of the function parametrizing the hypersurface
$\{ \denom = 0 \}$ in logarithmic coordinates.
\end{enumerate}
\end{thm}
Going back to Example~\ref{eg:pascal}, let us see what this says
about binomial coefficients.  It will turn out 
(see Section~\ref{ss:binomial}) that $\zz (\rr) =
(\frac{r}{r+s} , \frac{s}{r+s})$.  Evaluating the Hessian and the 
partial derivatives then leads to equation~(\ref{eq:binomial}) below:
$$
a_{rs} \sim \left ( \frac{r+s}{r} \right )^r 
   \left ( \frac{r+s}{s} \right )^s \sqrt{\frac{r+s}{2 \pi r s}}
$$
as $r,s \to \infty$ at comparable rates.
This agrees, of course, with Stirling's formula.  When we try this
with the Delannoy numbers from Example~\ref{eg:delannoy}, we find that
$$
\zz (r,s) = \left ( \frac{\sqrt{r^2+s^2} - s}{r} , 
   \frac{\sqrt{r^2+s^2} - r}{s} \right )
$$
and that  
\begin{equation} \label{eq:delannoy}
a_{rs} \sim \left (\frac{\sqrt{r^2 + s^2} - s}{r} \right)^{-r}
   \left ( \frac{\sqrt{r^2 + s^2} - r}{s} \right )^{-s}
   \sqrt{\frac{1}{2 \pi}} \sqrt{\frac{r s}{(r+s-\sqrt{r^2 + s^2})^2
   \sqrt{r^2 + s^2}}}
\end{equation}
uniformly as $r,s \to \infty$ with $r/s$ and $s/r$ remaining bounded.
Setting $r = s = n$ gives the estimate 
$$
a_{nn} \sim (\sqrt{2} - 1)^{-2n} \sqrt{\frac{1}{2 \pi n}}
\frac{2^{-1/4}} {2 - \sqrt{2}}
$$
for the central Delannoy numbers.  This estimate can also be computed, with more effort,
by using the methods of Section~\ref{ss:transfer} applied to the diagonal generating function, 
which can itself be derived using the  \Em{diagonal method} as in~\cite[Example 6.3.8]{Stan1999} (we discuss the diagonal method in Section~\ref{ss:diagonal}).

In the sections to follow, we work through a number of applications of
this and other, newer theorems to problems in combinatorial
enumeration.  While doing so, we develop companion results that
simplify the computations in certain classes of examples. A
prototypical companion result is the following result for implicitly
defined generating functions, which arise commonly in recursions on
trees.
\begin{quote} {\em
If $f(z) = z \phi (f(z))$ then the $n^{th}$ coefficient of $f$ is
given asymptotically by
\begin{equation} \label{eq:thm 2}
n^{-3/2} \frac{y_0 \phi' (y_0)^n}{\sqrt{2 \pi \phi'' (y_0) / \phi
(y_0)}}
\end{equation}
where $y_0$ is the unique positive solution to $y \phi' (y) = \phi
(y)$.} 
\end{quote}
In~\cite{FlSe} the computation is reduced via
Lagrange inversion to the determination of the $y^{n-1}$ coefficient of
$\phi^n (y)$, which they then carry out via complex integration methods.
By viewing this coefficient instead as the $(n,n)$ coefficient of $$F(x,y) :=
\frac{y}{1 - x \phi (y)} \, ,$$ one sees that the complex integration
step is already done, and that~(\ref{eq:thm 2}) is immediate upon
identifying that $\zz (n,n) = (1 / \phi (y_0) , y_0)$; see
Proposition~\ref{pr:lagrange}, where the hypotheses for~(\ref{eq:thm
2}) are stated more completely.

The organization of the remainder of this paper is as follows. The
next section gives the promised quick primer in univariate methods.
Section~\ref{ss:results} is the most theoretical section, outlining
results from various sources to give a brief but nearly complete
explanation as to how one goes from~(\ref{eq:F}), via a multivariable
Cauchy formula, to asymptotic formulae for $a_\rr$. We then quote the
precise theorems we will use from~\cite{PeWi2002} and~\cite{PeWi2004}.  While we
attempt here to provide short and accessible statements of results,
readers interested only in knowing enough to handle a particular
application may wish to skip this section, find the closest match
among the examples in Sections~\ref{ss:details}~--~\ref{sec:kernel},
and then refer back to the necessary parts of
Section~\ref{ss:results}. Section~\ref{ss:details} works in detail
through a number of diverse examples.  The subsequent three sections
discuss collections of examples arising from three combinatorial
methods: transfer matrices, Lagrange inversion and the kernel method, 
respectively; we present applications to enumerating various kinds of
words, paths, trees and graphs. In Section~\ref{ss:compare} we discuss
open questions and extensions of the material presented here, and
compare our results with those of other authors.

\setcounter{equation}{0}
\section{A brief review of univariate methods} 
\label{ss:univariate} 

Since our main subject is the asymptotics of coefficients of
multivariate generating functions, we will give only a quick 
overview of the univariate case touching on three widely used methods, namely saddle 
point analysis, Darboux' circle method and branch point contours.  
In each case we will give a brief summary and an application or two,
as well as some pointers to the literature.  For readers who 
wish to understand univariate asymptotic methods in greater detail, 
we recommend beginning with Sections~10--12 of Odlyzko's survey 
paper~\cite{Odly1995}.  

We begin by disposing of a trivial case.  If $f(z)$ is a rational 
function, written as $P(z) / Q(z)$ for polynomials $P$ and $Q$, 
then a partial fraction decomposition exists:
\begin{equation} \label{eq:partial fractions}
f(z) = \sum_{\alpha , k} \frac{P_{\alpha , k} (z)}{(1 - z / \alpha)^k}
\end{equation}
where $\alpha$ ranges over roots of $Q$ and the integer $k$
is between~1 and the multiplicity of the root $\alpha$.  
For $j \geq i$, the $z^j$ coefficient of 
$\displaystyle{\frac{z^i}{(1 - z/\alpha)^k}}$ is 
$\displaystyle{\alpha^{i-j} \binom{j-i+k}{k}}$, 
which leads to a representation 
$$a_n = \sum_\alpha p_\alpha (n) \alpha^{-n}$$
for easily computed polynomials $p_\alpha$.

We turn now to analyses via contour integrals.
To see what underlies all of these methods, we begin with a 
preliminary estimate.  Suppose that $f$ has radius of convergence 
$R$.  Taking the contour of integration in~(\ref{eq:cauchy0}) to be 
a circle of radius $R - \ee$ gives $a_n = O(R - \ee)^{-n}$ 
for any $0 < \ee < R$.  If $f$ is continuous on the closed disk of radius $R$ 
then $a_n = O(R^{-n})$.  All the methods we discuss
refine these basic estimates by pushing the contour out far
enough so that the resulting upper bound becomes asymptotically sharp.

\subsection{Saddle point methods} \label{ss:saddle}

We begin with a worked example, then comment on the method in general.

\begin{example}[set partitions] \label{eg:essential}
Let $a_n$ be the number of (unordered) partitions of $\{ 1 , \ldots ,
n \}$ into ordered sets.  This has exponential generating 
function~\cite[page~194]{Rior1968}
$$
f(z) := \sum_{n=0}^\infty \frac{a_n}{n!} z^n = e^{z/(1-z)} \, .
$$
To evaluate the Cauchy integral, we attempt to move the contour so it
passes through a point where the integrand is not rapidly oscillating.
Specifically, we let $I_n = \log (z^{-n-1} f(z))$ denote the logarithm
of the integrand and we find where $I_n'$ vanishes.  It is not hard to
see that 
$$
I_n' = \frac{-n-1}{z} + \frac{1}{(1-z)^2}
$$
vanishes at a value $1 - \beta_n$ where $\beta_n = n^{-1/2} +
O(n^{-1})$. Expand the contour to a circle passing through $1 -
\beta_n$, or equivalently and more cleanly, a piece of the line $1 -
\beta_n + i x$ near $x=0$. Replacing the integrand $\exp(I_n(z)) \,
dz$ by its two-term Taylor approximation, one obtains
\begin{align*}
\frac{a_n}{n!} & =  \frac{1}{2 \pi i} \int_{-\infty}^\infty
   \exp(I_n(1 - \beta_n + it)) (i \, dt) \\
& \sim  \frac{\exp(I_n(1-\beta_n))}{2 \pi} \int_{-\infty}^\infty 
   \exp \left ( \frac{1}{2} I_n'' (1 - \beta_n) (it)^2 \right ) \, dt \\
& =  \frac{\exp(I_n(1-\beta_n))}{2 \pi} 
   \sqrt{\frac{1}{2 \pi I_n''(1 - \beta_n)}} \, .
\end{align*}
The approximation may be justified by routine estimates.  Plugging 
in the value of $1 - \beta_n$ into $I$ and $I_n''$ shows the estimate 
to be asymptotically 
$$
(1 + o(1)) \sqrt{\frac{1}{4 \pi e n^{3/2}}} \exp (2 \sqrt{n}) \, .
$$
\end{example}

For this method to work either $f$ must be entire or the saddle 
point (where $I_n'$ vanishes) must be in the interior of the domain of 
convergence of $f$.  These conditions are frequently satisfied and
the method is widely applicable.  We do not know when saddle point
approximation was first used in the context of generating functions.
Hayman's seminal paper from 1956~\cite{Haym1956} 
defines a broad class of functions, called \Em{admissible functions}
for which the saddle point method can be shown to work and
the Gaussian approximation mechanized.  A variation of the generating
function $\exp(z/(1-z))$ is analyzed in Hayman's Theorem~XIII.

\subsection{Circle methods} \label{ss:circle}

Suppose $f$ has a positive, finite radius of convergence, $R$. If
there is no saddle point of $z^{-n-1} f(z)$ in the open disk of radius
$R$, one might try to push the contour of integration near or onto the
circle of radius $R$.  Darboux' method is essentially this, with the
refinement that if $f$ extends to the circle of radius $R$ and is $k$
times continuously differentiable there, then 
\begin{equation} \label{eq:darboux method}
\int z^{-n} f(z) \, dz = O \left ( n^{-k} R^{-n} \right ) \,.
\end{equation}
This follows from integration by parts.  

The following very old theorem may be found, among other places,
in~\cite[Theorem~11.10b: ``Theorem of Darboux'']{Henr1991}.
\begin{thm}[Darboux] \label{th:darboux}
Let $f(z) = \sum_n a_n z^n = (r - z)^\alpha L(z)$ where $r > 0$,
$\alpha$ is not a positive integer, and $L$ is analytic in a disk 
of radius greater than $r$.  Then
$$a_n \sim r^{\alpha-n} n^{-\alpha - 1} \frac{L(r)}{\Gamma (-\alpha)} \, .$$
\end{thm}
\begin{example}[2-regular graphs] \label{eg:darboux}
Let $b_n$ be the number of 2-regular graphs.  The
exponential generating function~\cite[equation~(3.9.1)]{Wilf1994}
is given by
$$
f(z) := \sum_n \frac{b_n}{n!} z^n = \frac{e^{-z/2 - z^2/4}}{\sqrt{1-z}} \, .
$$
Letting $\alpha = -1/2, r=1$ and $L(z) = e^{-z/2 - z^2/4}$ gives
$$
a_n := \frac{b_n}{n!} \sim \frac{e^{-3/4}}{\sqrt{\pi n}} \, .
$$
\end{example}

To prove the above version of Darboux' Theorem, write $L(z)$ as
$m$ terms of a Taylor series about $r$ plus a remainder.  This
expresses $f(z)$ as 
$$
\sum_{k=0}^{m-1} c_k (r-z)^{\alpha + k} + O \left ( (r-z)^{\alpha + m}
   \right )  \, .
$$
The known asymptotics for coefficients of $(r - z)^\nu$ together
with the estimate~(\ref{eq:darboux method}) is good enough to yield
$m-1$ terms in the conclusion provided one takes $m \geq 1 + 
(\Real \{ - \alpha \} )^+$, where $x^+$ denotes the maximum of $x$ and~0.
$\Cox$

There are a great many variations on this, depending on how badly
$f$ behaves near the circle of radius $r$.  One of the most classical
fruits of the circle method is Hardy and Ramanujan's 
estimate~\cite{HaRa1917} of the partition numbers. 
\begin{example} 
Let $p_n$ denote the number of partitions of $n$, that is,
representations of $n$ as a sum of positive integers with the summands
written in descending order.  The number of identities involving these
numbers and their relatives is staggering (see,
e.g.,~\cite{AnEr2004} for an elementary survey). Euler
observed that the generating function may be written as
$$
\lambda (z) := \sum_{n=0}^\infty p_n z^n = \prod_{k=0}^\infty
   \frac{1}{1-z^k} \, .
$$
 From this, Hardy and Ramanujan obtained
$$
p_n \sim \frac{\exp (\pi \sqrt{2n/3})}{4 \sqrt{3} n}
$$
as $n \to \infty$.  
\end{example}

\subsection{Transfer theorems} \label{ss:transfer}

A closer look at the proof of Darboux' Theorem shows that one
can do better.  Analyticity of $f/(r-z)^\alpha$ beyond the
disk of radius $r$ is used only to provide a series for $f$
in decreasing powers of $(r-z)$.  On the other hand, taking the
contour to be a circle of radius $r$ loses one power of $n$ in
the estimate, though this doesn't matter when the expansion of $f$
has sufficiently many terms.  

\begin{figure}[ht] 
\hspace{0.2in} \includegraphics[scale=0.33]{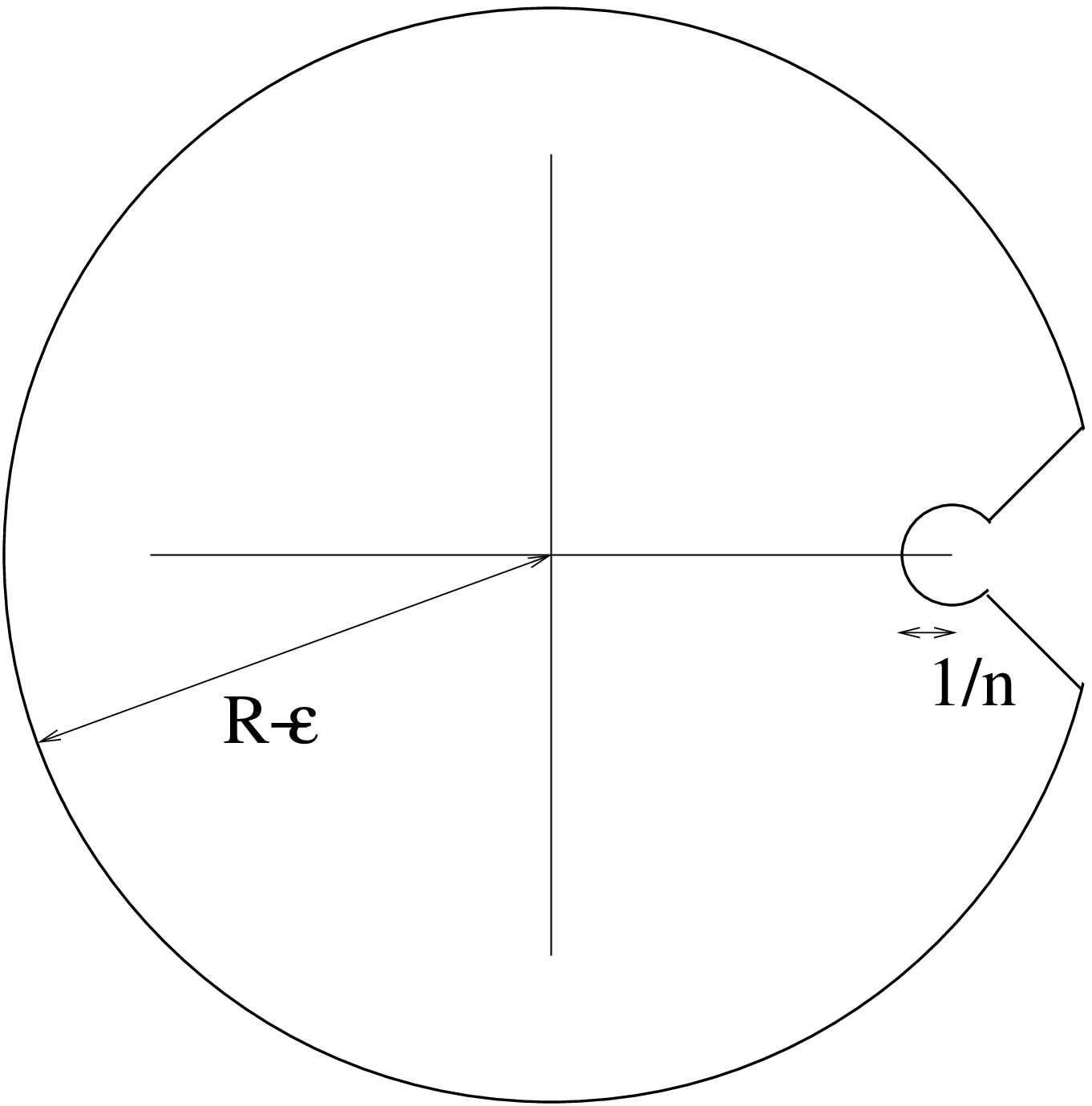}
\vspace{0.2in} \hspace{1.2in} \includegraphics[scale=0.78]{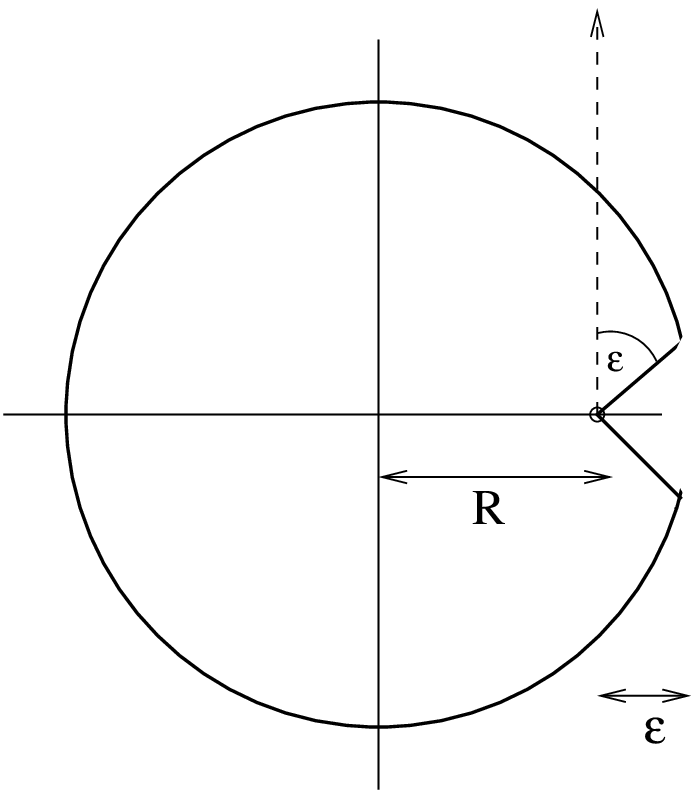}
\caption{contour for an algebraic singularity and the corresponding
region of analyticity}
\label{fig:contour} 
\end{figure}

In light of these observations, one sees that choosing a custom-built
contour can simultaneously weaken the analyticity assumption 
while strengthening the estimate.  The choice of contour will
depend on the nature of the singularity of $f$ on the boundary 
of its domain of convergence, but in many instances a good choice
is known.  Both a non-integral power of $(r-z)$ and a logarithm
or iterated logarithm of $1/(r-z)$ are branch singularities and
require similar contours.  In the simplest case, the contour,
pictured in Figure~\ref{fig:contour},
consists of an arc of a circle of radius $1/n$ around $r$, an
arc of a circle of radius $r - \ee$ centered at the origin, and
two line segments connecting corresponding endpoints of the two
circular arcs.

Among the many refinements of Darboux' Theorem via such contours
are the \Em{transfer theorems} of Flajolet and Odlyzko.  Their 
paper~\cite{FlOd1990} is masterfully written and
worth taking a couple of pages here to summarize.  They begin 
by isolating the gist of~(\ref{eq:darboux method}), namely that 
$a_n = O(n^{-\lfloor \alpha \rfloor})$ if $f(z) = O(1 - z)^\alpha$.  
Already this allows for $f$ not to be in any nice class of functions, 
as long as it is bounded by a class we understand.  Next, they use 
the above contour to improve this to a sharp estimate,
\begin{equation} \label{eq:FO1}
f(z) = O(1-z)^\alpha \;\; \Rightarrow \;\; a_n = O(n^{-\alpha - 1}) \,
.
\end{equation}
Let $\alg$ denote the class of functions that are a product of 
a power of $r-z$, a power of $\log(1/(r-z))$ and a power of 
$\log \log (1/(r-z))$.  Flajolet and Odlyzko estimate $a_n$
for any function in $\alg$.  
Their last crucial observation is that the implication
in~(\ref{eq:FO1}) gives an explicit constant for the bound on $a_n$
based on the constant in the bound $f(z) = O(1-z)^\alpha$,
which is sufficiently constructive to give
$$
f(z) = o \left ( (r-z)^\alpha \right ) \;\; \Rightarrow 
   \;\; a_n = o(n^{-\alpha - 1} ) \, .
$$
The surprisingly strong consequence of this is that any $f$ 
with an asymptotic expansion by functions $\{ g_k \}$ in $\alg$ 
has coefficients asymptotically given by summing the asymptotics 
of the coefficients of the functions $g_k$.

The main result of~\cite{FlOd1990} begins by 
stating what analyticity assumption is necessary in
order to use the contour in Figure~\ref{fig:contour}.
Given a positive real $R$ and $\ee \in (0,\pi/2)$, the
so-called \Em{Camembert-shaped region}, 
$$
\{ z : |z| < R + \ee , z \neq R, |\arg (z-R)| \geq \pi / 2 - \ee \} \,
,
$$  
shown on the right of Figure~\ref{fig:contour}, 
is defined so that it includes the contour in Figure~\ref{fig:contour}.
\begin{thm} \label{th:FO}
Let $f(z) = \sum a_n z^n$ be singular at $R$ but analytic in a
Camembert-shaped region.  For $g (z) = \sum b_n z^n$ in the class
$\alg$, asymptotic relations between $f$ and $g$ as $z \to R$
imply estimates on the coefficients as follows.
\begin{enumerate}
\romenumi
\item 
$$f(z) = O(g(z)) \Rightarrow a_n = O(b_n) \, ;$$
\item 
$$f(z) = o(g(z)) \Rightarrow a_n = o(b_n) \, ;$$
\item 
$$f(z) \sim g(z) \Rightarrow a_n \sim b_n \, .$$
\end{enumerate}
In particular, when $f(z) \sim C (r-z)^\alpha$, this result subsumes 
Theorem~\ref{th:darboux}.
\end{thm}

We have indicated the arguments deriving the second and third part of
the theorem from the first.  We now sketch, in the simplest case where
$g(z) = (1-z)^\alpha$, the upper estimate on $a_n$, which appears as
Corollary~2 of~\cite{FlOd1990}. The integral over the large
circular arc is small because $z^{-n-1}$ is exponentially small there.
The integral over the small circular arc is $O(n^{-\alpha - 1})$
because the integrand is $O(n^{-\alpha})$ and the arc has length
$O(n^{-1})$.  On the two radial line segments, the modulus of
$z^{-n-1}$ is small on the line segments except in a neighborhood of
size $O(1/n)$ of~1.  Here, changing variables to $(z-1)/n$, one has an
integral of the form $\int A(z) \exp (n \phi (z)) \, dz$ from which
one obtains (e.g., via the well known Watson-Doetsch
lemma~\cite[Theorem~11.5]{Henr1991}) a value \begin{equation}
\label{eq:darboux} A(1) \, \frac{n^{-\alpha - 1}}{\Gamma (-\alpha)} \,
. \end{equation}
$\Cox$

\begin{example}[Catalan numbers] \label{eg:FO}
Let $\displaystyle{a_n := \frac{1}{n+1} \binom{2n}{n}}$ be the 
$n^{th}$ Catalan number.  A great many naturally occurring
combinatorial classes are counted by this sequence: there is a
list of 66 of these in~\cite[Problem~6.19]{Stan1999}.  The
generating function for the Catalan numbers is
$$f(z) := \sum_{n=0}^\infty a_n z^n = \frac{1-\sqrt{1-4z}}{2z}
   = \frac{1 - 2 \sqrt{\frac{1}{4} - z}}{2z} \, .$$
There is an algebraic singularity at $r = 1/4$, near which the
asymptotic expansion for $f$ begins
$$f(z) = 2 - 4 \sqrt{\frac{1}{4}-z} + 8 (\frac{1}{4}-z) - 16 
   (\frac{1}{4} - z)^{3/2} + O(\frac{1}{4} - z)^2 \, .$$
Note that $f/\sqrt{1/4-z}$ is not analytic in any disk of radius 
$1/4  +\ee$, since both integral and half-integral powers appear in $f$,
but $f$ is analytic in a Camembert-shaped region.  Theorem~\ref{th:FO} 
thus gives (note that the integral powers of $(1-z)$ do not contribute):
\begin{eqnarray*}
a_n & \sim & \left ( \frac{1}{4} \right )^{1/2 - n} n^{-3/2}
   \frac{-4}{\Gamma (-1/2)} + \left ( \frac{1}{4} \right )^{3/2 - n} n^{-5/2}
   \frac{-16}{\Gamma (-3/2)} + O(n^{-7/2}) \\
& = & 4^n n^{-3/2} \frac{(-4) (\frac{1}{4})^{1/2}}{\Gamma (-1/2)} + 
   4^n n^{-5/2} \frac{(-16) (\frac{1}{4})^{3/2}}{\Gamma (-3/2)} + 
   O(n^{-7/2}) \\
& = & 4^n \left ( \frac{n^{-3/2}}{\sqrt{\pi}} 
   - n^{-5/2} \frac{3}{ 2 \sqrt{\pi}} + O(n^{-7/2}) \right ) \, .
\end{eqnarray*}
\end{example}

We have given only a brief view of what is known in the univariate case.  
We close this section with several more pointers to the literature.

One review of univariate asymptotics, that predates the work of 
Flajolet and Odlyzko but is still quite useful, 
is~\cite[part~II]{Bend1974}.  The 1995 survey
article by Odlyzko~\cite{Odly1995} to which we have already
referred is somewhat more extensive.  
As mentioned before, 
contours such as the one in Figure~\ref{fig:contour} occur
primarily in the univariate literature, but one recent
extension to the multivariate setting, via a product of
these contours, occurs in~\cite{GaWo2000}: their
Lemma~3 gives a Darboux-type estimate for the 
$(n , k_1 , \ldots , k_j)$-coefficient of a generating function 
asymptotic to $(1-z_1)^{-\alpha} \prod_{i=1}^j (1 - z_i)^{- \beta_i}$
uniform as long as $k_i = O(n)$ for all $i$.
Finally, the book~\cite{FlSe} is the most 
up to date, though it is only available in electronic preprint form 
at this time.  

\setcounter{equation}{0}
\section{New multivariate results} 
\label{ss:results} 

In~\cite{PeWi2002,PeWi2004,BaPe} the authors derived asymptotic formulae for $a_\rr$, as
$|\rr| \to \infty$, that are uniform as $\rr / |\rr|$ varies over some
compact set.  It is useful to separate $\rr$ into the scale parameter
$|\rr|$, which is a positive real number, and a direction parameter
$\rbar$, which is an element of real projective space.  Although $\rr$ is
always an element of the positive orthant of $\R^{d}$, it will sometimes
make sense to consider it as an element of $\C^{d}$.  When convenient, 
we identify $\rr$ with its class $\rbar$ in projective space or its
projection $\rhat$ in the real $(d-1)$-simplex $\simplex$, the set 
$\{\xx \in (\R^+)^d : |\xx| = 1 \}$ (where  $|\xx|:= \sum_{j=1}^d x_j$).  Thus $\rr$ may
appear anywhere in the following diagram, where $\mathcal{O}$ and
$\orthant$ denote the positive orthants of $\R^{d}$ and $\RP^{d-1}$
respectively.

\setlength{\unitlength}{1.5pt}
\begin{picture}(170,80)(0,-5)
\put(30,60){$\mathcal{O}$}
\put(80,60){$\R^{d}$}
\put(130,60){$\C^{d}$}
\put(30,30){$\orthant$}
\put(80,30){$\RP^{d-1}$}
\put(130,30){$\CP^{d-1}$}
\put(30,0){$\simplex$}
\put(49,62){\vector(1,0){27}}
\put(99,62){\vector(1,0){27}}
\put(49,32){\vector(1,0){27}}
\put(99,32){\vector(1,0){27}}
\put(34,56){\vector(0,-1){20}}
\put(84,56){\vector(0,-1){20}}
\put(134,56){\vector(0,-1){20}}
\put(34,26){\vector(0,-1){20}}
\put(24,13){$\widehat{\;\;}$}
\put(24,48){$\overline{\;\;}$}
\put(74,48){$\overline{\;\;}$}
\put(124,48){$\overline{\;\;}$}
\end{picture}

Quantities that depend on $\rr$ only through its direction will be 
denoted as functions of $\rbar$.  The results we will quote in this 
section may be informally summarized as follows.  Recall that
$F = \numer / \denom$.

\begin{enumerate}
\romenumi
\item Asymptotics in the direction $\rbar$ are determined by the 
geometry of the pole variety $\sing = \{ \zz : H (\zz) = 0 \}$ near a
finite set, $\crit_{\rbar}$, of \Em{critical points}
[Definition~\ref{def:crit}].
\item For the purposes of asymptotic computation, one may reduce this
set of critical points further to a set $\contrib_{\rbar} \subseteq
\crit_{\rbar}$ of \Em{contributing critical points}, usually a single
point [Formula~(\ref{eq:decomp}) and Definition~\ref{def:contrib}].
\item One may determine $\crit_{\rbar}$ and $\contrib_{\rbar}$ by 
a combination of algebraic and geometric criteria  
[Proposition~\ref{pr:crit} and Theorem~\ref{th:contrib} respectively]. 
\item Critical points may be of three types: smooth, multiple or bad
[Definition~\ref{def:classify}].
\item Corresponding to each smooth or multiple critical point, $\zz$, 
is an asymptotic expansion for $a_\rr$ which is computable in terms
of the derivatives of $\numer$ and $\denom$ at $\zz$
[Sections~\ref{ss:smooth} and~\ref{ss:multiple} respectively].
\end{enumerate}

The culmination of the above is the following meta-formula:
\begin{equation} \label{eq:meta}
a_\rr \sim \sum_{\zz \in \contrib_{\rbar}} \mbox{ \bf formula} (\zz)
\end{equation}
where $\mbox{\bf formula} (\zz)$ is one function of the
local geometry for smooth points and a different function 
for multiple points.  Specific instances of~(\ref{eq:meta}) are 
given in equations~(\ref{eq:smooth asym})~--~(\ref{eq:m<d});
the simplest case is
$$a_\rr \sim C |\rr|^{\frac{1-d}{2}} \zz^{-\rr}$$
where $C$ and $\zz$ are functions of $\rbar$.  
No general expression for \textbf{formula} (\zz) is yet known when 
$\zz$ is a bad point, hence the name ``bad point''.

Fundamental to all the derivations is the Cauchy integral representation 
\begin{equation} \label{eq:cauchy}
a_\rr = \frac{1}{(2 \pi i)^{d}} \int_T  \zz^{- \rr - \one} F(\zz) \, d\zz
\end{equation}
where $T$ is the product of sufficiently small circles around the
origin in each of the coordinates, $\one$ is the $d$-vector of all
ones, and $d\zz$ is the holomorphic volume form $dz_1 \wedge \cdots
\wedge dz_{d}$. The gist of the analyses in \cite{PeWi2002, PeWi2004, BaPe} is the
computation of this integral.  The key step is to replace the cycle
$T$ by a product $T_1 \times T_2$, where the inner integral over $T_1$
is a multivariate residue and the outer integral over $T_2$ is a
saddle point integral.  A complete asymptotic series may then be read
off in a straightforward manner by well known methods.  The first
subsection below explains what critical points and contributing
critical points are and gives the ``big picture'', namely the
topological context as outlined in~\cite{BaPe}.  The second subsection
gives the specific definitions from~\cite{PeWi2002} that we need to 
find the contributing critical points (note: the question of 
sorting these into smooth, multiple or bad critical points is 
addressed at the beginning of Section~\ref{ss:details}).
Sections~\ref{ss:smooth} and~\ref{ss:multiple} quote
results from~\cite{PeWi2002} and~\cite{PeWi2004} that give asymptotics for
smooth points if they are respectively smooth or multiple.
Section~\ref{ss:limit theorems} restates some of these asymptotics in
terms of probability limit theorems.

\subsection{Topological representation}

As $\rr \to \infty$, the integrand in~(\ref{eq:cauchy}) becomes
large. It is natural to attempt a saddle point analysis. That is,
we try to deform the contour, $T$, so as to minimize the maximum
modulus of the integrand.  If $\rhat$ remains fixed, 
then the modulus of the integrand is well approximated
by the exponential term $\displaystyle{\exp \left ( 
- |\rr| (\rhat \cdot \log |\zz|) \right )}$,
where $\log |\zz|$ is shorthand for the vector $(\log |z_1| , 
\ldots , \log |z_{d}|)$.  
This suggests that the real function 
\begin{equation} \label{eq:h}
h(\zz) := - \rhat \cdot \log |\zz|
\end{equation}
be thought of as a \Em{height function} in the Morse theoretic sense,
and that we try to deform $T$ so that its maximum height is as low as
possible.  Stratified Morse Theory~\cite{GoMa1988} solves
the problem of accomplishing this optimal deformation.  We now give a 
brief summary of this solution; for details, consult~\cite{BaPe}.

The variety $\sing$ may be given a Whitney stratification.  If $\sing$
is already a manifold (which is the generic case), there is just a
single stratum, but in general there may be any finite number of
strata, even in two dimensions. The set of smooth points of $\sing$
always constitutes the top stratum, with the set of singular points
decomposing into the remaining strata, each stratum being a smooth
manifold of lower complex dimension.  To illustrate some of the 
possibilities, we give two examples; for further details on 
Whitney stratification, see~\cite[Section~I.1.2]{GoMa1988}.  
%
\begin{example}[complete normal intersection] 
Let $\sing$ be the variety where $(1-x)(1-y)(1-z)$ vanishes.
The smooth points are those where exactly one of $x,y$ or $z$ 
is equal to 1.  There are three one-dimensional strata where two but
not all of $1-x, 1-y$ and $1-z$ vanish, and one zero-dimensional stratum 
at $(1,1,1)$.  
\end{example}
%
\begin{example}[isolated singular point]
Let $\sing$ be the variety $\{ 1-x-y-z+4xyz = 0 \}$.  This
has a zero-dimensional stratum at its unique singularity, $(\frac{1}{2}
, \frac{1}{2} , \frac{1}{2})$; everything else constitutes the
two-dimensional stratum. 
\end{example}

In each stratum, there is a finite set of critical points, 
namely points where the gradient of the height function restricted
to the stratum vanishes.  Since the height function depends on
$\rr$ via $\rhat$, the set of critical points may be defined
as a function of $\rhat$ or equivalently of $\rbar$:
\begin{defn}[critical points] \label{def:crit}
If $S$ is a stratum of $\sing$, define the set $\crit_{\rbar} (S)$
of critical points of $S$ for $\denom$ in direction $\rbar$ 
to be the set of $\zz \in S$ for which $\grad h |_S (\zz) = 0$.  
The set $\crit_{\rbar}$ of all critical points of $\sing$
is defined as the union of $\crit_{\rbar} (S)$ as $S$ varies over 
all strata.
\end{defn}

If $S$ has dimension $k$ then the condition $\grad h |_S = 0$
defines an analytic variety of co-dimension $k$.  Generically, then, 
$\crit_{\rbar} (S)$ is zero-dimensional, that is, finite.  It does 
happen sometimes that for a small set of values of $\rbar$ 
the critical set is not finite; for example this happens when $\sing$
is a binomial variety $\{ \zz : \zz^{\bf a} - \zz^{\bf b} = 0 \}$, 
in which case $\grad h |_S$ vanishes everywhere
on $S$ for a particular $\rbar$ and nowhere on $S$ for other values
of $\rbar$.  We address such a case in Section~\ref{ss:integer}.
Aside from that, such a degeneracy does not occur among our examples
and a standing assumption~(\ref{ass:hole} below) will ensure that
it does not cause trouble in any case.

If $F$ is a rational function, then the critical points are
a finite union of zero-dimensional varieties, so elimination theory
(see, e.g.,~\cite[Chapters~1~and~2]{CLO2005})
will provide, in an automatic way, minimal polynomials for the
critical points; see Proposition~\ref{pr:crit} below.  So as to 
refer to them later, we enumerate the critical points as 
$\{ \zz^{(1)} , \ldots , \zz^{(m)} \}$.  

Let $\M$ denote the domain of holomorphy of the integrand 
in~(\ref{eq:cauchy}), that is, 
$$
\M := \C^{d} \setminus \left \{ \zz : 
   \denom \cdot \prod_{j=1}^{d} z_j = 0 \right \} \, .
$$
The topology of $\M$ is determined by neighborhoods of the critical
points of $\sing$ (the complement of $\M$) in a manner we now
describe. Given a real number $c$, let $\M^c$ denote the set of points
of $\M$ with height less than $c$.  Let $\M^+$ denote the $\M^c$ for
some $c$ greater than the greatest height of any critical point and
let $\M^-$ denote $\M^c$ for some $c$ less than the least height of
any critical point. If the interval $[c,c']$ does not contain a
critical value of the height function, then $\M^c$ is a strong
deformation retract of $\M^{c'}$, so the particular choices of $c$
above do not matter to the topology and hence to the evaluation of the
integral~(\ref{eq:cauchy}).

Let $X$ denote the topological pair $(\M^+ , \M^-)$.
The homology group $H_d (X)$ is generated by the homology
groups $H_d (\M^{c'},\M^{c''})$ as $[c',c'']$ ranges over 
a finite set of intervals whose union contains all critical values.
These groups in turn are generated by quasi-local cycles. At any
critical point $\zz^{(i)}$ in the top stratum, the quasi-local cycle
is the Cartesian product of a patch $\patch_i$ diffeomorphic to
$\R^{d-1}$ inside $\sing$, and whose maximum height is achieved at
$\zz^{(i)}$, with an arbitrarily small circle $\gamma_i$ transverse to
$\sing$.

\begin{figure}[ht] \hspace{1.5in}
\includegraphics[scale=0.6]{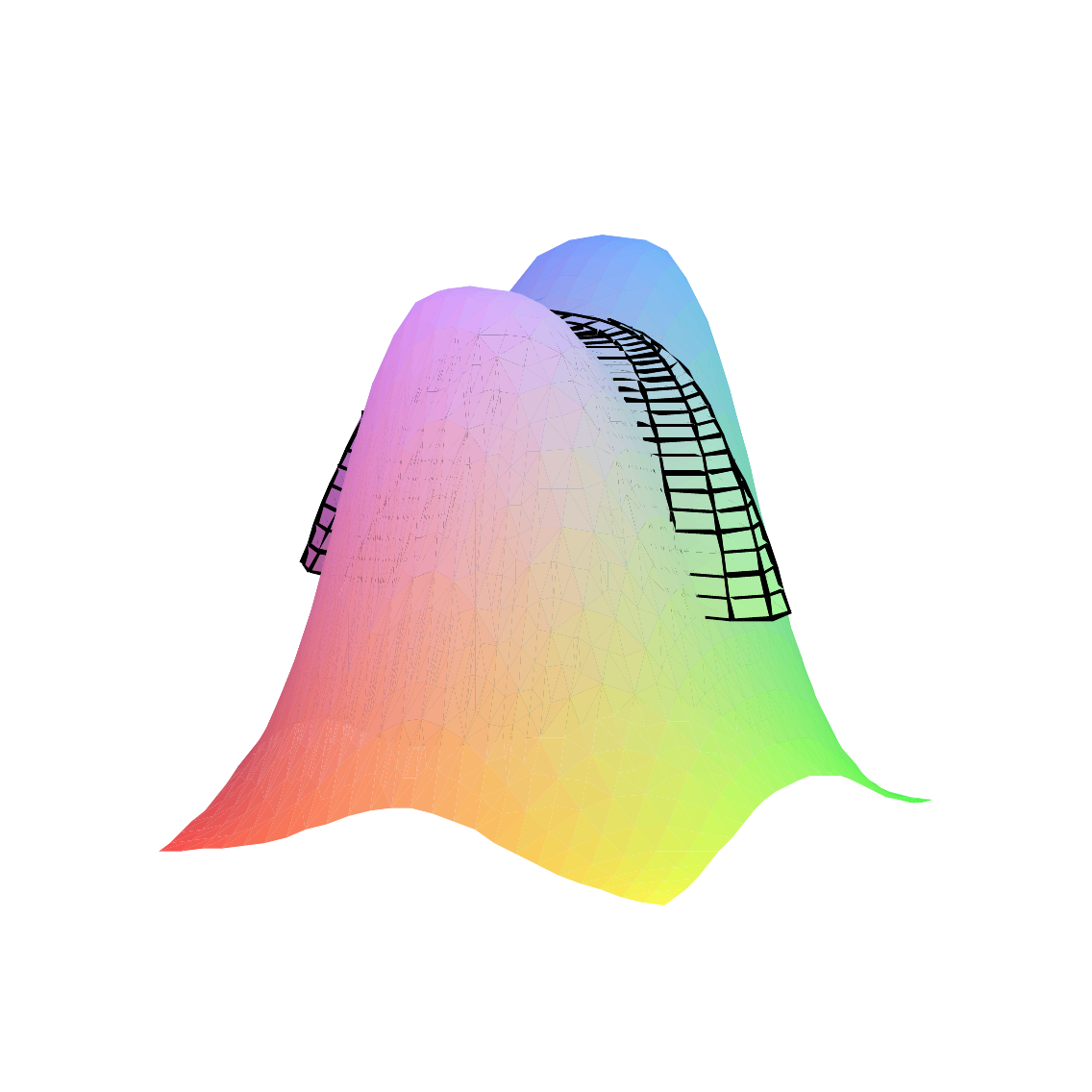}
\caption{a piece of a quasi-local cycle at top-dimensional critical point}
\label{fig:macaroni} 
\end{figure}

Thus for example when $d=2$, the quasi-local cycles near smooth points
look like pieces of macaroni: a product of a circle with a rainbow
shaped arc whose peak is at $\zz^{(i)}$.  We have not described the
quasi-local cycles for non-maximal strata, but for a description of
the quasi-local cycles in general, one may look in~\cite{BaPe}
or~\cite{GoMa1988}\footnote{For generic $\rbar$, the
function $h$ is Morse.  Each attachment map is a $d$-dimensional
complex, so the attachments induce injections on $H_d(X)$.  It may
happen for some $\rbar$ that $h$ is not Morse, but in this case one
may understand the local topology via a small generic perturbation.
Two or more critical points may merge, but the fact that the
attachments induce injections on $H_d(X)$ implies that such a merger
produces a direct sum in $H_d(X)$; in particular, it will be useful
later to know that a cycle in the merger is nonzero if and only if at
least one component was nonzero. \label{ft:1}}.

Since the quasi-local cycles generate $H_d (\M^+ , \M^-)$, it follows
that the integral~(\ref{eq:cauchy}) is equal to an integer linear
combination 
\begin{equation} 
\label{eq:decomp}
\sum n_i \int_{\CC_i} \zz^{-\rr - \one} F(\zz) \, d\zz
\end{equation}
where $\CC_i$ is a quasi-local cycle near $\zz$ and $T = \sum n_i
\CC_i$ in $H_d (\M^+ , \M^-)$.  When $\zz$ is in the top stratum,
i.e., $\zz$ is a smooth point of $\sing$, then $\CC_i = \patch_i
\times \gamma_i$, the product of a $d$-patch with a transverse circle,
so the integral may be written as $\int_{\patch_i} \int_{\gamma_i}
\exp [ -|\rr| h(\zz) ] F(\zz) \, d\zz $. The inner integral is a
simple residue and the outer integral is a standard saddle point
integral.  The asymptotic evaluation of $a_\rr$ is therefore solved if
we can compute the integers $n_i$ in the decomposition of $T$ into
$\sum n_i \CC_i$. Observe that the contribution from $\zz^{(i)}$ is of
exponential order no greater than $\exp (h(\zz^{(i)}))$.  It turns out
(Theorem~\ref{th:smooth asym} and the formulae of
Section~\ref{ss:multiple}) that for smooth and multiple points, this
is a lower bound on the exponential order as well, provided, in the
multiple point case, that $\numer (\zz^{(i)}) \neq 0$.  Thus if $n_i
\neq 0$, one may ignore any contributions from $\zz^{(j)}$ with
$h(\zz^{(j)}) < h(\zz^{(i)})$ and still obtain an asymptotic expansion
containing all terms not exponentially smaller than the leading term.
\begin{defn}[contributing critical points] \label{def:contrib}
The set $\contrib = \contrib_{\rbar}$ of contributing critical points
is defined to be the set of $\zz^{(i)}$ such that $n_i \neq 0$ and
$n_j = 0$ for all $j$ with $h(\zz^{(j)}) > h(\zz^{(i)})$.  In other
words, the contributing critical point(s) are just the highest ones
with nonzero coefficients, $n_i$.  Note that if there is more than one
contributing critical point, all must have the same height.
\end{defn}

The problem of computing the topological decomposition may in general
be difficult and no complete solution is known.  In two dimensions,
in the case where $\sing$ is globally smooth, an effective algorithm
is known for finding the contributing critical points. The algorithm
follows approximate steepest descent paths; the details are not yet
published~\cite{vdPe}.  In any dimension, if the variety $\sing$ is
the union of hyperplanes, it is shown in~\cite{BaPe} how to evaluate all
the $n_i$; see also some preliminary work on this
in~\cite{BeMc1993}. In the case where $\sing$ is not
smooth, while no general solution is known, we may still state a
sufficient condition for the critical point $\zz^{(i)}$ to be a
contributing critical point. It is to these geometric sufficient
conditions that we turn next.

\subsection{Geometric criteria} \label{ss:geom}

Let $F, \numer, \denom$ and $a_\rr$ be as in~(\ref{eq:F}).  
If $a_\rr \geq 0$ for all $\rr$, we say we are in the 
\Em{combinatorial case}.  We assume throughout
that $\numer$ and $\denom$ are relatively prime in the 
ring of analytic functions on the domain of convergence of $F$.
We wish to compute the function $\contrib$ which maps directions in
$\orthant$ to finite subsets of $\sing$ (often singletons).  The
importance of \Em{minimal points}, defined below, is that when
$\contrib_{\bar{\rr}}$ contains a minimal point, this point may be 
identified  by a variational principle.
\begin{defn}[notation for polydisks and tori] \label{def:polydisk}
Let $\poly (\zz)$ and $\torus (\zz)$ denote respectively the polydisk 
and torus defined by
\begin{eqnarray*}
\poly (\zz) & := & \{ \zz' : |z_i'| \leq |z_i| \mbox{ for all } 
   1 \leq i \leq d \} \; ; \\[1ex]
\torus (\zz) & := & \{ \zz' : |z_i'| = |z_i| \mbox{ for all } 
   1 \leq i \leq d \} \, .
\end{eqnarray*}
The open domain of convergence of $F$ is denoted $\dom$ and is
the union of tori $\torus (\zz)$.  The logarithmic domain 
of convergence, namely those $\xx \in \R^{d}$ with
$(e^{x_1}, \ldots , e^{x_{d}}) \in \dom$, is denoted $\logdom$ 
and is always convex~\cite{Horm1990}.   The
image $\{ \log |\zz| : \zz \in \sing \}$ of $\sing$ under the 
coordinatewise log-modulus map is denoted $\logsing$ (this is
sometimes called the \Em{amoeba} of $\sing$~\cite{GKZ1994}).  
\end{defn}
\begin{ass} \label{ass:hole}
We assume throughout that $\logdom$ is strictly convex, that is,
its boundary contains no line segment.  
\end{ass}
This ensures that $h_{\rhat}$ is always uniquely minimized on 
$\partial \logdom$.  The assumption is satisfied in all examples
in this paper. It is a consequence of $\logsing$ containing no 
line segments, which is true in all examples outside of 
Section~\ref{ss:integer}.  This may be checked by computer algebra;
we do not include the checking of this fact in our worked examples.
\begin{defn}[minimality] \label{def:minimal}
A point $\zz \in \sing$ is \Em{minimal} if all coordinates are nonzero
and the relative interior of the associated polydisk contains no
element of $\sing$, that is, $\sing \cap \poly (\zz) \subseteq
\partial \poly (\zz)$.  More explicitly, $z_i \neq 0$ for all $i$ and
there is no $\zz' \in \sing$ with $|z_i'| < |z_i|$ for all $i$.  The
minimal point $\zz$ is said to be \Em{strictly minimal} if it is the
only point of $\sing$ in the closed polydisk: $\sing \cap \poly (\zz)
= \{ \zz \}$.
\end{defn}

\begin{unremark}
This definition is equivalent to the apparently stronger definition
of minimality stated in~\cite{PeWi2002}, namely that $\sing \cap \disk
(\zz) \subseteq \torus (\zz)$ -- see Proposition~\ref{pr:linear}.
\end{unremark}

\begin{pr}[minimal points in $\contrib$] 
\label{pr:minimal}
Recall the standing assumption~\ref{ass:hole}.
\begin{enumerate}
\romenumi
\item A point $\zz \in \sing$ with nonzero coordinates is minimal 
if and only if $\zz \in \partial \dom$.
\item All minimal points in $\contrib_{\rbar}$ must lie on the
same torus and the height $h_{\rhat} (\zz)$ of the
point(s) in $\contrib_{\rbar}$ are at most the minimum, $c$, of the
height function $h_{\rr}$ on $\partial \dom$.  
\item If $\zz \in \contrib_{\rbar} \cap \partial \dom$, then the
hyperplane normal to $\rbar$ at $\log |\zz|$ is a support hyperplane
to $\logdom$.
\end{enumerate}
To summarize, if $\contrib_{\rbar}$ contains any minimal points, then
these points minimize $h_\rr$ on $\dom$ and they all project, via
coordinatewise log-modulus, to a single point of $\partial \logdom$
where the support hyperplane is normal to $\rr$.
\end{pr}

\begin{proof}
A power series converges absolutely on any open polydisk about the
origin in which it is analytic; this can be seen by using Cauchy's
formula~(\ref{eq:cauchy}) to estimate $a_\rr$.  Thus any minimal point
is in the closure of $\dom$.  A series cannot converge at a pole, so
$\zz \notin \dom$ and the first assertion follows.

To prove the second assertion, by strict convexity there is a 
unique $\zz \in \partial \dom$ that minimizes $h$.  Let $\homot$ 
be the homotopy $\{ t \torus (\zz) : \ee \leq t < 1-\ee \}$.  
This may be extended to a homotopy $\homot'$ that pushes
the height below $c - \ee$ except in small neighborhoods of $\sing
\cap \torus(\zz)$.  Thus the small torus $T$ in~(\ref{eq:cauchy}) is
homotopic to a cycle in the union of $\M^{c-\ee}$ and a neighborhood
of $\sing \cap \torus (\zz)$. Sending $\ee \to 0$, it follows that
$\contrib$ has height at most $c$. If there is a line segment in
$\partial \logdom$, it is possible that the minimum height on $\dom$
occurs on more than one torus, but in this case, if $\zz$ and $\ww$
are on different such tori, the above argument shows that neither
$\zz$ nor $\ww$ can be in $\contrib_{\rbar}$.

The third assertion is immediate from the fact that $\log |\zz|$
minimizes the linear function $h_{\rhat}$ on the convex set $\logdom$.
\end{proof}

Having characterized $\contrib_{\rbar}$ when it is in $\partial \dom$,
via a variational principle, we now relate this to the algebraic
definition of $\crit_{\rbar}$ which is better for symbolic
computation.  It is easy to write down the inverse of $\crit_{\rbar}$,
which we will denote by $\linear$.  
\begin{defn}[geometry of minimal points]
\label{def:classify}
Say that $\zz$ is a \Em{smooth point} of $\sing$ if $\sing$ is a
manifold in a neighborhood of $\zz$.  Say that $\zz$ is a \Em{multiple
point} of $\sing$ if $\sing$ is the union of finitely many manifolds
near $\zz$ intersecting transversely (that is, normals to any $k$ of
these manifolds span a space of dimension $\min \{ k,d \}$).  Points
that are neither smooth nor multiple are informally called bad points.
\end{defn}
\begin{unremark}
The transversality assumption streamlines the arguments but is
not really needed (see, e.g.,~\cite{BaPe}).
\end{unremark}
\begin{defn}[the linear space $\linear$] \label{def:linear}
Let $\zz \in \sing$ be in a stratum $\stratum$. 
Let $\linear (\zz) \subseteq \CP^{d-1}$ denote the span of the
projections of vectors $(z_1 v_1 , \ldots , z_{d} v_{d})$
as $\vv$ ranges over vectors orthogonal to the tangent space 
of $\stratum$ at $\zz$.  
\end{defn}

\begin{figure}[ht] \hspace{0.5in}
\includegraphics[scale=0.6]{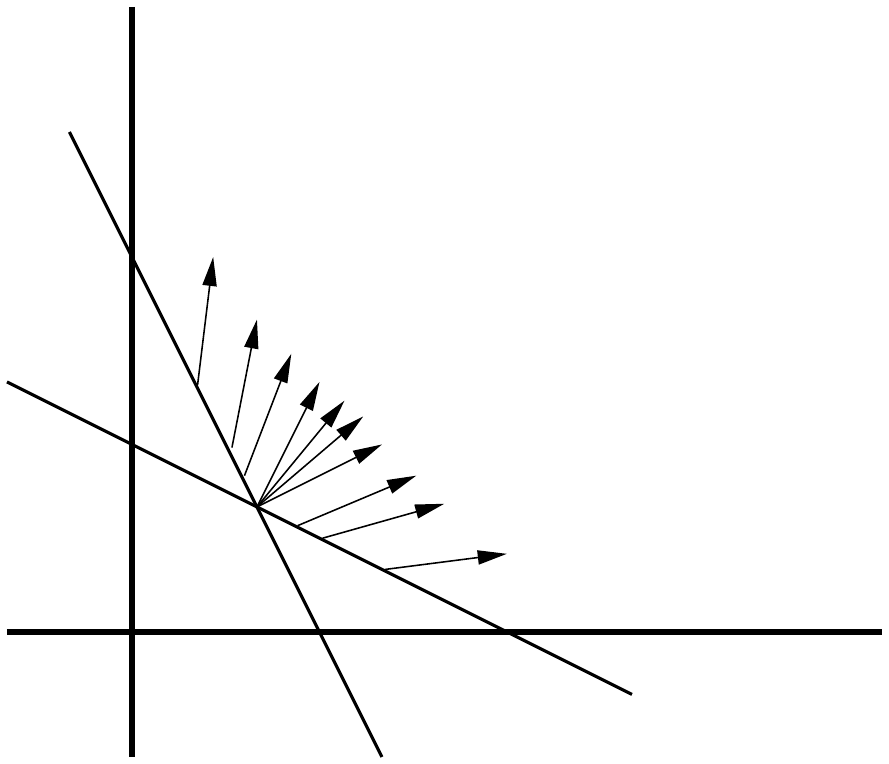}
\caption{direction $\linear (\zz)$ for positive real points of $\sing$ when
$1/F = (3-x-2y)(3-2x-y)$}
\label{fig:lines}
\end{figure}

The Figure~\ref{fig:lines} gives a pictorial example of the foregoing
definitions (some parts of the figure refer to Definition~\ref{def:cone} 
below).  The singular variety $\sing$ for the function
$F(x,y) = 1/[(3-x-2y)(3-2x-y)]$ is composed of two complex lines
meeting at $(1,1)$.  The crossing point $\zz_0 := (1,1)$ is a 
singleton stratum.  The following proposition, showing why
$\linear$ is important, follows directly from the definitions.

\begin{pr}[$\linear$ inverts $\crit$] \label{pr:crit}
The point $\zz$ is in $\crit_{\rbar}$ if and only if 
$\rbar \in \linear (\zz)$.  If $\zz$ is a smooth point, then 
(replacing $\denom$ by its radical if needed), 
$\linear (\zz)$ is the singleton set (in projective space)
$$
\linear (\zz) = \left \{ \loggrad \denom (\zz) \right \} := 
   \left \{ \left ( z_1 \frac{\partial \denom}{\partial z_1} (\zz) , 
   \ldots , z_{d} \frac{\partial \denom}{\partial z_{d}} (\zz) 
   \right ) \right \}
$$
while the set $\crit_{\rbar}$ is the (usually zero-dimensional) 
variety given by the $d$ equations
\begin{align} 
\denom & =  0 \nonumber \\[3ex]
r_{d} z_j \frac{\partial \denom}{\partial z_j} & =  r_j z_{d} 
   \frac{\partial \denom}{\partial z_{d}} \qquad (1 \leq j \leq d-1) 
   \label{eq:zero-dim}
\end{align}
If $\zz$ is a multiple point, then $\linear (\zz)$ is the 
span of the vectors
$$
\loggrad \denom_k (\zz) := \left ( z_1 \frac{\partial\denom_k}{\partial z_1} 
   (\zz) , \ldots , z_{d} \frac{\partial\denom_k}{\partial z_{d}} (\zz)
   \right )
$$
where $\denom = \prod_k \denom_k$.
$\Cox$
\end{pr}

Connecting this up to the variational principle, we have:
\begin{pr}[smooth minimal points are critical]
\label{pr:linear}
If $\zz$ is minimal and smooth then $\linear (\zz) \in \orthant$
and is normal to a support hyperplane to $\logdom$ at $\log |\zz|$.
Consequently, a minimal smooth point $\zz$ is a critical point for
some outward normal direction to $\logdom$ at $\log |\zz|$.
\end{pr}

\begin{proof} Assume first that none of the logarithmic partial 
derivatives $h_j := z_j \partial \denom / \partial z_j$ vanishes 
at $\zz$.  Suppose the arguments of two of the partials $h_k$ and $h_l$ 
are not equal.  Since we have assumed no $h_j$ vanishes, there
is a tangent vector $\vv$ to $\sing$ with $v_j = 0$ for 
$j \notin \{ k , l \}$ and $v_k \neq 0 \neq v_l$.  We may 
choose a multiple of this so that $\overline{v_k} h_k$ and 
$\overline{v_l} h_l$ both have negative real parts.  Perturbing
slightly, we may choose a tangent vector $\vv$ to $\sing$
at $\zz$ with all nonzero coordinates so that $\overline{v_j} h_j$
has negative real part for all $j$.  This implies there is a
path from $\zz$ on $\sing$ such that the moduli of all
coordinates strictly decrease, which contradicts minimality.
It follows that all arguments of $h_j$ are equal, and therefore 
that $\linear (\zz)$ is a point of $\orthant \subseteq \CP^d$;
hence it is normal to a support hyperplane to $\logdom$ at $\log |\zz|$.

Now, suppose $h_j = 0$ for $j$ in some nonempty set, $J$.
Let $\xx = \log |\zz|$.  Vary $\{ z_j : j \notin J \}$ by  
$\Theta(\ee)$ so as to stay on $\sing$, varying $\{ z_j : j \in J \}$ 
by $O(\ee^2)$.  We see that the complement of $\logdom$ contains
planes arbitrarily close to the $|J|$-dimensional real plane through 
$\xx$ in directions $\{ e_j : j \in J \}$, whence the closure of
the complement contains this plane.  But by monotonicity, the orthant 
of this plane in the $-e_j$ directions must be in the closure of
$\logdom$.  Thus this orthant is in $\partial \logdom$, whence there
is a lifting of this orthant to $\sing$.  Now the argument by
contradiction in the previous paragraph shows that if the arguments 
of $h_k$ and $h_l$ are unequal for some $k,l \notin J$, then we may move 
on $\sing$ so the moduli of the $z_j$ for $j \notin J$ strictly decrease,
while moving down the lifting of the orthant allows us also to
decrease the moduli of the $z_j$ for $j \in J$, and we have again
contradicted minimality.  
\end{proof}

Having related critical points to $\linear (\zz) \subseteq \CP^{d-1}$, 
we now relate contributing critical points to a set $\cone (\zz) 
\subseteq \RP^{d-1}$.  An open problem is to find the right definition
of $\cone (\zz)$ when $\zz$ is not minimal.  

\begin{defn} [the cone $\cone$ of a minimal point] \label{def:cone}
If $\zz$ is a smooth point of $\sing$ in $\partial \dom$, let
$\dir (\zz)$ be defined as $\linear (\zz)$,
which is in $\orthant$ by Proposition~\ref{pr:linear}.
For minimal $\zz \in \sing$, not necessarily smooth, define 
$\cone (\zz)$ to be the convex hull of the set of limit points 
of $\dir (\ww)$ as $\ww \to \zz$ through smooth points.
\end{defn}

For an example of this, see Figure~\ref{fig:lines}, in which
arrows are drawn depicting $\cone (\zz)$ 
for smooth points $\zz$ in the positive real quadrant and 
the two-dimensional cone $\cone (\zz_0)$ is shaded.
The next proposition then follows from Proposition~\ref{pr:crit}
and basic properties of convex sets.
\begin{pr}[description of $\cone$]
\label{pr:cone}
Let $\zz$ be a minimal point.
\begin{enumerate}
\romenumi
\item $\cone (\zz) \subseteq \linear (\zz) \cap \orthant$.
\item If $\zz$ is a smooth or multiple point, then $\cone (\zz)$ 
is the cone spanned by the vectors $\loggrad \denom_k$ of 
Proposition~\ref{pr:crit}.
\item If a neighborhood of $\zz$ in $\sing$ covers a neighborhood
of $\log |\zz|$ in $\partial \logdom$, then $\cone (\zz)$ is the 
set of outward normals to $\logdom$ at $\log |\zz|$.
\end{enumerate}
\end{pr}

The first nontrivial result we need is:
\begin{thm}[$\cone$ inverts $\contrib$] \label{th:contrib}
Let $\zz$ be a minimal point and either 
smooth or multiple.  Then $\rbar \in \cone (\zz)$ if and only if 
$\zz \in \contrib_{\rbar}$, provided, in the multiple point case,
that $\numer (\zz) \neq 0$.
\end{thm}

\begin{unremark}
Here, rather than global meromorphicity, it is only required
that $F$ be meromorphic in a neighborhood of $\poly(\zz)$; 
see~\cite{PeWi2002} and~\cite{PeWi2004} for details. 
\end{unremark}

\begin{proof} Assume first that $\zz$ is strictly minimal,
that is, $\torus (\zz) \cap \sing = \{ \zz \}$.  

Suppose for now that $\rbar$ is in the relative interior of $\cone
(\zz)$. Theorems~\cite[Theorem~3.5]{PeWi2002}
and~\cite[Theorems~3.6~and~3.9]{PeWi2004} give expressions for $a_\rr$ which
are of the order $|\rr|^\beta \zz^{-\rr}$, relying, in the multiple
point case, on the assumption of transverse intersection, which we have
built into our definition of multiple point, and on $\numer (\zz) \neq
0$.  

On one hand, it is evident from~(\ref{eq:cauchy}) that 
$$
a_\rr = O \left ( \exp \left [ |\rr| (h_{\rhat} (\ww) + \ee) 
   \right ] \right )
$$
for any $\ee > 0$, where $\ww \in \contrib_{\rbar}$.  On the
other hand, we know by Proposition~\ref{pr:minimal} that no 
points higher than $\zz$ are in $\contrib_{\rbar}$.  Since
$\contrib_{\rbar}$ is nonempty, we conclude that it has points
at height $h_{\rhat} (\zz)$, whence from the expressions for $a_\rr$
again, we see that $\zz \in \contrib_{\rbar}$.  

In the case $\rbar \in \partial \cone (\zz)$ 
consider $\rbar' \to \rbar$ through the interior of $\cone (\zz)$
(here we mean relative boundary and relative interior).
There will be a contributing point $\zz'$ which converges to $\zz$.
The coefficient of the quasi-local cycle at $\zz$ in the limit 
must be nonzero (see footnote~\ref{ft:1}).  The theorem is now
proven for strictly minimal points.  

To remove the assumption of strict minimality, one must verify that 
this was not necessary for the formulae we quoted from~\cite{PeWi2002,PeWi2004}.
These formulae were proved by reducing the Cauchy integral to
an integral over a neighborhood $\nbd$ of a $(d-1)$-dimensional 
subset $\Eta$ of $\sing$.  It is pointed out~\cite[Corollary~3.7]{PeWi2002}
that it is sufficient to assume \Em{finite minimality}, that is,
finiteness of $\torus (\zz) \cap \sing$.  In fact, one needs only
finiteness of $\torus (\zz) \cap \Eta$, since the truncation of
the integral to $\nbd$ incurs a boundary term that is sufficiently
small as long as $\Eta$ avoids $\torus (\zz)$.  The remaining
case, where $\Eta \cap \torus (\zz)$ is infinite, can be handled
by contour rotation arguments, but since that work is not yet published,
we point out here that in the case $d=2$ (the only case used in this
survey), the set $\Eta$ is one-dimensional.  Thus
when $\Eta \cap \torus (\zz)$ is infinite, the set $\Eta$ must be 
a subset of $\torus (\zz)$; the integrals in~\cite{PeWi2002,PeWi2004} may
then be taken over all of $\Eta$ with no truncation. 
\end{proof}

We have remarked earlier that the main challenge in computing
asymptotics is to identify $\contrib_{\rbar}$.  Our progress to this
point is that we may find all strictly minimal smooth or multiple points
in $\contrib_{\rbar}$ by solving the equation $\rbar \in \cone (\zz)$
for $\zz$.  This equation may be solved automatically when $F$ is a
rational function and is often tractable in other cases, for example in
the case $F = (e^x - e^y)/(x e^y - ye^x)$ of Section~\ref{ss:euler}.  

To complement this, we would like to know when solving $\rbar \in \cone
(\zz)$ and checking for minimality does indeed find all points of
$\contrib_{\rbar}$.  It cannot, for instance, do so if there are no
minimal contributing points.  Thus we ask (a) are there any minimal
points $\zz$ with $\cone (\zz) = \rbar$, (b) are all the contributing
points minimal, and (c) is there more than one minimal point?  It turns
out that the answers are, roughly: (a) yes, in the combinatorial case,
(b) yes as long as (a) is true, by part~(ii) of
Proposition~\ref{pr:minimal}, and (c) rarely (we know they are never on
different tori).

Let $\dirset \subseteq \orthant$ denote the set of all normals
to support hyperplanes of $\logdom$.  If the restriction of $F$ 
to each coordinate hyperplane is not entire, then $\dirset$ is 
the entire nonnegative orthant.  This is because $\logdom$ has 
support hyperplanes parallel to each coordinate hyperplane.  When
$\dirset$ is not the whole orthant, then in directions $\rbar \notin
\orthant$, the quantity $a_\rr$ is either identically zero or 
decays faster than exponentially (this follows, for example, 
from~\cite[(1.2)]{PeWi2004}, since $- \rhat \cdot \log |\zz|$ is not
bounded from below on $\domain$).  
\begin{thm}[existence of minimal points in the combinatorial case]
\label{th:universal}
Suppose that we are in the combinatorial case, $a_\rr \geq 0$.  
Under the standing assumption~\ref{ass:hole},
for every $\rbar$ in the interior of $\dirset$, 
there is a minimal point $\zz \in \sing$ which lies in the positive
orthant of $\R^{d}$ and has $\rbar \in \cone(\zz)$.  
\end{thm}

\begin{unremark} 
At this point, it is clear how far one can generalize beyond
rational functions.  Given a compact set $K \subseteq \dirset$
in which one wishes to find asymptotics.  Let $\logdom(K)$ denote
the set of all $\xx$ such that $\xx \leq \yy$ coordinatewise 
for some $\yy \in \partial \logdom$ whose normal direction 
is in $K$.  Let $\dom (K)$ be the inverse image under the
coordinatewise log-modulus map of $\logdom (K)$.  Then, in the
combinatorial case, for the existence of a minimal point $\zz$ 
in the positive orthant with $\rbar \in \cone (\zz)$, 
it is sufficient that $F$ be meromorphic in a neighborhood 
of $\dom (K)$.  
\end{unremark}

\begin{proof}
We follow the proof of~\cite[Theorem~6.3]{PeWi2002}.  For any $\rbar$
in the interior of $\dirset$ there is a point $\xx \in \partial \logdom$ 
with a support hyperplane normal to $\rbar$.  Let $z_i = e^{x_i}$
so $\zz$ is a real point in $\partial \dom$.  We claim that 
$\zz \in \sing$.   To see this, note that there is a singularity on 
$\torus (\zz)$, which must be a pole by the assumption of meromorphicity 
on a neighborhood of $\poly(\zz)$.  Together, $a_\rr \geq 0$ and lack 
of absolute convergence of the series on $\torus (\zz)$ imply
that $F(\ww)$ converges to $+\infty$ as $\ww \to \zz$ from beneath.
By meromorphicity, $\zz$ is therefore a pole of $F$, so $\zz \in \sing$.

We conclude that there is a lifting of $\partial \logdom$ to 
the real points of $\sing$, that is, a subset of the real points of
$\sing$ maps properly and one to one onto $\partial \logdom$.  
By the last part of Proposition~\ref{pr:cone}, it follows that 
$\rbar \in \cone (\zz)$, and hence from Theorem~\ref{th:contrib} 
that $\zz \in \contrib_{\bar{\rr}}$.  
\end{proof}

\begin{unremark}
By convexity of $\logdom$, the point $\zz \in \orthant^d \cap \crit(\rbar)$
in Theorem~\ref{th:universal} is unique unless there is a line segment
in $\partial \logdom$. In this case there is a continuum of such $\zz$.
Since in almost every example the critical point variety is finite,
there will be precisely one critical point in the positive orthant. Thus
in practice we will have no difficulty in determining the point $\zz$. 
\end{unremark}

We say that a power series $P$ is \Em{aperiodic} if the sublattice of
$\Z^{d}$ of integer combinations of exponent vectors of the 
monomials of $P$ is all of $\Z^{d}$.  By a change of variables,
we lose no generality from the point of view of generating functions
if we assume in the following proposition that $P$ is aperiodic.

\begin{pr}[often, every minimal point is strictly minimal]
\label{pr:classify}
If $\denom = 1 - P$ where $P$ is aperiodic with nonnegative
coefficients, then every minimal point is strictly minimal and lies 
in the positive orthant.
\end{pr}

\begin{proof}
Suppose that $\zz = \xx e^{i \btheta}$, with $\btheta
\neq \zero$, is minimal. Since $\zz \in \sing$, we have 
\begin{eqnarray*}
1 & = & \left| \sum_{\rr \neq \zero} a_\rr \xx^\rr
\exp(i \rr \cdot \btheta) \right| \\
& \leq & \sum_{\rr \neq \zero} a_{\rr} \xx^\rr \, .
\end{eqnarray*}
Thus $\zz$ can be strictly minimal only if $\btheta = \zero$.

Equality holds in the above inequality if and only if either there is
only one term in the latter sum, or, in case there is more than one such
term, if $\exp(i \rr \cdot \btheta) = 1$ whenever $a_\rr \neq
0$.  The former case occurs precisely when $P$ has the form
$c \xx^\rr$ and the latter can occur on when $P$ has the form
$P(\zz) = g(\zz^{\bf b})$ for some ${\bf b} \neq 0$, both of which
are ruled out by aperiodicity. 
\end{proof}

We restate the gist of Theorem~\ref{th:universal} and 
Proposition~\ref{pr:classify} as the following corollary.
\begin{cor} \label{cor:procedure}
If $\denom = 1 - P$ with $P$ aperiodic and having nonnegative
coefficients, then the contributing critical points as $\rbar$ varies
are precisely the points of $\sing \cap \mathcal{O}^{d}$ that are minimal
in the coordinatewise partial order.  The point $\zz$ is in
$\contrib_{\rbar}$ exactly when $\rbar \in \cone (\zz)$.

Without the assumption $\denom = 1-P$, assuming only $a_\rr \geq 0$, the
same holds except that there might be more contributing critical points
on the torus $\torus (\zz)$. 
\noproof
\end{cor}

\subsection{Asymptotic formulae for minimal smooth points}
\label{ss:smooth} 

It is shown in~\cite{PeWi2002} how to compute the integral~(\ref{eq:cauchy})
when $\rbar$ is fixed and $\contrib_{\rbar}$ is a finite set of smooth
points on a torus.  Given a smooth point $\zz \in \sing$, let $\ft_\zz$ 
be the map on a neighborhood of the origin in $\R^{d-1}$ taking the origin 
to zero and taking $(\theta_1 , \ldots , \theta_{d-1})$ to $\log w$ 
for $w$ such that 
$$\left ( z_1 e^{i \theta_1} , \ldots , z_{d-1} e^{i \theta_{d-1}} , 
   z_{d} w \right ) \in \sing \, .$$
The following somewhat general result is shown in~\cite[Theorem~3.5]{PeWi2002}.

\begin{thm}[smooth point asymptotics] \label{th:smooth asym}
Let $K \subset \dirset$ be compact, and suppose that for
$\rbar \in K$, the set $\contrib_{\rbar}$ is a single smooth
point $\zz (\rbar)$ and that $\ft_\zz$ has nonsingular Hessian
(matrix of second partial derivatives).  
Then there are effectively computable functions $b_l (\rbar)$
such that 
\begin{equation} \label{eq:smooth asym}
a_\rr \sim \zz (\rbar)^{-\rr} \sum_{l \geq 0} b_l (\rbar) (r_d)^{-(l+d-1)/2}
\end{equation}
as an asymptotic expansion when $|\rr| \to \infty$, uniformly 
for $\rbar \in K$.
\noproof
\end{thm}

\begin{unremark}
In connection with Theorem~\ref{th:smooth asym}, the following points should be noted.
\begin{enumerate}
\romenumi
\item
The coefficients $b_l$ depend on the derivatives of $\numer$ 
and $\denom$ to order $k+l-1$ at $\zz (\rbar)$, and $b_0 (\rbar) = 0$
if and only if $\numer (\zz (\rbar)) = 0$. 
\item  The sum on the right of equation~\ref{eq:smooth asym} may be rewritten as a sum of 
$b_l^* (\rbar) |\rr|^{-(l+d-1)/2}$ where $b_l^* = \rhat_d^{-(l+d-1)/2} b_l$,
in order to see what depends on $\rhat$ and what depends on $|\rr|.$
\item
Suppose $\contrib_{\rbar}$ is a finite set of points $\zz^{(1)} , 
\ldots \zz^{(k)}$, satisfying the hypotheses of the theorem.  Then
one may sum~(\ref{eq:smooth asym}) over these, obtaining
$$a_\rr \sim \sum_k \zz^{(k)} (\rbar)^{-\rr} 
   \sum_{l \geq 0} b_{k,l}^* (\rbar) |\rr|^{-(l+d-1)/2} \, .$$
\end{enumerate}
\end{unremark}
Theorem~\ref{th:smooth asym} is somewhat messy: even when $l=1$, 
the combination of partial derivatives of $\numer$ and $\denom$ is 
cumbersome, though the prescription for each $b_l$ in terms of 
partial derivatives of $\numer$ and $\denom$ is completely algorithmic
for any $l$ and $d$.  In the applications herein, we will confine our 
asymptotic computations to the leading term.  The following expression 
for the leading coefficient is more explicit than~(\ref{eq:smooth asym}) 
and is proved in~\cite[Theorems~3.5 and~3.1]{PeWi2002}.

\begin{thm}[smooth point leading term] \label{th:smooth nice}
When $\numer (\zz(\rbar)) \neq 0$, the leading coefficient is
\begin{equation} \label{eq:smooth}
b_0 = (2 \pi)^{(1-d)/2} \hess^{-1/2} \frac{\numer (\zz(\rbar))}{- z_{d} 
   \partial \denom / \partial z_{d}}
\end{equation}
where $\hess$ denotes the determinant of the
Hessian of the function $\ft_\zz$ at the origin.  
\end{thm}
\begin{cor} \label{th:smooth d=2}
In particular, when $d=2$ and $\numer (\zz(\rbar)) \neq 0$, we have
\begin{equation} \label{eq:d=2}
a_{rs} \sim \frac{\numer (x,y)}{\sqrt{2 \pi}} x^{-r} y^{-s}
   \sqrt{\frac{- y \denom_y}{s Q(x,y)}}
\end{equation}
where $(x,y) = \zz (\overline{(r,s)})$ and $Q(x,y)$ is defined 
to be the expression
\begin{equation} \label{eq:Q}
- (x\denom_x)(y\denom_y)^2  - (y\denom_y)(x\denom_x)^2 -  \left[ 
   (y\denom_y)^2 x^2\denom_{xx} + (x\denom_x)^2 y^2\denom_{yy} - 2
   (x\denom_x) (y \denom_y) xy\denom_{xy} \right]. 
\end{equation}
$\Cox$
\end{cor}

\begin{unremark}
In the combinatorial case, the expression in the radical will be 
positive real (this is true for more general $F$ with correct choice 
of radical).  The identity $r y H_y = s x H_x$ (see Proposition~\ref{pr:crit})
shows that the given expression for $a_{rs}$, though at first 
sight asymmetric in $x$ and $y$, has the expected symmetry. 
\end{unremark}

\subsection{Multiple point asymptotics} 
\label{ss:multiple}

Throughout this section, $\zz$ will denote a minimal
multiple point.  We will let $m$ denote the number of sheets 
of $\sing$ intersecting at $\zz$ and let 
\begin{equation} \label{eq:multi}
\denom = \prod_{j=1}^m \denom_k^{n_k}
\end{equation}
denote a local representation of $\denom$ as a product of nonnegative
integer powers of functions whose zero set is locally smooth.  
For $1 \leq k \leq m$, let $\vb^{(k)}$ denote the vector 
whose $j^{th}$ component is $z_j \partial \denom_k / \partial z_j$.  

We divide the asymptotic analysis of $a_\rr$ into two cases, namely $m
\geq d$ and $m < d$.  In the former case, since we have assumed
transverse intersection at multiple points, the stratum of $\zz$ is
just the singleton $\{ \zz \}$.  The asymptotics of $a_\rr$ are
simpler in this case.  In fact it is shown
in~\cite[Theorem~3.1]{Pema2000} that 
\begin{equation} \label{eq:highpole}
a_\rr = \zz^{-\rr} \left ( P(\rr) + E(\rr) \right )
\end{equation}
where $P$ is piecewise polynomial and $E$ decays exponentially
on compact subcones of the interior of $\cone (\zz)$.  We begin with 
this case.

We will state three results in decreasing order of generality. 
The first is completely general, the second holds in the special 
case $m = d$ and $n_k = 1$ for all $k \leq m$, and the last holds
for $m=d=2$ and $n_1 = n_2 = 1$.  Our version of the most general
result provides a relatively simple formula from~\cite{BaPe} but
under the more general scope of~\cite{PeWi2004}.  

\begin{defn} \label{def:poly}
Let $M = \sum_{k=1}^m n_k$ and for $1 \leq j \leq M$ let $t_j$
be integers so that the multiset $\{ t_1 , \ldots , t_M \}$ 
contains $n_k$ occurrences of $k$ for $1 \leq k \leq m$.
Define a map $\map : \R^M \to \R^d$ by $\map (e_j) = \vb^{(t_j)}$.
Let $\lambda^M$ be Lebesgue measure on $\mathcal{O}^M$ and let $P(\xx)$
be the density at $\xx$ of the pushforward measure $\lambda^M \circ 
\map$ with respect to Lebesgue measure under $\map$.  The function 
$P$ is a piecewise polynomial of degree $M - d$, the regions of 
polynomiality being no finer than the common refinement of 
triangulations of the set $\{ \vb^{(k)} : 1 \leq k \leq m \}$ 
in $\RP^{d-1}$~\cite[definition~5]{BaPe}.
\end{defn}

The above definition of $P$ is a little involved but is made clear 
by the worked example in Section~\ref{ss:integer}.  Armed with this
definition, we may say what happens when there are $d$ or more sheets 
of $\sing$ intersecting at a single point. 

\begin{thm}[isolated point asymptotics] \label{th:m>=d}
If $m \geq d$ and $\zz$ is a minimal point in a 
singleton stratum, with $\numer (\zz) \neq 0$, 
then uniformly over compact subcones of $\cone (\zz)$,
\begin{equation} \label{eq:M>=d}
a_\rr \sim \numer (\zz) P \left ( \frac{r_1}{z_1} , \ldots ,
\frac{r_d}{z_d} \right ) \zz^{-\rr} \, ,
\end{equation}
provided that $\contrib_{\rbar}$ contains only $\zz$.  This formula
may be summed over $\contrib_{\rbar}$ as long as $\contrib_{\rbar}$
has finite cardinality. 

In the case $\zz = \one$,~(\ref{eq:M>=d}) reduces further to
$$a_\rr \sim \numer (\one) P(\rr) \, .$$
If, furthermore, $a_\rr / \numer (\one)$ are integers, 
then $a_\rr / \numer (\one)$ is actually a piecewise polynomial whose 
leading term coincides with that of $P(\rr)$.  
\end{thm}

\begin{proof} First assume $\zz$ is strictly minimal.  
Theorem~3.6 of~\cite{PeWi2004} gives an asymptotic expression for $a_\rr$
valid whenever $\zz$ is a strictly minimal multiple point in a
singleton stratum.  The formula, while not impossible to use, is not
as useful as the later formula given in~\cite[equation~(3.8)]{BaPe}.
This latter equation was derived under the assumption that $\denom$ is
a product of linear polynomials.  Noting that the formula
in~\cite[Theorem~3.6]{PeWi2004} depends on $\denom$ only through the
vectors $\{ \vb^{(k)} : 1 \leq k \leq m \}$, we see it must agree
with~\cite[equation~(3.8)]{BaPe} which is~(\ref{eq:M>=d}).  The last
statement of the theorem follows from~(\ref{eq:highpole}).

To remove the assumption of strict minimality, suppose that
$\zz$ is minimal but not strictly minimal and that there 
are finitely many points $\zz^{(1)} , \ldots , \zz^{(k)}$
in $\contrib_{\rbar}$, necessarily all on a torus.  The torus 
$T$ in~(\ref{eq:cauchy}) may be pushed out to $(1 + \ee) \torus (\zz)$
except in neighborhoods of each $\zz^{(j)}$.  The same sequence
of surgeries and residue computation in the proofs 
of~\cite[Proposition~4.1, Corollary~4.3 and Theorem~4.6]{PeWi2004}
now establish Theorem~3.6 of~\cite{PeWi2004} without the hypothesis of
strict minimality, and as before one obtains~(\ref{eq:M>=d})
from~(3.8) of~\cite{BaPe}.
\end{proof}

While this gives quite a compact representation of the leading term,
it may not be straightforward to compute $P$ from its definition as
a density.  When $M = m = d$, $P$ is a constant and the computation 
may be reduced to the following formula.

\begin{cor}[simple isolated point asymptotics] \label{cor:m=d}
If $m=d$ and each $n_k = 1$ in~(\ref{eq:multi}), then if 
$\numer (\zz) \neq 0$, and $\zz$ is a minimal point in a 
singleton stratum, then
$$
a_\rr \sim \numer (\zz) \zz^{-\rr} \left | J \right |^{-1}
$$ 
where $J$ is the Jacobian matrix $(\partial \denom_i / \partial z_j)$.
This formula may be summed over finitely many points as in the remark following Theorem~\ref{th:smooth asym}.
\noproof
\end{cor}

Occasionally it is useful to have a result that does not depend on
finding an explicit factorization of $\denom$.  The matrix $J$ can
be recovered from the partial derivatives of $\denom$.  The result
of this in the case $m = d = 2$ is given by the following 
formula~\cite[Theorem~3.1]{PeWi2004}.

\begin{cor}[simple isolated point asymptotics, dimension $2$]
\label{cor:m=d=2}
If $M=m=d=2$, and if $\numer (\zz) \neq 0$, then setting $\zz =
(x,y)$,
$$
a_{rs} \sim x^{-r} y^{-s} \frac{\numer (x,y)}{\sqrt{-x^2 y^2 \hess}}
$$
where $\hess := \denom_{xx} \denom_{yy} - \denom_{xy}^2$ is the 
determinant of the Hessian of $H$ at the point $(x,y)$.  
In the special but frequent case $x=y=1$ we have simply
$$
a_{rs} \sim \frac{\numer (1,1)}{\sqrt{-\hess}} \, .
$$
This formula may be summed over finitely many points as in the 
second remark following Theorem~\ref{th:smooth asym}.
\noproof
\end{cor}

Finally, we turn to the case $m < d$.  In this case, $\zz$ is in a
stratum containing more than just one point.  The leading term
of $a_\rr$ is obtained by doing a saddle point integral of a formula 
such as~(\ref{eq:M>=d}) over a patch of the same dimension as the
stratum.  Stating the outcome takes two pages in~\cite{PeWi2004}.  Rather
than give a formula for the resulting constant here, we state the
asymptotic form and refer the reader to~\cite[Theorem~3.9]{PeWi2004} for 
evaluation of the constant, $b_0$.

\begin{thm} \label{th:m<d}
Suppose that as $\rbar$ varies over a compact subset $K$ of
$\orthant$, the set $\contrib_{\rbar}$ is always a singleton varying
over a fixed stratum of codimension $m < d$.  If also $\zz (\rbar)$
remains a strictly minimal multiple point, with each smooth sheet
being a simple pole of $F$, and if $\numer (\zz) \neq 0$ on $K$, then
\begin{equation} \label{eq:m<d}
a_\rr \sim \zz (\rbar)^{-\rr} b_0 (\rbar) |\rr|^{\frac{m-d}{2}}
\end{equation}
uniformly as $|\rr| \to \infty$ in $K$.
\noproof
\end{thm}

\subsection{Distributional limits}
\label{ss:limit theorems}

Of the small body of existing work on multivariate asymptotics, a good
portion focuses on limit theorems.  Consider, for example, the point
of view taken in~\cite{Bend1973,BeRi1983} and the sequels to
those papers~\cite{GaRi1992,BeRi1999}, where one
thinks of the the numbers $a_\rr$ as defining a sequence of
$(d-1)$-dimensional arrays, the $k^{th}$ of which is the
\Em{horizontal slice} $\{ a_\rr : r_d = k \}$
(cf.~Section~\ref{ss:GF-seq}). Often this point of view is justified
by the combinatorial application, in which the last coordinate, $r_d$,
is a size parameter, and one wishes to understand rescaled limits of
the horizontal slices. Distributional limit theory assumes nonnegative
weights $a_\rr$, so for the remainder of this section we assume we are
in the combinatorial case, $a_\rr \geq 0$.

If ${C_k := \sum_{\rr \, : \,  r_d = k} a_\rr < \infty}$ 
for all $k$, we may define the slice distribution $\mu_k$ on 
$(d-1)$-vectors to be the probability measure giving mass 
$a_\rr / C_k$ to the vector $(r_1 , \ldots , r_{d-1})$ in the $k^{th}$ 
slice.  There are several levels of limit theorem that may be of interest.  
A weak law of large numbers (WLLN) is said to hold if the measures
$A \mapsto \mu_k (\frac{1}{k} A)$ converge to a point mass at 
some vector $\mm$ (here, division of $A$ by $k$ means division of
each element by $k$).  Equivalently, a WLLN holds if and only if
for some mean vector $\mm$, for all $\ee > 0$,
\begin{equation} \label{eq:WLLN}
\lim_{k \to \infty} \mu_k \left \{ \rr : \left | \frac{\rr}{k} - \mm 
   \right | > \ee \right \} = 0 \, .
\end{equation}

Stronger than this is Gaussian limit behavior.  As $\rr$
varies over a neighborhood of size $\sqrt{k}$ about $k \mm$
with $r_d$ held equal to $k$, formula~(\ref{eq:smooth asym})
takes on the following form:
$$a_\rr \sim C |\rr|^{(-(d-1)/2} \exp \left [ \phi (\rr) \right ] \, .$$
Here the factor $C |\rr|^{-(d-1)/2}$ varies only by $o(1)$ when
 $\rr$ varies by $O(\sqrt{k})$, and 
$$\phi (\rr) := -\rr \cdot \log |\zz (\rr) | \, ,$$
or just $- \rr \cdot \log \zz (\rr)$ when $\zz (\rr)$ is real.
The WLLN identifies the location $\rr_{\max}$ at which 
$\phi$ takes its maximum for $r_d$ held constant at $k$.
A Taylor expansion then gives 
\begin{equation} 
\label{eq:gaussian}
\exp \left [ \phi (\rr_{\max} + \sss) - \phi (\rr_{\max}) \right ] 
   =  \exp \left [ - \frac{B(\sss)}{2k} \right ] 
\end{equation}
for some nonnegative quadratic form, $B$.  Typically $B$ will
be positive definite, resulting in a local central limit theorem
for $\{ \mu_k \}$.  

Unfortunately, the computation of $B$, though straightforward from
the first two Taylor terms of $\denom$, is very messy.  Furthermore,
it can happen that $B$ is degenerate, and there does not seem to
be a good test for this.  As a result, existing limit theorems 
such as~\cite[Theorems~1 and~2]{BeRi1983} 
and~\cite[Theorem~2]{GaRi1992} all contain the hypothesis ``if
$B$ is nonsingular''.  Because of the great attention that has been
paid to Gaussian limit results, we will state a local central limit
result at the end of this section.  But since we cannot provide a
better method for determining nonsingularity of $B$, we will not
develop this in the examples.  

Distributional limit theory requires the normalizing constants $C_k$
be finite.  We make the slightly stronger assumption that $a_\rr = 0$ 
when $\rr / r_d$ is outside of a compact set $K$.  We define the map 
$$
\rr \mapsto \rr^* := \left ( \frac{r_1}{r_d} , \ldots , 
   \frac{r_{d-1}}{r_d} \right ) \, ,
$$
which will be more useful for studying horizontal slices than 
was the the previously defined projection $\rr \mapsto \rhat$.
\begin{thm}[WLLN] \label{th:WLLN}
Let $F$ be a $d$-variate generating function with nonnegative
coefficients, and with $a_\sss = 0$ for $\sss / s_d$ outside some
compact set $K$. Let $f(x) = F(1, \ldots , 1 , x)$.  Suppose $f$ has a
unique singularity $x_0$ of minimal modulus which is a simple pole.
Let $\rbar = \dir (\xx)$, where $\xx = (1 , \ldots , 1 , x_0)$.  If
$F$ is meromorphic in a neighborhood of $\disk (\xx)$, then the measures
$\{ \mu_k \}$ satisfy a WLLN with mean vector $\rr^* := (r_1 / r_d ,
\ldots , r_{d-1} / r_d)$.
\end{thm}

\begin{proof} 
We begin by showing that $\logdom$ has a unique normal at $\log
|\xx|$, which is in the direction $\rbar$.  We know that $x_0$ is
positive real, since $f$ has nonnegative coefficients. The proof of
Theorem~\ref{th:universal} showed that a neighborhood $\nbd$ of $\xx$
in $\sing$ maps, via log-moduli, onto a neighborhood of $\log |\xx|$
in $\partial \logdom$.  Since $x_0$ is a simple pole of $f$, we see
that $\nbd$ is smooth with a unique normal $\grad \denom (\xx)$,
whence $\logdom$ has the unique normal direction $\rbar$ at $\log
|\xx|$.

We observe next that $\sum_k C_k x^k = f(x)$, whence, via
Theorem~\ref{th:darboux}, $C_k \sim c x_0^{-k}$. Thus
Theorem~\ref{th:WLLN} will follow once we establish:
\begin{equation} \label{eq:exp decay}
\limsup_k \frac{1}{k} \log \left ( \sum_{\sss : s_d = k , \sss^* \in K
, |\sss^* - \rr^*| > \ee} a_\sss \right ) < - \log x_0  \, . 
\end{equation}

Each $\sss^* \in K$ with $|\sss^* - \rr^*| \geq \ee$ is not
normal to a support hyperplane at $\log |\xx|$, so for each such
$\sss^*$ there is a $\yy \in \partial \logdom$ for which 
\begin{equation} 
\label{eq:-log x}
h_{\sss^*} (\yy) < h_{\sss^*} (\xx) = - \log x_0 \, .
\end{equation}
By compactness of $\{ \sss^* \in K : |\sss^* - \rr^*| \geq \ee\}$
we may find finitely many $\yy$ such that one of these, denoted 
$\yy (\sss^*)$, satisfies~(\ref{eq:-log x}) for any $\sss^*$.  
By compactness again, 
$$
\sup_{\sss^*} h_{\sss^*} (\yy (\sss^*)) < - \log x_0 \, .
$$

It follows from the representation~(\ref{eq:cauchy}), choosing 
a torus just inside $\torus (\yy)$ for $\yy$ as chosen above 
depending on $\sss$, that $|\sss|^{-1} \log a_\sss$ is at most 
$- \log x_0 - \delta$ for some $\delta > 0$ and all sufficiently 
large $\sss$ with $|\sss^* - \rr^*| \geq \ee$.  Summing over
the polynomially many such $\sss$ with $s_d = k$ 
proves~(\ref{eq:exp decay}) and the theorem.   
\end{proof}

Coordinate slices are not the only natural sections on which 
to study limit theory.  One might, for example study limits 
of the sections $\{ a_\rr : |\rr| = k \}$, that is,
over a foliation of $(d-1)$-simplices of increasing size.  The
following version of the WLLN is adapted to this situation.
The proof is exactly the same when slicing by $|\rr|$ or
any other linear function as it was for slicing by $r_d$.

\begin{thm}[WLLN for simplices] \label{th:WLLN2}
Let $F$ be a $d$-variate generating function with nonnegative 
coefficients.  Let $f(x) = F(x, \ldots , x)$.  Suppose $f$ has a 
unique singularity $x_0$ of minimal modulus which is a simple
pole.  Let $\rbar = \dir (\xx)$, where $\xx = (x_0 , \ldots , x_0)$.  
If $F$ is meromorphic in a neighborhood of $\disk (\xx)$, then 
the measures $\{ \mu_k' \}$, defined analogously to $\mu_k$ but 
over the simplices $\{ |\rr| = k \}$, satisfy a WLLN with 
mean vector $\rhat$. 
$\Cox$
\end{thm}

We end with a statement of a local central limit theorem.
\begin{thm}[LCLT] \label{th:LCLT}
Let $F, f, x_0, \xx , \rr$ and $\rr^*$ be as in Theorem~\ref{th:WLLN}.
Suppose further that $\contrib_{\rbar}$ is the singleton, $\{ \xx \}$.
If the quadratic form $B$ from~(\ref{eq:gaussian}) is nonsingular,
then
\begin{equation} \label{eq:LCLT}
\lim_{k \to \infty} k^{(d-1)/2} \sup_\rr |\mu_k (\rr) - \normalB (\rr)| = 0
\end{equation}
where 
\begin{equation} \label{eq:normal}
\normalB (\rr) := (2 \pi)^{-(d-1)/2} \det (B)^{-1/2} 
   \exp \left ( -\frac{1}{2} B(\rr) \right ) \, .
\end{equation}
is the discrete normal density.   
$\Cox$
\end{thm}

\setcounter{equation}{0}
\section{Detailed examples} 
\label{ss:details}

In Sections~\ref{ss:details}~--~\ref{sec:kernel} we work a multitude
of examples.  We aim to cover a sufficient variety of examples 
so that a user with a new application is likely to find 
a worked example that is similar to his or her own.  With luck,
the derivation of the result will be adaptable, and the user can
thereby avoid original sources as well as the bulk of 
Section~\ref{ss:results} of the present survey.  The use of
computer algebra is essential to many of these examples, and
we preface the worked examples with a brief explanation of this.

To apply any of the theorems of 
Sections~\ref{ss:smooth}~--~\ref{ss:limit theorems}, one must
solve for the point $\zz \in \sing$ as a function of the
direction $\rbar$ and then plug this into a variety of formulae.
This may always be done numerically, but there are advantages
to doing it analytically when possible: one may then differentiate,
compute asymptotics near a point of interest, solve in terms of
other data, and so forth.  Even when one is interested chiefly in
a numerical approximation for a single $a_\rr$, one often obtains
better numerics by simplifying analytically as much as possible.

When the generating function is rational or algebraic, as is quite
frequent in applications, if one attempts to solve for $\zz$ and plug
in the results, one often finds that a computer algebra system cannot
simplify an expression that one suspects should be simpler.
Consequently, the computer cannot tell whether such an expression is
zero, may have trouble plotting the expression, and so forth.  As is
well known in computational commutative algebra, these problems can be
avoided by working directly with polynomial ideals.  This survey
incorporates the Groebner package in the computer algebra system Maple
to illustrate the required computational algebra.  We hope that these
commands from Maple version~10 will remain in the Maple platform, but
in any case, we include platform-independent explanations of the
algebra.  Our treatment is necessarily quite brief, and readers who
require more explanation should consult a source such as~\cite{CLO2005,
Stur2005}.

Suppose a point $\zz \in \C^d$ is the solution to polynomial equations
$P_1 (\zz) = \cdots = P_k (\zz) = 0$.  The set of all polynomials $P$
for which $P(\zz) = 0$ is an \Em{ideal}, and, if the number of common
solutions is finite, is said to be a \Em{zero-dimensional} ideal.
Such an ideal has many generating sets or \Em{bases}, of which some
are particularly useful to know.  In particular, given a \Em{term
order}, that is, an order on the monomials $\zz^\rr$ that obeys
certain properties, each ideal has a unique \Em{Gr\"obner basis}.  The
lexicographic term order, known to Maple as {\tt plex}, has the
property that it forces the Gr\"obner basis of a zero-dimensional
ideal to contain a univariate polynomial in the variable designated to
come last in the lexicographical order.  Thus, if $p_1, p_2$ and $p_3$
are polynomials in $x, y$ and $z$ with finitely many common solutions,
then the commands (the first just loads the Gr\"obner basis package
into Maple)
\begin{verbatim}
> with(Groebner);
> L := Basis([p1,p2,p3],plex(x,z,y));
\end{verbatim}
will produce a basis for the ideal of all polynomials vanishing
on the common solutions to $p_1 , p_2$ and $p_3$.  The last variable
in the variable list is $y$, so the choice of term order {\tt plex}
causes the last polynomial of the basis $L$ to be a polynomial in 
$y$ alone, whose roots are the possible $y$-coordinates of the
common solution points.  

A few words may prove helpful concerning the classification 
of critical points into smooth, multiple and bad points.  A
stratum always consists of points of a single topological type.
The critical point equations are different for different strata.
Therefore, sorting critical points by type occurs as a preliminary
step of sorting strata by type.  If one assumes the denominator
$\denom$ of $F$ is square free (if not, a single Maple command allows
one to pass to the radical), then the smooth points of $\sing$
are precisely those where $\grad \denom$ does not vanish.  The 
first step is always to check whether there are any non-smooth 
points; generically, there are not, though interesting
applications tend not to be generic.  This step is accomplished 
by the single command
\begin{verbatim}
> Basis ([H , diff(H,x) , diff(H,y) , diff(H,z)] , tdeg(x,y,z));
\end{verbatim}
Maple returns the trivial ideal {\tt [1]} if and only if there
are no non-smooth points on $\sing$.  We have used {\tt tdeg}
rather than the slower {\tt plex} to check whether we get {\tt [1]};
if not, then we can go back and use {\tt plex} to produce a more
useful basis for the ideal corresponding to the variety of 
singular points.  Further tests may be done to determine whether
a stratum of singular points consists of multiple points: one
must check whether the radical of the ideal is \Em{reduced};
for definitions and algorithms, the reader is referred to a text 
such as~\cite{Mumf1995}.

The algebra above is often sufficient to compute quantities such as
the quantity $\crit_{\rbar}$, which is defined by the polynomial
equations~(\ref{eq:zero-dim}).  One may, however, go further: the
substitution of one of these finitely many points into a polynomial
expression such as~(\ref{eq:Q}) results in an algebraic number which
Maple usually cannot write in its simplest form.  Later, we discuss
two reasonable ways to handle this.  Given a quantity $x$ defined by
algebraic equations $p_1 , \ldots , p_k$, the most straightforward
way to simplify the polynomial expression $Q(x)$ is to reduce
$Q$ modulo the ideal generated by $\{ p_1 , \ldots , p_k \}$
(see Section~\ref{ss:delannoy}).  An alternative is to 
obtain directly the minimal polynomial satisfied by $Q(x)$.
This may be done by elimination (see Section~\ref{ss:distinct subseq})
or by matrix representation (see Section~\ref{ss:fine}).
These techniques are useful for almost every example, but to avoid
being repetitive, we often refer back to these fully worked
computations.
 
\subsection{Binomial coefficients} \label{ss:binomial}

To allow the reader a chance to check the application of 
Corollary~\ref{th:smooth d=2} in a familiar setting, we begin with an
example where the numbers $a_{rs}$ are explicitly known.
Let 
$\displaystyle{a_{rs} = \binom{r+s}{r,s} = \frac{(r+s)!}{r! \, s!}}$
be the binomial coefficients and let 
$$
F(x,y) = \sum_{r,s \geq 0} a_{rs} x^r y^s \, .
$$
The binomial coefficients satisfy the recurrence $a_{r,s} = a_{r-1,s}
+ a_{r,s-1}$ and this holds for all $(r,s) \neq (0,0)$, provided that
we take $a_{rs}$ to be zero when either $r$ or $s$ is negative. Linear
recursions for $a_{rs}$ in terms of values $a_{r',s'}$ with $(r',s')
\leq (r,s)$ in the coordinatewise partial order lead easily to
rational generating functions (as discussed in
Section~\ref{sec:kernel}, more general linear recursions do not yield
rational functions). To find $F$, we observe from the recursion that
all the coefficients of $(1 - x - y) F$ vanish except the
$(0,0)$-coefficient, which is~1. Thus $$ F(x,y) = \frac{1}{1-x-y} \,
.$$

Let us compute the set $\contrib_{\rbar}$ and the quantities appearing
in Corollary~\ref{th:smooth d=2}.  The singular variety $\sing$ is the
complex line $x + y = 1$.  The numerator of $F$ is~1, and in
particular, it never vanishes.  For any $\zz \in \sing$, the space
$\linear (\zz)$ is the linear span of $\zz$.  To see this,
either use Proposition~\ref{pr:crit} or note that the tangent space to
$\sing$ is everywhere orthogonal to $(1,1)$ and plug this into
Definition~\ref{def:linear}.

For each direction $\rbar$ in the positive real orthant, there is thus
a unique solution $\zz \in \sing$ to $\rbar \in \linear (\zz)$,
namely $\zz = \rhat = (\frac{r}{r+s} , \frac{s}{r+s})$ for any
representative $(r,s)$ of $\rbar$.  One may apply
Theorem~\ref{th:universal} to conclude that 
$\contrib_{\rbar} = \{ \rhat \}$.  

To apply Corollary~\ref{th:smooth d=2}, we need only compute the
quantity $Q$ in~(\ref{eq:Q}) and verify that it is nonzero. We have
$\denom = 1 - x - y, \denom_x = -1 , \denom_y = -1$ and all other
partial derivatives of $\denom$ are zero, whence $Q = -xy(x + y)$ and
plugging in $(x,y) = (\frac{r}{r + s} , \frac{s}{r + s})$ gives
$$
\frac{-y \denom_y}{s Q(x,y)} = \frac{1}{s x (x + y)} = \frac{r +
s}{rs} \, .
$$
Substituting this and $\numer \equiv 1$ into~(\ref{eq:d=2}) gives
\begin{equation} 
\label{eq:binomial}
a_{rs} \sim \left ( \frac{r+s}{r} \right )^r 
   \left ( \frac{r+s}{s} \right )^s \sqrt{\frac{r+s}{2 \pi r s}} \, .
\end{equation}
This asymptotic expression is valid as $(r,s) \to \infty$, uniformly
if $r/s$ and $s/r$ remain bounded -- see the uniformity conclusion in
Theorem~\ref{th:smooth asym}.  In fact, we know from Stirling's
formula that this holds uniformly as $\min \{ r , s \} \to \infty$;
see~\cite{Llad} to obtain the latter result in the present
framework.

Although we have not presented precise error bounds in our asymptotic
approximations, the errors for such smooth bivariate problems can be
shown to be of the order $1/s$ and our approximations are good even
for moderate values of $r$ and $s$. As a randomly chosen  example, we
observe that the approximation above yields a relative error of about
$0.8\%$ when $r = 25, k = 12$, while the relative error has decreased
to about $0.4\%$ when $r = 50, k = 24$ and $0.08\%$ when $r = 250, s =
120$.

\subsection{Delannoy numbers} \label{ss:delannoy}

Recall the Delannoy numbers $a_{rs}$ that count paths from
the origin to the point $(r,s)$ with each step having 
displacement $(1,0)$, $(0,1)$ or $(1,1)$.  In exactly the
same way that we obtained the generating function for
the binomial coefficients, we may use the recursion
$a_{r,s} = a_{r-1,s} + a_{r,s-1} + a_{r-1,s-1}$ valid 
for all $(r,s)$ except $(0,0)$ to obtain the generating function
$$
F(x,y) := \sum_{r,s \geq 0} a_{rs} x^r y^s = \frac{1}{1-x-y-xy}
$$
so the denominator is given by $\denom = 1 - x - y - xy$.  

Solving for $\contrib_{\rbar}$ is only a little more involved for this
generating function than it was for $1/(1 - x - y)$. However, performing
all calculations by hand, as we did in the case of binomial
coefficients, is tedious for this example and completely impractical for
later examples. The step-by-step computation below illustrates how a
computer algebra system such as Maple can carry the derivation to completion at a
symbolic level.

The linear space $\linear (x,y)$ is the one-dimensional complex
vector space spanned by $(x \denom_x , y \denom_y) = (- x(1 + y) , 
- y (1 + x))$. 
The most convenient way in which to solve the equations
$\denom = 0$ and $\rbar \in \linear (x,y)$ is to obtain a 
Gr\"obner basis.  We load Maple's Groebner
package using the command \texttt{with(Groebner)} and then execute the
following command. 
\begin{verbatim}
> GB:=Basis([H, s*x*diff(H, x) - r*y*diff(H, y)], plex(x,y));
\end{verbatim}
This computes a basis for the ideal corresponding to the common
solutions to $\denom = 0$ and $\overline{(r,s)} \in \linear(x,y)$. The
answer is $GB:=[-s+sy^2+2ry, s-sy-r+rx] =:[p_1(y), p_2(x,y)] $. The
first element of the basis is an elimination polynomial in $y$ alone.
Solving $p_1 = 0$ yields the two values $y = (-r \pm \sqrt{r^2 + s^2})
/ s$.  For each of the two $y$-values there is a unique $x$-value
obtained by solving $p_2$ for $x$.  One might observe that symmetry
tells us these will be $(-s \pm \sqrt{r^2 + s^2}) / r$; however we
need $p_2$ to tell us which $x$-value goes with which $y$-value.  We
find that there is a unique positive solution
\begin{eqnarray*}
x (\rbar) & = & \frac{\sqrt{r^2 + s^2} - s}{r} \, ; \\
y (\rbar) & = & \frac{\sqrt{r^2 + s^2} - r}{s} \, .
\end{eqnarray*}
By Theorem~\ref{th:universal} we see that $\contrib_{\rbar} =
(x(\rbar) , y(\rbar))$ for any positive direction $\rbar$. The
numerator, $\numer$, is again identically~1. It remains only to plug
the values for $x$ and $y$ into expressions for $Q$ and $a_{rs}$ and
simplify.

Using the definition of $Q$ and simplifying, we obtain $Q = x y (1 +
x) (1 + y) (x + y)$. Let us reduce $Q$ modulo $H$ by the command 
\begin{verbatim}
NormalForm (Q, [H, s*x*diff(H ,x) - r*y*diff(H, y)], plex(x,y));
\end{verbatim}
to get
\begin{equation} 
\label{eq:Q delannoy}
Q = -2x^2 - 2y^2 + 6x + 6y - 4 \, .
\end{equation}
Finally, substituting $x = x(\rbar)$ and $y = y(\rbar)$ into
$-y \denom_y / (sQ)$ yields
$$
\frac{- y \denom_y}{s Q} = \frac{rs}{\sqrt{r^2 + s^2}(r+s-\sqrt{r^2 +
s^2})^2 
}
$$

and putting this all together yields  the
expression~(\ref{eq:delannoy}) for the Delannoy numbers:
$$
a_{rs} \sim \left ( \frac{\sqrt{r^2 + s^2} - s}{r} \right)^{-r} \left (
\frac{\sqrt{r^2 + s^2} - r}{s} \right )^{-s} \sqrt{\frac{1}{2 \pi}}
\sqrt{\frac{r s}{\sqrt{r^2 + s^2}(r+s-\sqrt{r^2 + s^2})^2 }} \, .
$$
For example, putting $r = s = n$ shows that the central Delannoy
numbers $a_{nn}$ have first order asymptotic approximation 
$$
a_{nn} \sim \frac{\cosh (\frac{1}{4} \log 2)}{\sqrt{\pi}}(3 + 2
\sqrt{2})^n n^{-1/2}.
$$
\subsection{Powers, quasi-powers, and generalized Riordan arrays} 
\label{ss:riordan}

It often happens that we wish to estimate $[z^n] v(z)^k$, that is,
the $n^{th}$ coefficient of a large power of a given function $v(z)$. 
Clearly, this is equal to the $x^n y^k$ coefficient
of the generating function
\begin{equation} \label{eq:riordan 1}
F(x,y) := \frac{1}{1 - y v(x)} \, .
\end{equation}
One place where this arises is in the enumeration of a combinatorial class 
whose objects are strings built from given blocks. Let $v(z) :=
\sum_{n=1}^\infty a_n z^n$ count the number $a_n$ of blocks of size
$n$. Then the generating function~(\ref{eq:riordan 1}) counts objects
of a given size by the number of blocks in the object. 
\begin{example} \label{eg:prefixes}
A long sequence of zeros and ones may be divided into blocks
by repeatedly stripping off the unique initial string that 
is a leaf of $T$, a given \Em{prefix tree}.  Lempel-Ziv coding,
for instance, does this but with an evolving prefix tree.  
When $v(x)$ is the generating function for the number $a_n$ 
of leaves of $T$ at depth $n$, then $1/(1- y v(x))$ generates 
the numbers $a_{rs}$ of strings of length $r$ made of $s$ blocks 
(the final block must be complete).
\end{example}

Another place where generating functions of this form arise is in the
Lagrange inversion formula. This application is discussed at length in
Section~\ref{sec:lagrange}, but briefly, if $h$ solves the equation
$h(z) = z v(h(z))$, then even if we cannot explicitly solve for $h$,
its coefficients are given by 
$$[z^n] h (z) = \frac{1}{n} [z^{n-1}] v(z)^n \, .$$
This identity has been very profitable in the analysis of 
planar graphs and maps, cf.\ the discussion of results
of~\cite{GaWo1999,BFSS2001} in Section~\ref{ss:core}.

A third place where coefficients of powers arise is in sums of
independent random variables.  Let $v(z) = \sum_{n = 0}^\infty a_n
z^n$ be the probability generating function for a distribution on the
nonnegative integers, that is, $a_n = \P (X_j = n)$ where $\{ X_j \}$
are a family of independent, identically distributed random variables.
Then $v(z)^n$ is the probability generating function for the partial
sum $S_n := \sum_{j=1}^n X_j$, and hence 
$$
\P (S_n = k) = [z^k] v(z)^n \, .
$$
A \Em{Riordan array} is defined to be an array $\{ a_{nk} \,:\, n,k
\geq 0\}$ whose generating function $F(x,y) := \sum_{n,k \geq 0}
a_{nk} x^n y^k$ satisfies 
\begin{equation} 
\label{eq:riordan 2}
F(x,y) = \frac{\phi (x)}{1 - y v(x)}
\end{equation}
for some functions $\phi$ and $v$ with $v(0) = 0$ and $\phi (0) \neq
0$.  If in addition $v'(0) \neq 0$ the array is called a \Em{proper}
Riordan array.  Just as~(\ref{eq:riordan 1}) represents sums of independent, identically distributed
random variables when $v$ is a probability generating function, the
format~(\ref{eq:riordan 2}) generalizes this to \Em{delayed renewal}
sums (see, e.g.~\cite[Section~3.4]{Durr2004}), where an initial summand
$X_0$ may be added that is distributed differently from the others.
The quasi-powers~(\ref{eq:quasi-power}) arising in GF-sequence
analysis, which were described in Section~\ref{ss:GF-seq}, are
asymptotically of this form as well. Thus~(\ref{eq:riordan 2})
approximately encompasses most of the known results leading to
Gaussian behavior in multivariate generating functions.

Riordan arrays have been widely studied.  In addition to enumerating a
great number of combinatorial classes, Riordan arrays also behave in
an interesting way under matrix multiplication (note that the
condition $v(0) = 0$ implies $a_{nk} = 0$ for $k < n$, and, by
triangularity of the infinite array, that multiplication in the
Riordan group is well defined).  Surveys of the Riordan group and its
combinatorial applications may be found
in~\cite{Spru1994,SGWW1991}.

As we will see in this section, the asymptotic analysis of these
arrays is relatively simple, or at least, is no more difficult than
analyses of the functions $\phi$ and $v$ that define the array. One
should note, however, that Riordan arrays are often defined by data
other than $\phi$ and $v$.  Commonly, one has a linear recurrence for
$a_{n,k+1}$ as a sum $\sum_{s=1}^{k-n} c_s a_{n+s,k}$; the generating
function $A(t) := \sum_{j=1}^\infty c_j t^j$ is known and is fairly
simple, but the function $v(x)$ in~(\ref{eq:riordan 2}) is known only
implicitly through an equation~(\cite[Equation
6]{Roge1978})
$$
v(x) = x A(v(x))
$$
which is reminiscent of Lagrange inversion.  
The paper~\cite{Wils2005}, which is devoted to bivariate
asymptotics of Riordan arrays, discusses at length how to proceed when
the known data includes $A(t)$ rather than $v(x)$.  Thus there are
versions of~(\ref{eq:riordan-asymp}) below, such as
Proposition~\ref{pr:lagrange-bivariate}, which state asymptotics in
terms of $A$ without explicit mention of $v$.  The discussion in this
section, however, will be limited to deriving asymptotics in terms of
$\phi$ and $v$.  For our analyses it is not important to require $v(0)
= 0$ (for example, neither the binomial coefficient nor Delannoy
number examples above satisfies that condition), so we drop this
hypothesis and consider generalized Riordan arrays that
satisfy~(\ref{eq:riordan 2}) but may have $v(0) \neq 0$.

The following computations show that the quantities $\linear (x,y)$
and $Q(x,y)$ turn out to be relatively simply expressed in terms
of the function $v(x)$.  Define the quantities
\begin{align}
\mu (v; x) & := \frac{x v' (x)}{v(x)} \label{eq:mu} \, ; \\[2ex]
\sigma^2 (v; x) & := \frac{x^2 v''(x)}{v(x)} + \mu (v; x) - \mu(v;
x)^2 = x \frac{d \mu(v; x)}{d x} \, . \label{eq:sigma}
\end{align}
It is readily established that for $(x,y) \in \sing$, we have $\linear
(x,y) = \overline{(\mu(v; x), 1)}$.  In other words, $(x, 1/v(x)) \in
\crit_{(r,s)}$ if and only if $s \mu (v; x) = r$. Furthermore, when
this holds, $Q(x, 1/v(x)) = \sigma^2(v; x)$. Provided that $\phi$ and
$\sigma^2$ are nonzero at a minimal point, the leading term of its
asymptotic contribution in~(\ref{eq:d=2}) then becomes
\begin{equation}
\label{eq:riordan-lt}
a_{rs} \sim x^{-r} v(x)^{s}\frac{\phi(x, 1/v(x))}
   {\sqrt{2 \pi s \sigma^2(v; x)}}
\end{equation}
where $\mu(v; x) = r/s$.

The notations $\mu$ and $\sigma^2$ are of course drawn from
probability theory.  These quantities are always nonnegative when $v$
has nonnegative coefficients. To relate this to the limit theorems in
Section~\ref{ss:limit theorems}, observe that setting $x = 1$ gives
\begin{eqnarray}
\mu (v ; 1) & = & \frac{v'(1)}{v(1)} \, ; \label{eq:riordan-mu}\\
\sigma^2 (v ; 1) & = & \frac{v v'' - (v')^2 + v v'}{v^2} (1) \, . 
\label{eq:riordan-sigma}
\end{eqnarray}
Thus, under the hypotheses of Theorem~\ref{th:WLLN}, a WLLN will
hold with mean $\mm = \mu (v ; 1)$.  Of course we see here that
$\mu (v ; 1)$ is simply the mean of the renormalized distribution
on the nonnegative integers with probability generating function $v$.
Similarly, we see in Theorem~\ref{th:LCLT} that 
$B(r,s) = (s - \mu (v;1) r)^2 / \sigma^2 (v;1)$
is the Gaussian term corresponding to the variance $\sigma^2 (v;1)$
of this renormalized distribution.  

Suppose we are in the combinatorial, aperiodic case: the coefficients
of $v$ are nonnegative and $v(x)$ cannot be written as $x^b g(x^d)$
for any power series $g$ and $d > 1$.  The set $\sing$ is the union of
the set $\sing_0$, parametrized as $(x, 1/v(x))$, with the union of
horizontal lines where the value of $x$ is a singular value for $v$ or
$\phi$.  As $x$ increases from 0 to $R$, the (possibly infinite)
minimum of the radii of convergence of $v$ and $\phi$, it is easy to
verify that all the points of $\sing_0$ encountered are strictly
minimal, and that $\sigma^2 (v;x) > 0$.  Since the derivative of $\mu
(v;x)$ is $\sigma^2 (v;x)/x$, this shows that $\mu (v;x)$ is strictly
increasing on $(0,R)$.  Thus $A := \mu (v;0)$ and $B := \mu (v;R)$ are
well defined as one-sided limits and for $A < \lambda < B$ there is a
unique solution to $\mu (v;y) = \lambda$. In fact, $A$ is the
order of $v$ at 0, so equals~1 for proper Riordan arrays, but may be~0
for generalized Riordan arrays or an integer greater than~1 for
improper Riordan arrays.  From this it is evident that $a_{rs} = 0$
for $r/s < A$, so we should not have expected a solution to $\mu(v;x)
= \lambda$ when $\lambda < A$. If $R = \infty$, then $B$ is the
(possibly infinite) degree of $v$, and one has again that $a_{rs} = 0$
for $r/s > B$.  If $v$ or $\phi$ is not entire, then one cannot say
without further analysis what one expects for $a_{rs}$ when $r/s \to
\lambda > B$.  To summarize:
\begin{pr}
\label{pr:riordan}
Let  $(v(x), \phi(x))$ determine a generalized Riordan array.
Suppose that $v(x)$ has radius of convergence $R > 0$ and is 
aperiodic with nonnegative coefficients, and that $\phi$ has 
radius of convergence at least $R$.  Let $A, B$ be as above.
Then for $A < r/s < B$,
\begin{equation} 
\label{eq:riordan-asymp}
a_{rs} \sim v(x)^{s} x^{-r} s^{-1/2} 
\sum_{k = 0}^\infty b_k (r/s) s^{-k}
\end{equation} 
where $x$ is the unique positive real solution to $\mu(v; x) = r/s$. 
Here $b_0 = \frac{\phi(x)}{\sqrt{2\pi \sigma^2(v; x)}} \neq 0$. 
The asymptotic approximation is uniform as $r/s$ varies within
a compact subset of $(A, B)$, whereas $a_{rs} = 0$ for $r/s < A$. 
\noproof
\end{pr}

We note that if the combinatorial restriction is lifted, much more
complicated behavior can occur.  The generating function 
$3/(3 - 3x - y + x^2)$ is of Riordan type with $\phi(x) = v(x) = 
(3 - 3x + x^2)^{-1}$, and even though $v$ is aperiodic,
$\contrib_{\rr}$ has cardinality $2$. Furthermore, at the unique
contributing point for the diagonal direction, $\sigma^2$ vanishes.

The condition on the radius of convergence of $\phi$ is satisfied 
in most applications.  One way in which it may fail is when $F$ is a 
product of more than one factor.  We will see an example of this
in Section~\ref{ss:distinct subseq}.  The next two subsections 
consider applications of Proposition~\ref{pr:riordan} to 
combinatorial applications.

\subsection{Maximum number of distinct subsequences}
\label{ss:distinct subseq}

Flaxman, Harrow and Sorkin~\cite{FHS2004} consider
strings over an alphabet of size $d$ which we take to be  $\{1, 2,
\dots ,d\}$ for convenience.  They are interested in strings of length
$n$ which contain as many distinct subsequences (not necessarily
contiguous) of length $k$ as possible.  Let $a_{nk}$ denote the
maximum number of distinct subsequences of length $k$ that can be
found in a single string of length $n$.  Initial segments $S|_n$ of
the infinite string $S$ consisting of repeated blocks of the string
$12\cdots d$ turn out always to be maximizers, that is, $S|_n$ has
exactly $a_{nk}$ distinct subsequences of length $k$.  The generating
function for $\{ a_{nk} \}$ is given
by~\cite[equation~(7)]{FHS2004}:
$$
F(x, y) = \sum_{n,k} a_{nk} x^n y^k = \frac{1}{1 - x - xy(1 - x^d)}\, .
$$
This is of Riordan type with $\phi(x) = (1 - x)^{-1}$ and $v(x) = x +
x^2 + \dots + x^d$.

\begin{figure}[ht] \hspace{1.5in}
\includegraphics[scale=0.45]{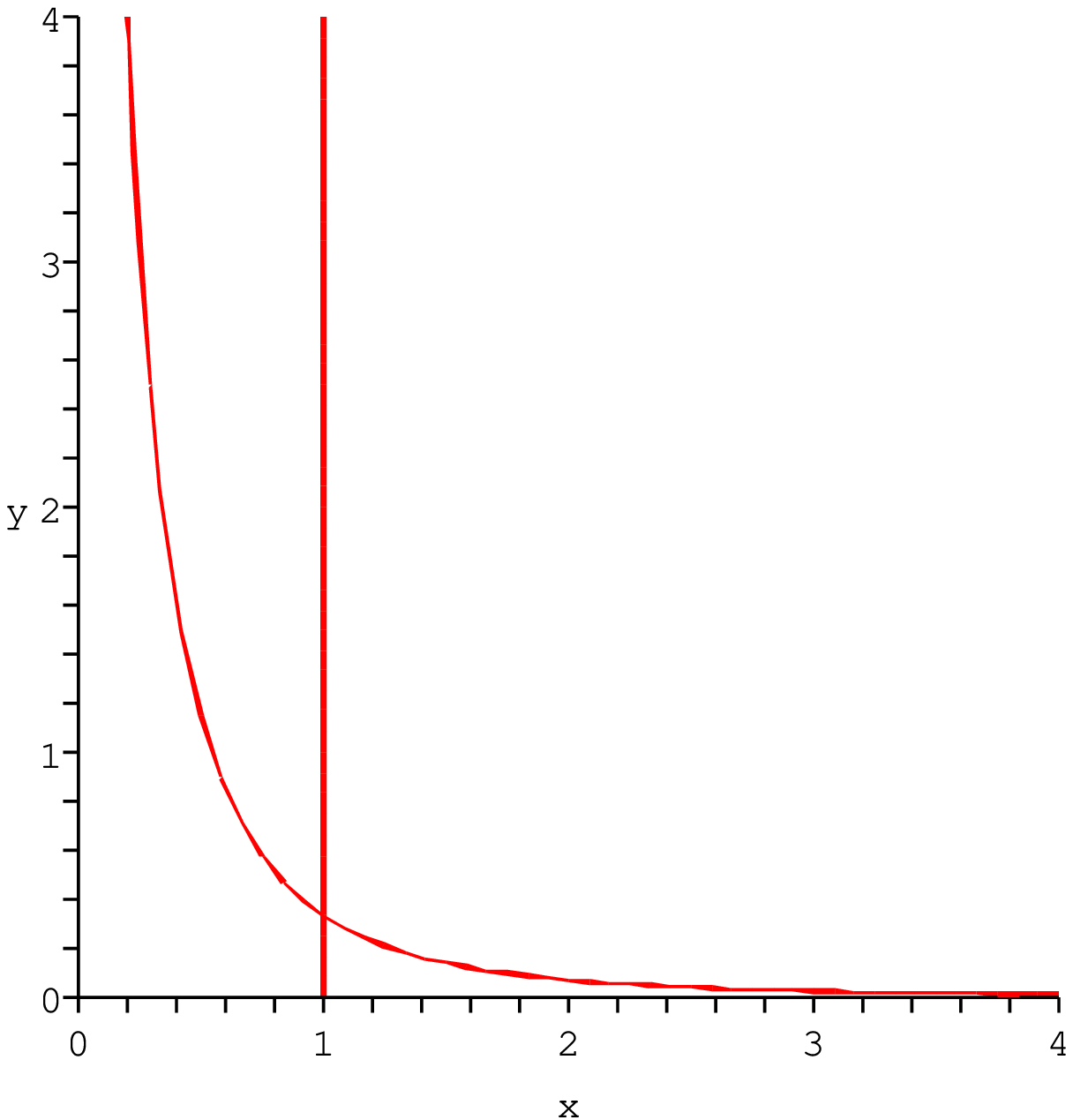}
\caption{$\sing$ in the case $d=3$}
\label{fig:flaxman} 
\end{figure}

The case $d = 1$ is uninteresting. Suppose that $d \geq 2$.  The
singular variety $\sing$ is the union of the line $x = 1$ and the
smooth curve $y = 1/v(x)$ and they meet transversely at the double
point $(1, 1/d)$; see Figure~\ref{fig:flaxman} for an illustration of
this when $d = 3$.

This is a case where the radius of convergence of $\phi$ is less than
the radius of convergence of $v$, the former being~1 and the latter
being infinite.  We have $\mu(v; x) = 1/(1-x) - dx^d/(1 - x^d) = (1 +
2x + 3x^2 + \dots + dx^{d-1})/(1 + x + x^2 + \dots + x^{d-1})$.  As
$x$ increases from~0 to~1 (the radius of convergence of $\phi$, which
is the value of $x$ at the double point), $\mu$ increases from~1 to
$(d + 1)/2$. Thus when $\lambda:=n/k$ remains in a compact
sub-interval of $(1,\frac{d + 1}{2})$, the Gaussian asymptotics of
equation~(\ref{eq:riordan-asymp}) hold.  

To compute these in terms of $\lambda$, we solve for $x$ in 
\begin{equation} \label{eq:poly 1}
\mu (v;x) = \lambda:= \frac{n}{k}
\end{equation}
and plug this into~(\ref{eq:riordan-lt}).  One can do this 
numerically, but in the case where $v$ is a polynomial, one
can do better.  

Solving $\mu(v;x) = \lambda$ by radicals and plugging
into~(\ref{eq:sigma}), which worked in
Section~\ref{ss:delannoy}, will not be possible when $d \geq 5$ and is
not practical even for smaller $d$.  However, we see that $\sigma^2$
is algebraic, in the same degree $d - 1$ extension of the rationals
that contains the value of $x$ solving~(\ref{eq:poly 1}).  We look
therefore for a polynomial with coefficients in $\mathbb{Q}(\lambda)$,
of degree $d - 1$, which annihilates the $\sigma^2$ in
Proposition~\ref{pr:riordan}.  To find this polynomial, the best
tactic is to work directly with generators of polynomial ideals, and
we give the details below.

When $d = 2$, the solution
$$
x (\lambda) = \frac{\lambda - 1}{2 - \lambda}
$$ 
is a rational function of $\lambda$ and nothing fancy is needed to
arrive at $\sigma^2 = (\lambda - 1)(2 - \lambda)$.  We therefore
illustrate with $d = 3$, though this procedure is completely general
and will work any time $v$ is a polynomial.

Plugging in the expression~(\ref{eq:mu}) for $\mu(v;x)$
in~(\ref{eq:poly 1}) and clearing denominators gives a polynomial
equation for $x$:
$$
x \frac{dv}{dx} - \lambda v = 0 \, .
$$
In our example,
\begin{equation} 
\label{eq:flaxman x}
x \left( 1+2\,x+3\,{x}^{2} \right) - \lambda \, \left(
x+{x}^{2}+{x}^{3} \right) = 0 \, .
\end{equation}
We now need to evaluate 
\begin{equation} 
\label{eq:x poly}
\sigma^2(v;x) = x \frac{d\mu}{dx} = {\frac {x \left( 1+4\,x+{x}^{2} 
   \right) } { \left( 1+x+{x}^{2} \right) ^{2}}}
\end{equation}
at the value $x$ that solves~(\ref{eq:flaxman x}).  

To do this we compute a Gr\"{o}bner basis of the ideal in
$\mathbb{Q}(\lambda)[x, S]$ generated by $\mu(v;x) - \lambda$ and
$\sigma^2(v;x) - S$ (after clearing denominators). The commands
\begin{quote}
\begin{verbatim}
p1:=(1+2*x+3*x^2)-lambda*(1+x+x^2):
p2:=x*(1+4*x+x^2)-S*(1+x+x^2)^2:
Basis([p1, p2], plex(x, S));
\end{verbatim}
\end{quote}
produce the elimination polynomial
$$
p(S;\lambda) = 3S^2 + (6\lambda^2 - 24\lambda + 16) S + 3\lambda^4 -
24\lambda^3 + 65\lambda^2 - 68\lambda + 24
$$
which is easily checked to be irreducible (using Maple's
\texttt{factor} command, for example), and hence is generically the
minimal polynomial for $\sigma^2$. 
It is easy to choose the right branch of the curve: the function
$x(\lambda)$ increases in $(0,1)$ as $\lambda$ increases in $(1,2)$, and
$\sigma^2$ is given in~\eqref{eq:x poly} as an explicit function of $x$
that is easily checked to be increasing. It follows from these that
$\sigma^2$ increases from $0$ to $2/3$ as $\lambda$ goes from $1$ to
$2$.

To finish describing the asymptotics, we first note that values of
$\lambda$ greater than $d$ are uninteresting.  It is obvious that any
prefix of $S$ of length at least $dk$ will allow all possible
$k$-subsequences to occur. Thus $a_{nk} = d^k$ when $\lambda \geq d$.

We already know that as $\lambda := n/k \to (d+1)/2$ from below, the
asymptotics are Gaussian and the exponential growth rate approaches
$d$.  For slopes $\lambda \geq (d+1)/2$, we use
Theorem~\ref{th:universal} to see that for each such $\lambda$, there
is a minimal point in the positive quadrant controlling asymptotics in
direction $\lambda$.  The only minimal point of $\sing$ we have not
yet used is the double point $(1, 1/d)$.  It is readily computed that
this cone has extreme rays corresponding to $\lambda = (d + 1)/2$ and
$\lambda = \infty$, and thus asymptotics in the interior of the cone
will be supplied by the double point. Using Corollary~\ref{cor:m=d=2}
we obtain $a_{\lambda k, k} \sim d^{k}$.

Although Section~\ref{ss:results} did not discuss what happens
when $\lambda$ approaches but is not equal to $(d+1)/2$, 
some results are known.  A refinement of
Corollary~\ref{cor:m=d=2}, given in~\cite[Theorem~3.1~(ii)]{PeWi2004}, 
shows that asymptotics in the boundary direction $\lambda = (d+1)/2$ 
are smaller by a factor of $2$ than in the interior of the cone.
In fact an examination of the proof there
(see~\cite[Lemma~4.7~(ii)]{PeWi2004}) shows that one has Gaussian behavior, 
$$
a_{nk} \sim d^k \Phi (x) \mbox{ for } n = \frac{d+1}{2} k 
+ x \sqrt{\frac{d^2-1}{12} \, k}
$$
where $\Phi (x) = (2 \pi)^{-1/2} \int_{-\infty}^x e^{-t^2/2} \, dt$ is
the standard normal CDF and the constant $\sqrt{(d^2-1)/12}$ is
obtained in a manner similar to $Q$.  This can also be obtained from a
probabilistic analysis as follows.  The quantity $d^{-k} a_{nk}$ is
the probability that a uniformly chosen sequence of length $k$ is a
subsequence of $S|_n$. The length of an initial substring of $S$
required to contain a given sequence $u_1 u_2 \cdots u_k$ is
$\sum_{j=1}^k (u_j - u_{j-1}) \; {\rm mod} \; d$ where $0 \; {\rm mod}
\; d$ is taken to be $d$ and where $u_0 := 0$.  Thus the probability
of a uniformly chosen word of length $k$ being a subsequence of $S|_n$
is equal to
$$
\P \left ( \sum_{j=1}^k U_j \leq n \right )
$$
where $U_j$ are independent uniform picks from $\{ 1 , \ldots , d \}$.
The central limit theorem now gives
$$
\P \left ( \sum_{j=1}^k U_j \leq \frac{d+1}{2} k 
+ x \sqrt{\frac{d^2-1}{12} \, k} \right ) \sim \Phi (x) \, .
$$

\subsection{Paths, hills and Fine numbers} 
\label{ss:fine}
\Em{Dyck paths} are paths from the origin whose steps are in the set
$\{ (1,1), (1,0), (1, -1) \}$ and that never go below the $x$-axis.
In Section~\ref{ss:paths} we will count Dyck paths by their final
point, which may be anywhere in the bottom half of the first quadrant.
Often in the literature, the term Dyck path is reserved for a path
constrained to end on the $x$-axis.  In~\cite{DeSh2001}, a
number of combinatorial interpretations are found for these
constrained Dyck paths when counted by final $x$-value and by the
number of \Em{hills}, a hill being a peak of height 1 (a peak is an occurrence of the step $(1, 1)$ immediately followed by  $(1, -1)$).
Let $a_{nk}$ denote the number of Dyck paths from the origin to $(2n,0)$ that have
$k$ hills.  Of particular interest are the values $a_{n0}$ which count
hill-free Dyck paths and are called \Em{Fine numbers}.

\begin{figure}[ht] \hspace{1.5in}
\includegraphics[scale=0.6]{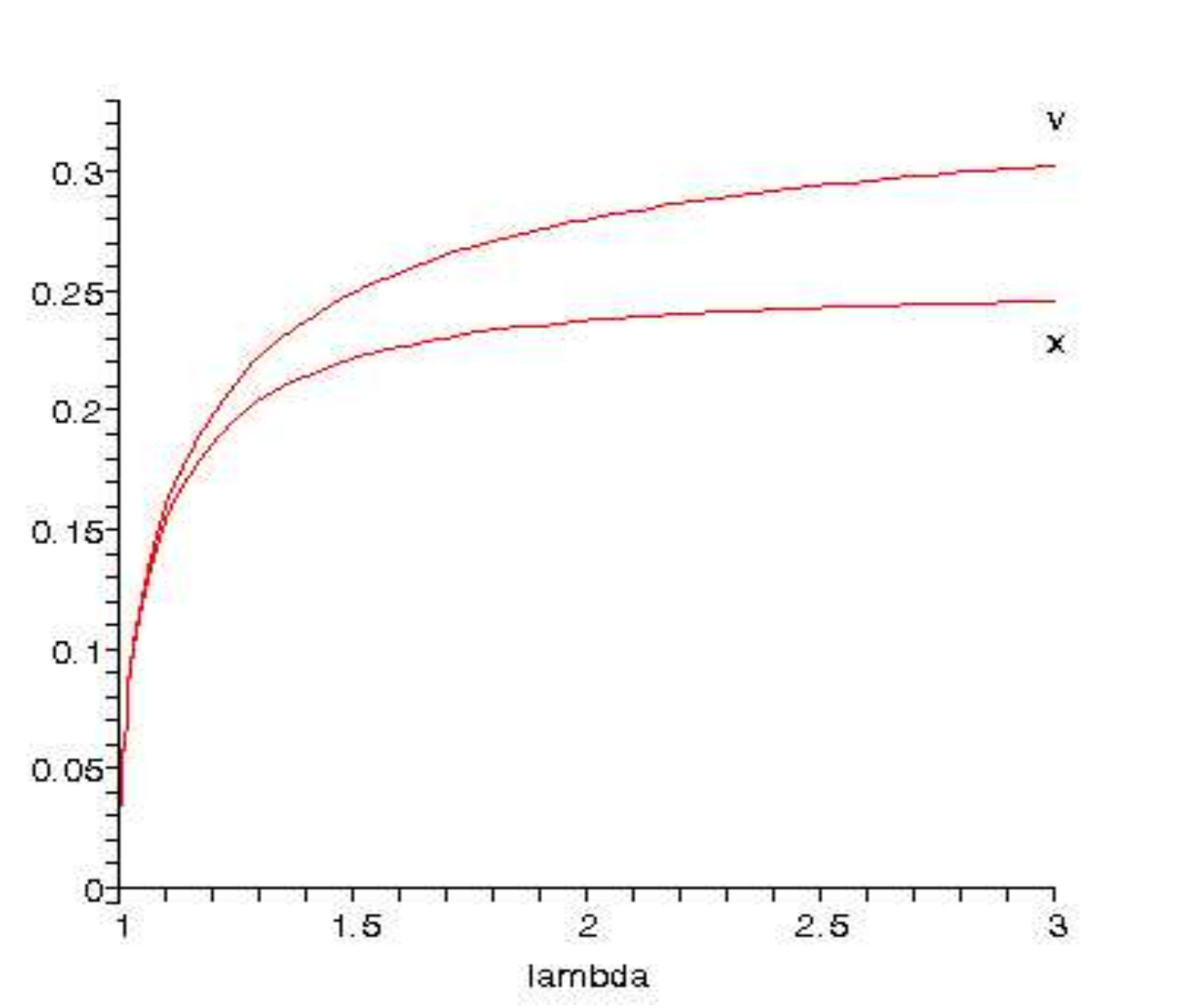}
\caption{$x$ and $v$ plotted against $\lambda$}
\label{fig:fine}
\end{figure}

The generating function for $\{ a_{nk} \}$ is derived 
in~\cite[Proposition~4]{DeSh2001}:
\begin{equation} \label{eq:fine}
F(x, y) := \sum_{n,k} a_{nk} x^n y^k 
   = \frac{2} {1 + 2x + \sqrt{1-4x} - 2xy }
   = \frac{v(x)/x}{1 - yv(x)} \, , 
\end{equation}
with
$$
v(x) = \frac{1 - \sqrt{1 - 4x}}{3 - \sqrt{1 - 4x}} ={\frac{2x}{1 +
\sqrt{1 - 4x} + 2x}}\, .
$$
Rationalizing the denominator, we have 
\begin{equation} 
\label{eq:v x}
v(x) = \frac{1 + 2x - \sqrt{1 - 4x}}{4 + 2x}  \, ,
\end{equation} 
and thus the domain of convergence is given by $|x| < 1/4, |yv(x)| <
1$. For $1 < \lambda < \infty$ the critical points are determined by 
$$ 
\mu(v;x) = \frac{4x}{(1 - \sqrt{1 - 4x}) (\sqrt{1 - 4x}) 
(3 - \sqrt{1 - 4x})} = \lambda \, . 
$$
The left side of the above equality increases strictly from 1 at $x =
0$ to $\infty$ as $x \uparrow 1/4$, and there is a unique positive
real solution $x_\lambda$ for each $\lambda > 1$.  Figure~\ref{fig:fine}
shows $x$ increasing from~0 to $1/4$ and $v$ increasing from~0 to $1/3$
as $\lambda$ increases from~1 to~$\infty$.  

We now want to complete the computations by simplifying
$\sigma^2(v;x)$ at the value of $x$ satisfying $\mu(v;x) = \lambda$.
We sketch here how the computations of Section~\ref{ss:distinct
subseq} generalize to the algebraic case. The algebra is only a little
more involved than it is for rational $v$ but some changes are clearly
required.

In order to perform Gr\"{o}bner basis computations, we must first pass
from~\eqref{eq:v x} to the implicit form $\alpha(x, v) = 0$ where 
$\alpha$ is a polynomial in $\Q [x,v]$.  It is easy in this case
to see by inspection that 
$$
\alpha (v , x) = (2+x) v^2 - (1+2x) v + x \, ;
$$
at worst, implicitizing will require going back to the
original derivation and following the computation at the
level of ideals \footnote{Implicitization techniques for polynomial 
or rational parametrizations $x = f(t) , v = g(t)$, though not needed
here, may be found in~\cite[Sections 3.2 and 6.4]{CLO2005}.}.
Next, according to~(\ref{eq:mu}) and~(\ref{eq:sigma}), differentiating 
$\alpha$ implicitly gives $\mu$ and $\sigma^2$ as rational functions 
of $v$ and $x$:  
\begin{equation} \label{eq:S fine}
\begin{array}{rclcl}
v' & = & \disp{- \frac{\alpha_x}{\alpha_v}} 
   & = & \disp{\frac{v^2 - 4x}{2v(2+x)} } \\[2ex]
v'' & = & \disp{ - \frac{\alpha_v^2 \alpha_{xx} + \alpha_x^2 \alpha_{vv} 
   - 2 \alpha_{xv} \alpha_x \alpha_v}{\alpha_v^3}} & = & 
   \disp{\frac{3v^4+16v^2-16x^2}{4v^3(2+x)^2}} \\[2ex]
\mu & = &  \disp{- \frac{x \alpha_x}{v \alpha_v}} & = & 
   \disp{\frac{x(4x-v^2)}{2v^2(2+x)}} \\[2ex] 
\sigma^2 & = &  \disp{- \frac {x \left( xv \alpha_v^2 \alpha_{xx}
   +xv \alpha_x^2 \alpha_{vv} - 2xv \alpha_{xv} \alpha_x \alpha_v 
   + v \alpha_v^2 \alpha_x + x \alpha_x^2 \alpha_v \right) }
   {\alpha_v^3 v^2} }
&=& \disp{\frac{-x(8x^3 - 4x^2 v^2 - 8xv^2 + v^4)}{v^4 (2+x)^2}} \, .
\end{array}
\end{equation}
Setting $\mu = n/k$ and clearing denominators gives a polynomial
$\beta = n v \alpha_v + k x \alpha_x \in \Q [n,k] \, [x,v]$ that
vanishes when $\mu(v;x) = n/k$. In our present example, we have
$$
\beta(x, v) = \left( 4\,n+2\,nx+kx \right) {v}^{2}+ \left(
   -2\,kx-n-2\,nx \right) v +kx \, . 
$$
 
To find the critical point corresponding to the direction $n/k$, we
could simply solve the equations $\alpha = 0, \beta = 0$ for $x, v$.
Then to simplify $\sigma^2$, we add in the equation stating that
$\sigma^2 = S$ (with denominator cleared), and find a Gr\"{o}bner
basis using the ``plex" order with $S$ as the last variable, just as
in the previous section. Unfortunately, as is well known, such
computations can be very slow, and in the present instance we have
difficulty in obtaining an answer with Maple in a reasonable time.

Thus we use the following alternative method, as described in
\cite[Section~2.2]{CLO2005}. The polynomials $\alpha$ and $\beta$ 
have finitely many
common solutions $(x,v)$ for any fixed $n,k$, so they define a
zero-dimensional ideal, $J$.  In other words, there are finitely many
linearly independent monomials in $x$ and $v$ over $\Q [n,k] \, / \,
\langle \alpha,\beta \rangle$.  A convenient choice is the set of
monomials not divisible by any leading term of a fixed Gr\"obner basis
for $J$. With respect to this basis, one may express $x$ and $v$ as
matrices over $\Q [n,k]$ for multiplication operators, then compute
the matrix for $\sigma^2$ as a rational function of these matrices,
and finally, find the minimal polynomial $\gamma$ for this matrix. The
polynomial $\gamma \in \Q [n,k]$ will vanish at $\sigma^2$.

To carry this out in Maple, we use the following commands:
\begin{verbatim}
sys:=[alpha, beta]:
monomialorder := tdeg(x, v):
gb := Basis(sys, monomialorder):
ns, rv := NormalSet(gb, monomialorder):
Mx := MultiplicationMatrix(y, ns, rv, gb, monomialorder):
Mv := MultiplicationMatrix(v, ns, rv, gb, monomialorder):
\end{verbatim}
Now evaluate the rational expression for $\sigma^2$ given above,
with $x, v$ replaced by $M_x, M_v$ (making  liberal use of the
simplification capabilities of Maple). We then compute the minimal
polynomial using the MinimalPolynomial command. The result is 
\begin{equation*}
\begin{split}
& z^3 \left (32 n^4 k^5 \right ) \\
+ & z^2 \left( -144 n^7 k^2 + 160 n^6 k^3 - 16 n^5 k^4 - 6 n^4 k^5 
   + 4n^2 k^7 + 2 k^9 \right ) \\
+ & z \left( 144 n^9 - 304 n^8 k + 203 n^7 k^2 - 46 n^6 k^3 - 15 n^5
   k^4 + 20 n^4 k^5 - 11 k^6 n^3 + 10 n^2 k^7 - n k^8 \right) \\
+ & \;\; \left ( - 15n^7k^2 + 31n^5k^4 - 27n^9 - 35n^6k^3 - 21n^4k^5 
   + 57n^8k + 11k^6n^3 - n^2k^7 \right )
\end{split}
\end{equation*}
and computing eigenvalues (using \texttt{solve}) we obtain three
possible values for $\sigma^2$, namely 
\begin{eqnarray*}
S_1 & := & \frac {(n - k)(n + k)(3n^2 + k^2)}{16n^4} \\[1ex]
S_2 & := & \frac {\left(9n - k + 3\sqrt{9n^2 - 10nk + k^2} \right) 
   \left(n - k \right) n}{4k^3} \\[1ex] 
S_3 & := & \frac {\left(9n - k - 3\sqrt{9n^2 - 10nk + k^2} \right) 
   \left(n - k \right) n}{4k^3} \, .
\end{eqnarray*}
Seeing which of these three gives the correct value for $\sigma^2$
may be done similarly to the way it was in the previous section.
Equation~\eqref{eq:S fine} gives $\sigma^2$ as a univalent function
of $x$ and $v$.  The correct branches for $x$ and $v$ as functions 
of $\lambda = n/k$ is shown in Figure~\ref{fig:fine}.  Plugging
these into~\eqref{eq:S fine} and comparing to $S_1, S_2$ and $S_3$
(to avoid a numerical comparison, one may compare limits at 0) shows
that second expression is correct:
$$
\sigma^2 = S_3 := \frac {\left(9n - k + 3\sqrt{9n^2 - 10nk + k^2}
\right) \left(n - k \right) n}{4k^3} \, .
$$   

To write $x (\lambda)$ and $v(\lambda)$ in their simplest forms,
observe that the values of $x$ and $v$ are simply eigenvalues of the
multiplication matrices $M_x, M_v$ and so we obtain their respective
minimal polynomials: 
\begin{align*}
4n^2 x^2 + (7n - k)(n - k) x - 2n(n - k) \, , \\
\mbox{ and } \hspace{0.30in} 2k v^2 +  3(n - k) v - (n - k) \, . 
\end{align*}
Since we know that the relevant point has positive coordinates, we have
the explicit forms
\begin{align*}
x & = \frac {(3n - k) \sqrt{(n - k)(9n - k)} - (7n - k)(n - k)}{8n^2}
\\
v & = \frac{\sqrt{(n - k)(9n - k)} - 3(n - k)}{4k} \, .
\end{align*}

As can be seen from this example, the computations can be rather
messy, though still routine and automatable.  In the end we have
verified that $\sigma^2 > 0$ for $\lambda > 1$, and we conclude
that Gaussian asymptotics given by Proposition~\ref{pr:riordan}
hold uniformly as $n/k$ varies over any compact subinterval
of $(1,\infty)$.  

Another way to simplify computations such as this one will be discussed
in Section~\ref{sec:lagrange}. Of course, for any given value of $n, k$,
we may shortcut the above computation and solve numerically for $x$ to
obtain $\sigma^2$.  For example, with $n/k = 2$, we have $x =
(5\sqrt{17} - 13)/32 \approx 0.2379852541$, so that $y = 1/v(x) = (3 +
\sqrt{17})/2 \approx 3.561552813$.  Hence 
$$ 
a_{2k, k} \sim (0.1228255460\cdots) (4.957474791\cdots)^k k^{-1/2} \,
.
$$
As far as we know, asymptotics such as these have not previously been
computed.  Using 10 significant figure floating point approximations, 
we obtain (using Maple) an answer accurate to within $0.8\%$ already 
for $k = 30$.

\subsection{Horizontally convex polyominoes}
\label{ss:polyomino}

A \Em{horizontally convex polyomino} (HCP) is a union of cells
$[a,a+1] \times [b,b+1]$ in the two-dimensional integer lattice 
such that the interior of the figure is connected and every row
is connected.  Formally, if $S \subseteq \Z^2$ and 
$P = \bigcup_{(a,b) \in S} [a,a+1] \times [b,b+1]$ then $P$ is
an HCP if and only if the following three conditions hold:
$B := \{ b : \exists a, (a,b) \in S \}$ is an interval; the
set $A_b := \{ a : (a,b) \in S \}$ is an interval for
each $b \in \Z$; and whenever $b , b+1 \in B$, the sets $A_b$ 
and $A_{b+1}$ intersect.
\begin{figure}[ht] \hspace{2.1in}
\includegraphics[scale=0.6]{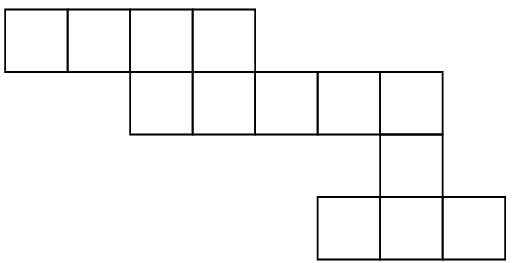}
\caption{an HCP with 13 cells and 4 rows}
\label{fig:hcp2} 
\end{figure}
Let $a_n$ be the number of HCP's with $n$ cells, counting two
as the same if they are translates of one another.  P\'olya~\cite{Poly1969} 
proved that 
\begin{equation} \label{eq:HCP by n}
\sum_n a_n x^n = \frac{x (1-x)^3}{1 - 5x + 7 x^2 - 4 x^3} \, .
\end{equation}
Further discussion of the origins of this formula and its accompanying
recursion may be found in~\cite{Odly1995}
and~\cite{Stan1997}.
The proof in~\cite[pages~150--153]{Wilf1994} shows in fact that
\begin{equation} \label{eq:HCP by n,k}
F(x,y) = \sum_{n,k} a_{nk} x^n y^k = \frac{xy (1 - x)^3}{(1 - x)^4 
   - xy (1 - x - x^2 + x^3 + x^2 y)} \, ,
\end{equation}
where $a_{nk}$ is the number of HCP's with $n$ cells and $k$ rows.
Let us find an asymptotic formula for $a_{rs}$.

All the coefficients of $F(x,y)$ are nonnegative; they vanish when $s >
r$ but otherwise are at least~1.  By Corollary~\ref{cor:procedure} (the
last part, which requires only $a_\rr \geq 0$), we know that all points
of $\sing$ in the first quadrant that are on the southwest facing part
of the graph (that is, that are minimal in the coordinatewise partial
order) are contributing critical points. We do not know yet but will see
later that there are no other critical points on each torus.  
As $\rbar$ varies over $\dirset = \{ \rbar : 0 < s/r < 1 \}$
from the horizontal to the diagonal, the point $\contrib_{\rbar}$ 
moves along this graph from $(1, 0)$ to $(0, \infty)$.

\begin{figure}[ht] \hspace{1.5in}
\includegraphics[scale=0.4]{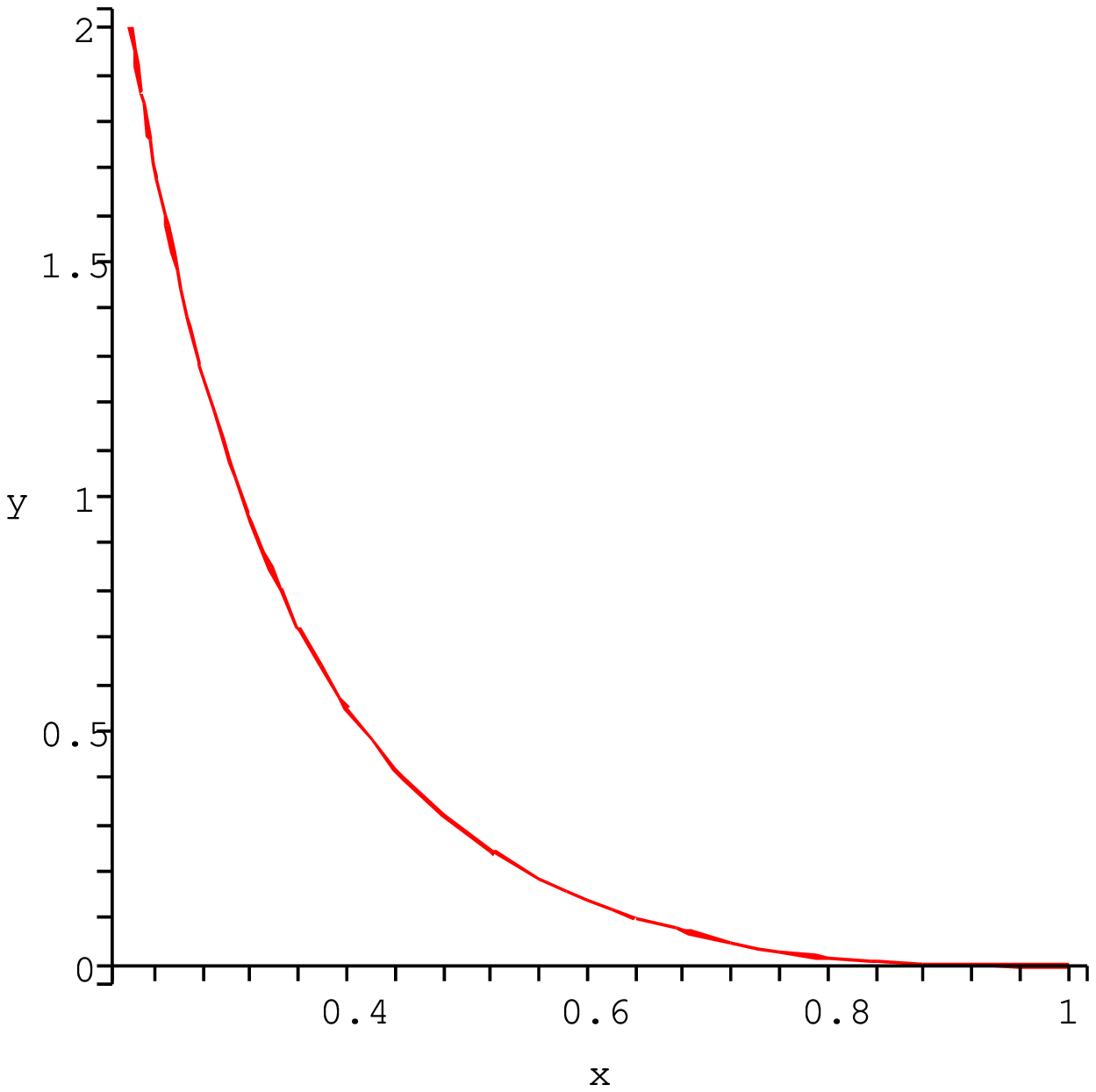}
\caption{minimal points of $\sing$ in the positive real quadrant}
\label{fig:hcp}
\end{figure}

To make the mapping from $\rbar$ to $\zz$ explicit, we use the fact that
$\rbar \in \cone (\zz) \subseteq \linear (\zz)$
(Theorem~\ref{th:universal} and part~(i) of Proposition~\ref{pr:cone}),
so $\zz$ may be gotten from $\rr$ by~(\ref{eq:zero-dim}). It is readily
computed that $\grad \denom \neq \zero$ except at $(1,0)$. Thus all
minimal points are smooth.  The only solution to $\numer = \denom = 0$
is at $(1,0)$, so the numerator is nonvanishing at any minimal point as
well.  Checking whether the quantity $Q$ defined in
equation~(\ref{eq:Q}) of Theorem~\ref{th:smooth d=2} ever vanishes, we
find that the only solutions to $\denom = Q = 0$ are at $(1,0)$ and at
complex locations that are not minimal because $\linear (\zz)$ is not
real there (Proposition~\ref{pr:linear}). Lastly, we must check that
there are no contributing points other than the minimal points in the
positive quadrant. We can ascertain that there are not by checking the
``extraneous'' critical points, which will lie on three other branches
of a quartic (see below), and seeing that they lie on the wrong torus.
We then know that the asymptotics for $a_{rs}$ are uniform as $s/r$
varies over a compact subset of the interval $(0,1)$ and given by
$$
a_{rs} \sim C x^{-r} y^{-s} r^{-1/2} \, .
$$

We will use Maple to determine $x,y$ and $C$ as explicit functions 
of $\lambda := s/r$, giving asymptotics for the number of HCP's whose 
shape is not asymptotically vertical or horizontal, but first
we see what we may tell without much computation.

A crude approximation at the logarithmic level is
$$
a_{rs} \approx \exp \left [ r (- \log x - (s/r) \log y \right ]
$$
where $x$ and $y$ of course still depend on $s/r$.  
We first compute the average row length of a typical HCP.  We may
apply Theorem~\ref{th:WLLN} to find the limit in probability
of $h_k / k$, where $h_k$ is the height of an HCP chosen uniformly
from among all HCP's of size $k$.  Setting $y = 1$ in the bivariate 
generating function recovers the univariate generating
function~(\ref{eq:HCP by n}).  The point $(x,1)$, where $x = x_0$ is
the smallest root of the denominator of~(\ref{eq:HCP by n}) controls
asymptotics in this direction; we compute $\dir (x,1)$ there to be
$(x \denom_x , y \denom_y) (x_0 , 1)$, which simplifies to $r/s = 
4 (5 - 14x_0 + 12x_0^2)/(5 - 9x_0 + 11x_0^2)$ or still further to 
$\alpha:=\frac{1}{47} (147 - 246x_0 + 344x_0^2) \approx 2.207$.  
We conclude that for large $k$, $k / h_k \to \alpha$.
Thus we see that the average row length in a typical large HCP 
is around $2.2$.
Finally, let us see how the computer algebra for general $\rbar$ turns
out.  To find $(x,y)$ given $\rbar$, solve the pair of equations
$\denom = 0, s x \denom_x = r y \denom_y$; this find points on $\sing$
with $\rbar = \overline{(r,s)} \in \linear (x,y)$. Explicitly, we set
$\lambda = s/r$ and ask Maple for a Gr\"obner basis for the ideal
generated by $\denom$ and $\lambda x \denom_x - y \denom_y$. 
The following Maple code fragment is useful in such cases. 

\begin{verbatim}
Hx := diff(H,x): Hy := diff(H,y): X:=x*Hx: Y:=y*Hy:
Hxx := diff(Hx,x): Hxy := diff(Hx,y): Hyy := diff(Hy,y):
Q := -X^2*Y-X*Y^2-x^2*Y^2*Hxx-X^2*y^2*Hyy+2*X*Y*x*y*Hxy:
L := [H,lambda*X-Y]:
gb := Basis(L, plex(y,x)):
\end{verbatim}

Maple returns a basis consisting of polynomials $\alpha$ 
and $(x - 1)^5 \beta$, where the quartic
$$
\beta := (1 + \lambda) x^4 + 4(1 + \lambda)^2 x^3+ 10 (\lambda^2 +
\lambda - 1) x^2 + 4 (2k - 1)^2 x + (1 - \lambda) (1 - 2 \lambda)
$$
is the elimination polynomial for $x$ and is generically irreducible.
Furthermore $\alpha$ is linear in $y$. Rather than express $y$ in terms
of $x$, it is perhaps easier to compute the elimination polynomial for
$y$, which is $y^3$ times the following polynomial:
\begin{eqnarray*}
&& (4 \lambda^4 - 4 \lambda^3 - 3 \lambda^2 + 4 \lambda - 1) y^4 
   + (40 \lambda^4 - 44 \lambda^3 - 20 \lambda^2 +48 \lambda - 16) y^3 \\
& + & (-172 \lambda^4 + 128 \lambda^3 + 160 \lambda^2 -256 \lambda + 64) y^2
   + (1152 \lambda^4 - 1024 \lambda^3 - 512 \lambda^2) y -1024 \lambda^4 \, .
\end{eqnarray*}
Note that when $H=0$, $x = 1$ if and only if $y = 0$. The point $(1,0)$ is a
solution to the critical point equations for every value of $\lambda$.
However it is never a contributing point, because $h_{\rhat} = \infty$ 
when any coordinate vanishes (recall $\rhat$ is in a compact subset
of the positive orthant).

Thus generically we have $4$ candidates for contributing points.
Precisely one of these is minimal and in the first quadrant. The 
others do not contribute for generic $\rbar$, which may easily
be checked for any given $\rbar$; a verification for all $\rbar$ 
simultaneously would require more computer algebra.

Finally, from Corollary~\ref{th:smooth d=2} we see that
$$
C = \frac{x y (1-x)^3}{\sqrt{2 \pi}} \sqrt{\frac{y \left( -x \left(
1 - x - x^2 + x^3 + x^2 y \right) - x^3 y \right)} {Q}} \, .
$$ 
The minimal polynomial for $\sqrt{2\pi} C$ can be computed as in
previous sections; it turns out to have degree $8$ for generic
$\lambda$. Of course, given floating point approximations for $x$ and
$y$, we may simply compute an approximation for $C$ directly.

As an example, suppose that $n = 2k$ so that $\lambda = 1/2$. In this 
case we have simplification and the minimal polynomials for $x$ and $y$ 
respectively are $3 x^2 + 18x - 5$ and $75 y^2 - 288 y + 256$.  Note that 
there is a single element 
$$(x_0 , y_0) := \left ( \sqrt{\frac{32}{3}} - 3 , \frac{48 - \sqrt{512}}{25}
   \right ) \approx (0.265986, 1.397442)$$
of $\crit$ in the positive 
quadrant.  By Theorem~\ref{th:universal} this belongs to 
$\contrib$, so by Proposition~\ref{pr:minimal}, $\contrib$ 
is the singleton $\{ (x_0 , y_0) \}$.
We obtain (by naive floating point computation in Maple) 
$a_{n, n/2} \sim (0.237305\dots) (3.18034\dots)^n n^{-1/2}$. 
For $n = 60$, the relative error in this first order approximation 
is about $1.5\%$.  Contrast the exponential growth rate with the 
exponential growth rate $3.20557$ for the number of HCPs with 
$n$ cells and any number of rows.

\subsection{Symmetric Eulerian numbers}
\label{ss:euler}

\begin{figure}[ht] 
\hspace{0.4in}
\includegraphics[scale=0.4]{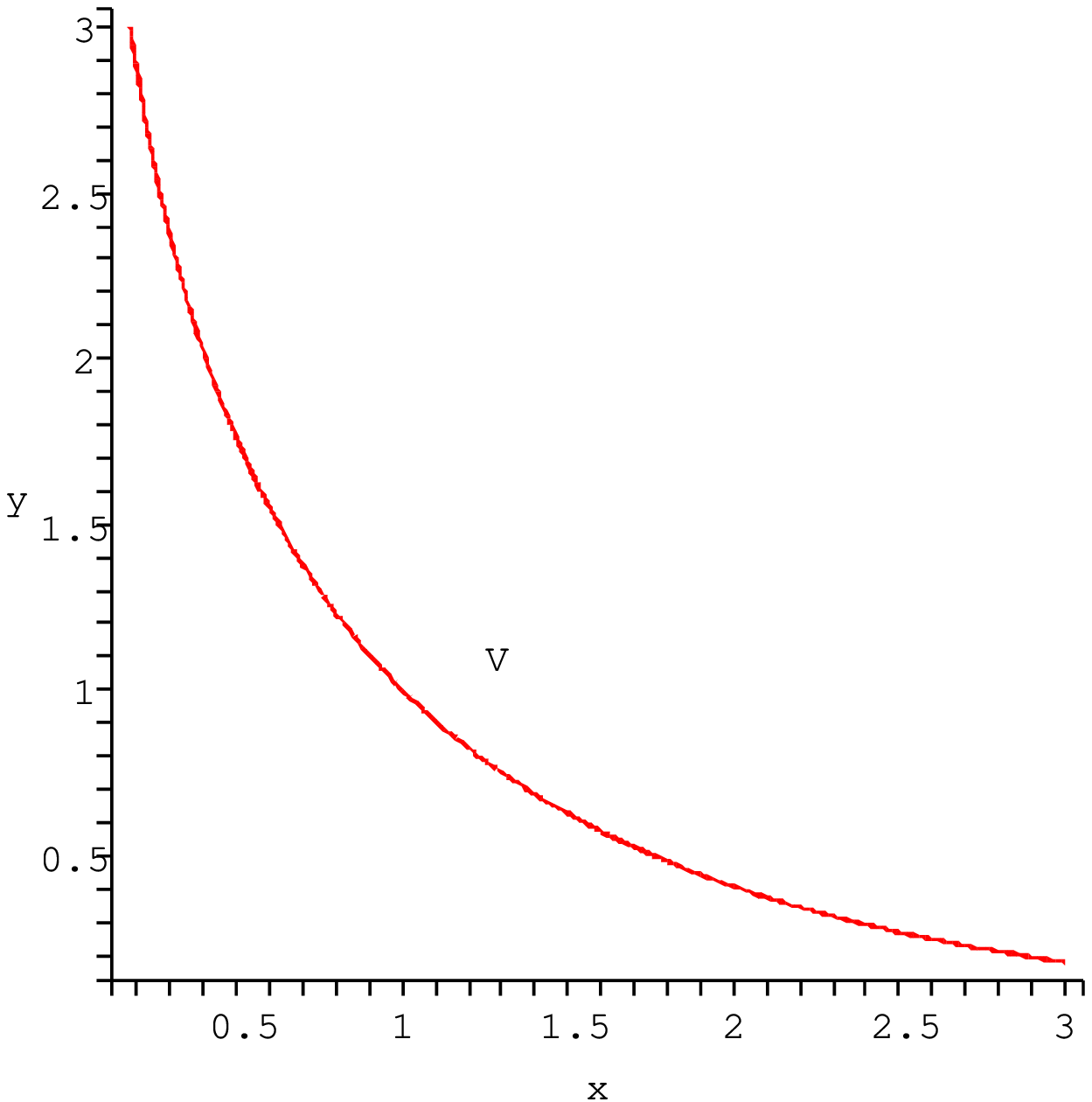}
\hspace{0.8in}
\includegraphics[scale=0.4]{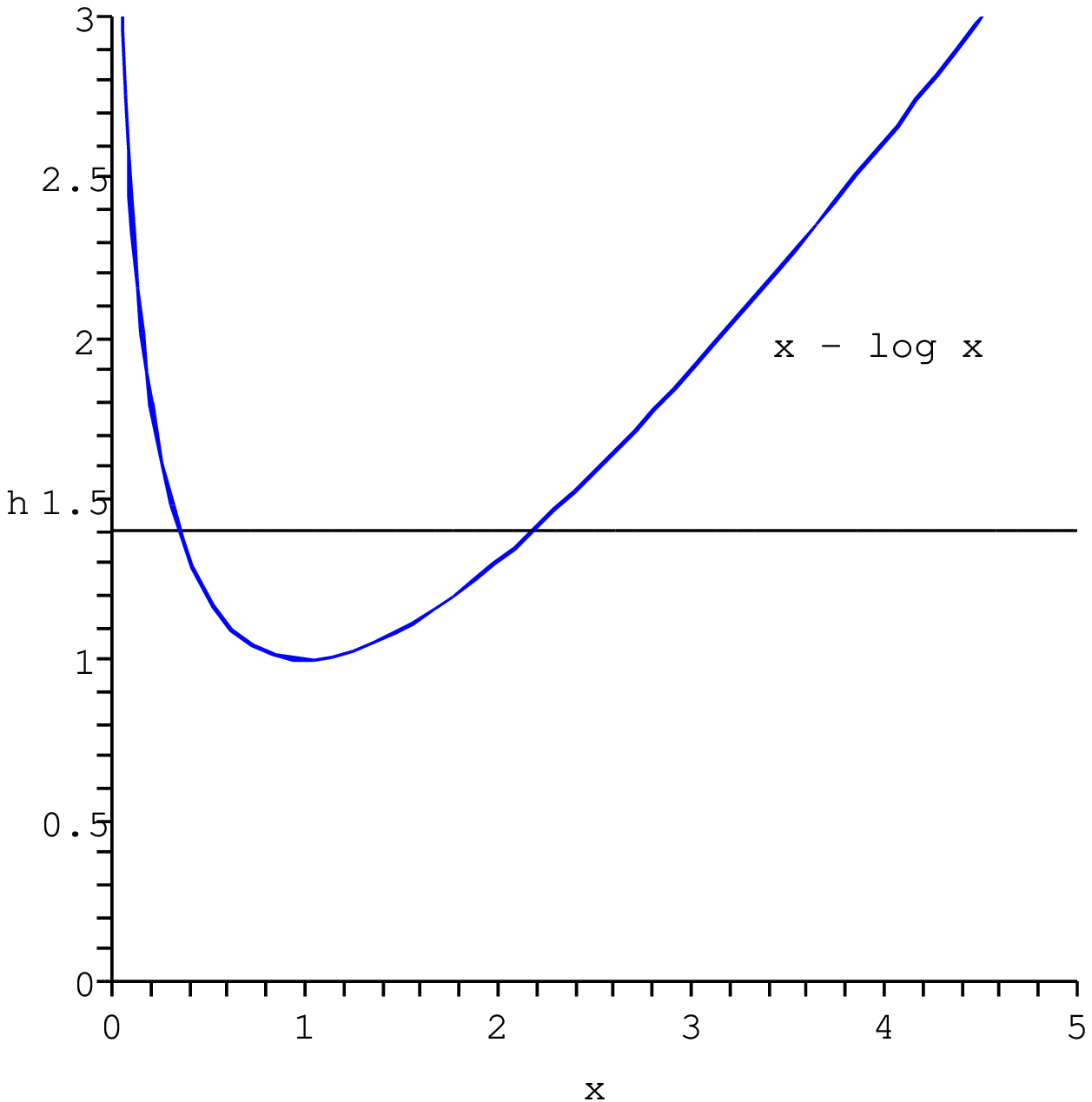}
\caption{the positive real points of $\sing$ and their
description by height-pairs}
\label{fig:euler}
\end{figure}

The symmetric Eulerian numbers
$\eul(r,s)$~\cite[page~246]{Comt1974} count the number of
permutations of the set $[r + s + 1]:=\{1, 2, \dots, r + s + 1\}$ with
precisely $r$ descents. By reading backwards, this equals the number
with exactly $s$ descents, hence the symmetry in $r$ and $s$.  The
symmetric Eulerian numbers have exponential generating
function~\cite[2.4.21]{GoJa2004} 
\begin{equation} 
\label{eq:eulerian}
F(x,y) = \frac{e^x - e^y}{x e^y - y e^x} = \sum_{r,s} \frac{\eul
(r,s)} {r! \, s!} =: \sum_{r,s} a_{rs} x^r y^s \, .
\end{equation}
The denominator in this representation is singular at the origin and
has a factor of $(x-y)$ in common with the numerator.  We factor this
out of the numerator and denominator, writing $F = \numer / \denom$
with $\numer = (e^x - e^y) / (x-y)$ and $\denom = (x e^y - y e^x) / (x
- y)$. In this representation, both the numerator and the denominator
are entire functions.  We quickly check that in the positive real 
quadrant, $\sing$ is the graph of a monotone deceasing function 
as in the left-hand side of Figure~\ref{fig:euler} and that the
quantity $Q$ (see the Maple code fragment in Section~\ref{ss:polyomino})
does not vanish on $\sing$:
\begin{itemize}
\item[] first check that $F$ is not entire but is analytic in a
neighborhood of the origin;  
\item[] next, L'H\^opital's rule shows us that $\denom$ never vanishes 
on $\{ x=y \}$ except at $(1,1)$;  
\item[] rearranging terms then shows that $\denom$ vanishes when
$x$ and $y$ are two positive real numbers with the 
same height $h(t) := t - \log t$ (see the right-hand side of 
Figure~\ref{fig:euler});
\item[] the exact symbolic expression for $Q$ has limit $e^3/12$
at $(1,1)$ and this is the minimum value of $Q$ on $\sing$.
\end{itemize}

The symmetric Eulerian numbers are nonnegative, so we can find
minimal points of $\sing$ in the first quadrant.  It is easy to 
check that the gradient of $\denom$ never vanishes on $\sing$, 
so $\sing$ is smooth.  At the point $x = y = 1$, we have 
$\partial \denom / \partial x = \partial \denom / \partial y = 1$, 
so $\dir(1,1) = \one$. For any other point of $\sing$, we may work 
with the representation of $\denom$ having $(x - y)$ in the denominator.  
We set
$$
\alpha = \frac{x \partial \denom / \partial x}
{x \partial \denom / \partial x + y \partial \denom / \partial y}
$$ 
so that $(x,y) \in \rbar \Leftrightarrow \rhat = (\alpha , 1-\alpha)$.
The expression for $\alpha$ does not look all that neat but on $\sing$
we may substitute $(y/x) e^x$ for $e^y$, after which the expression
for $\alpha$ reduces to $(1 - x)/(y - x)$.  We see therefore from 
Theorem~\ref{th:smooth asym} that ${\hat{A}} (r,s)$ is asymptotic 
to $C_\alpha (r+s)^{-1/2} r! s! \gamma^{r+s}$, where for a given value of
$\alpha = r/(r+s)$, the value of $\gamma$ is given by $x^{-\alpha}
y^{-(1 - \alpha)}$ after solving the following transcendental equations
for $(x,y)$:
\begin{eqnarray}
\label{eq:eulercrit}
\frac{1 - x}{y - x} & = & \alpha \, ; \\
x e^y & = & y e^x \, . \nonumber
\end{eqnarray}
We do not know whether there is a closed form expression for $\gamma$ in
terms of $\alpha$, though we note that by using~\eqref{eq:eulercrit}
several times, one may simplify somewhat: $\gamma = (y/x)^\alpha y^{-1}
= \exp(\alpha(y - x)) / y = \exp(1 - x)/y$.

The fact that the equations~\eqref{eq:eulercrit} have a unique positive 
real solution is possible to verify directly, but also follows from 
nonnegativity of the coefficients and the fact that the positive real 
part of $\sing$ lies along the boundary of the domain of convergence.
One must still check that $\torus (x(\alpha) , y(\alpha))$
contains no other points of $\sing$.  This is true for generic $\rbar$ 
and can easily be checked for a given $\rbar$, but a general proof that 
it works for all $\rbar$ is not obvious.

\subsection{Smirnov words}
\label{ss:smirnov}

Given an integer $d \geq 3$ we define a Smirnov word in the alphabet
$\{ 1 , \ldots , d \}$ to be a word in which no letter repeats
consecutively.  The number of Smirnov words of length $n$ is of course
easily seen by a direct counting argument to be $d (d-1)^{n-1}$. If we
count these words according to the number of occurrences of each
symbol, we get the multivariate generating function 
\begin{equation} 
\label{eq:smirnov GF}
F(\zz) = \sum a_\rr \zz^\rr = \frac{1}{1 - \sum_{j=1}^d 
 \frac{z_j}{1 + z_j}}
\end{equation}
where $a_\rr$ is the number of Smirnov words with $r_j$ occurrences of
the symbol $j$ for all $1 \leq j \leq d$.  This generating function
may be derived from the generating function $$E(\yy) = \frac{1}{1 -
\sum_{j=1}^d y_j}$$ for all words as follows. Note that collapsing all
consecutive occurrences of each symbol in an arbitrary word yields a
Smirnov word; this may be inverted by expanding each symbol of a
Smirnov word into an arbitrary positive number of identical symbols,
whence $\displaystyle{F(\frac{y_1}{1-y_1} , \ldots ,
\frac{y_d}{1-y_d}) = E(\yy)}$. Substituting $y_j = z_j / (1 + z_j)$
yields the formula~(\ref{eq:smirnov GF}) for $F$; see~\cite[II.7.3, Example 23]{FlSe}
 for more on this old result.  

When we write $F$ as a rational function, we see that
$$
\begin{array}{rclcl}
\numer & = & \prod_{j=1}^d (1+z_j) & = & \sum_{j=0}^d e_j (\zz)  \\
\denom & = & \sum_{j=0}^d (1-j) e_j (\zz) && 
\end{array}
$$
where $e_j$ is the $j^{th}$ elementary symmetric function of 
$z_1 , \ldots , z_d$ and $e_0 = 1$ by definition.  Thus when 
$d = 3$ for example, one obtains $\denom = 1 - (xy+yz+xz) - 2 xyz$.
The expression for $\denom$ may be denoted quite compactly by
$$
\denom = \left. g - u \frac{\partial g}{\partial u} \right |_{u=1}
$$
where $g (\zz , u) := \prod_{j=1}^d (1 + u z_j)$ is the polynomial in
$u$ whose coefficients are the elementary symmetric functions in
$\zz$. When $d \geq 3$, the denominator is of the form $1-P$ for an aperiodic
polynomial $P$ with nonnegative coefficients, so by
Proposition~\ref{pr:classify}, $\contrib_{\rbar}$ will always consist
of one strictly minimal point in the positive orthant.

The typical statistics of a Smirnov word of length $n$ are not in
doubt, since it is clear that each letter will appear with frequency
$1/d$.  This may be formally deduced from Theorem~\ref{th:WLLN2},
which also shows the variation around the mean to be Gaussian. We may,
perhaps, be more interested in the so-called large-deviation
probabilities: the exponential rate at which the number of words
decreases if we alter the statistics.  Let $\shat$ denote some
frequency vector.  Then as $|\rr| \to \infty$ with $\rhat \to \shat$,
we have $$\frac{1}{|\rr|} \log a_{\rr} \to -\rr \cdot \log |\zz|$$
where $\zz \in \mathcal{O}^d$ satisfies $\dir(\zz) = \rbar$.

To solve for $\zz$ in terms of $\rbar$, we use the symmetries of the
problem.  It is evident that $z_j$ is symmetric in the variables
$\{r_i : i \neq j \}$.  For example, when $d = 3$, if we refer to
$\zz$ as $(x,y,z)$ and $\rhat$ as $(r,s,t)$, then the
equations~(\ref{eq:zero-dim}) have the solution 
$$
x = \frac{r^2 - (s - t)^2}{2 r (s + t - r)}
$$ 
where the values of $y$ and $z$ are given by the same equation with
$r, s$ and $t$ permuted.  When $d \geq 4$, the solution is not a
rational function and it is difficult to get Maple to halt on a
Gr\"obner basis computation.  This points to the need for
computational algebraic tools better suited to working with symmetric
functions.

\subsection{Alignments of sequences}
\label{ss:alignments}

The problem of \Em{sequence alignment} is of interest in molecular
biology \cite{Wate1995, ReTa}, since a given string may
evolve via substitutions, insertions, and deletions.  We seek to place
several strings of varying lengths in parallel, which may necessitate
adding some spaces.  A basic problem underlying design of algorithms
which seek the best alignment in the sense of minimizing some score
function is simply to count such alignments. Note that for this problem
the elements of the string are irrelevant --- we only care whether there
is a letter or a space.

Mathematically, given positive integers $k$ and a $k$-tuple $\nn = (n_1,
\dots, n_k)$, we may define a $(k;\nn)$-alignment as a $k\times n$
binary matrix for some $n$, such that no columns are identically zero
and the $i$th row sum is $n_i$.  Column $j$ is \Em{aligned} if it has no
zero entry, and a \Em{block} of size $b$ is a $k\times b$ submatrix,
with contiguous columns, all of which are aligned.  For example,
Figure~\ref{fig:alignment} shows an alignment of $3$ strings, and the
corresponding matrix. The generating function giving the number of 
$(k, \nn)$-alignments with the $i$th row having sum $n_i$ is given by 
$$ F(z_1, \dots ,z_k) = \frac{1}{2 - p(\zz)} $$ 
where $p(\zz) = \prod_{j=1}^k (1+z_j)$.

\begin{figure}[ht]
\begin{center}
\begin{tabular}{|c|c|c|c|c|c|c|c|c|c|c|}
\hline
C & C & A & G & T & C & A & G & C & T & A \\ \hline
C &   & A & G & T & C &   & C & C & T & G \\ \hline
G & T & A & G & T & C &   &   & T & T & C \\ \hline
&&&&&&&&&& \\ \hline
1 & 1 & 1 & 1 & 1 & 1 & 1 & 1 & 1 & 1 & 1 \\ \hline
1 & 0 & 1 & 1 & 1 & 1 & 0 & 1 & 1 & 1 & 1 \\ \hline
1 & 1 & 1 & 1 & 1 & 1 & 0 & 0 & 1 & 1 & 1 \\ \hline
\end{tabular}
\end{center}
\caption{A $(3;11,9,9)$-alignment and its corresponding binary matrix.}
\label{fig:alignment}
\end{figure}

The special case where all sequences involved have the same length
(all rows have the same sum) is easily dealt with.  By
Proposition~\ref{pr:classify}, there is a single strictly minimal
point in the positive orthant that controls asymptotics in the
diagonal direction.  By symmetry of $F$ and of the equation $\dir(\zz)
= \one$, this point has the form $\zz = z \one$ for some positive $z$.
Thus we have $z = 2^{1/k} - 1$ and hence the asymptotic has the form
$C (z^{-k})^n n^{-\frac{k-1}{2}}$.  This result was derived by Griggs,
Hanlon, Odlyzko and Waterman in \cite{GHOW1990} using a
saddle point analysis.  Note that for large $k$, $z^{-k} \sim 2^{-1/2}
k^k (\log 2)^{-k}$.  

To compute the constant $C$, we use the formula
of Theorem~\ref{th:smooth nice}.  It is readily computed that $-z_k
\partial H /\partial z_k = 2z/(1 + z)$.  To compute the Hessian
determinant, we need only consider the Hessian in $\theta_1, \dots
,\theta_{d-1}$ of $- \log \prod_{j<k} (1+z \exp(i\theta_j)) + \log (2
- \prod_{j<k} (1+z \exp(i\theta_j))$.  Direct computation using the
fact that $\prod_{j<k} (1+z \exp(i\theta_j))$ takes the value
$2/(1+z)$ when $\btheta = \zero$ shows that the Hessian matrix has
diagonal elements $2/(1+z)$ and off-diagonal elements $1/(1+z)$.  Thus
$\hess = k (1+z)^{-(k-1)} = k 2^{1/k}/2$ and this yields $$C =
\frac{2^{(1-k^2)/2k}}{(2^{1/k} - 1) \sqrt{k \pi^{k-1}}} \, .$$ This
yields the result of \cite{GHOW1990} (note that there is an
error in the displayed formula on \cite[p. 1989]{Wate1995} --
the factor $r^k$ should read $r$, and the factor $2^{(k^2-1)/(2k)}$
should be $2^{(1-k^2)/(2k)}$).

A more important case for biological applications is when the minimum
block size is bounded below, by $b$, say.  The generating function in
this case is shown in \cite{ReTa}, by means of standard
operations on generating functions of formal languages, to have the form
$$
F(z_1, \dots, z_k) = \frac{A}{1 + (1 - p) A + (A - 1) t} 
= \frac{1 - t + t^b}{(1 - t)(1 - g) - t^b g} 
= \frac{1 + \frac{t^b}{1 - t}}{1 - g\left(1 + \frac{t^b}{1 - t}\right)} 
$$
where $t =\prod_j z_i, g = p - 1 - t$ and $A = 1 - t + t^b$. Note that
when $b = 1$ then $A = 1$, and we recover the unrestricted case analysed
above. In this case, asymptotics for the case where all sequences have
equal length  have been derived (as far as we know) only for the case $k
= 2$, using the ``diagonal method" discussed in
Section~\ref{ss:compare}.  However, we can readily deal with general $k$
using the methods of this paper.

Again, by symmetry we need only look for a contributing minimal point
of the form $\zz = z \one$.  Let $H = 1 + (1 - p) A + (A - 1) t$.
Since $t = z^k$ and $p = (1+z)^k$ we must find the root(s) $\rho$ of
smallest modulus of $h(z):= H(z, z, \dots, z)$. Note that the third
formula above shows that we may take $h = 1 - P$ where $P$ is
aperiodic with nonnegative coefficients, so there will be a unique
root $\rho$ of smallest modulus, and it will be positive real. The
exponential growth rate is then $\rho^{-k}$ with polynomial correction
of order $n^{(1 - k)/2}$ as above. 

We note that in the case $k = 2$, such a result was proved in
\cite{GHW1986}, but stated slightly differently. The value
$\tau = \rho^2$ is given as the minimal positive root of the polynomial
$(1 - x)^2- 4x (1 - x + x^b)^2$. Setting $x = z^2$ and factoring the
difference of squares yields $\rho$ as the minimal positive root of
$H(z, z) = 1 - 2z -z^2 + 2z^3 - 2z^{2b+1}$.

To better estimate $\rho:=\rho_b$, note that  $\rho_b < 1$ so it is
reasonable to consider the approximation obtained by setting all
powers involving $z^b$ to zero.  The minimal real zero $\rho_1$ of $(1
- t)(1 - g)$ is asymptotically of order $1/k$ (and equals $2^{1/k} -
1$). Now $\rho_b > \rho_1$ but $\rho_b$ should be close to $\rho_1$.
Indeed, it appears that an iteration scheme based on the fixed point
equation
$$
z = \left(1 + z^k + \frac{1 - z^k}{1 -z^k + z^{kb}}\right)^{1/k} - 1
$$
given by $h(z) = 0$ converges rapidly to $\rho_b$ from starting point
$\rho_1$.

\subsection{Probability that there is an edge in an induced subgraph}
\label{ss:edge}

 From the $n$-set $[n] := \{ 1 , \ldots , n \}$, a collection
of $t$ disjoint pairs is named.  Then a $k$ element subset, 
$S \subseteq [n]$, is chosen uniformly at random.  What is the
probability $p(n,k,t)$ that $S$ fails to contain as a subset any 
of the $t$ pairs?  This question is posed in~\cite{LPSS} 
as a step in computing the diameter of a random Cayley graph 
of a group of cardinality $n$ when $k$ elements are chosen at random 
(the diameter is infinite if the $k$ elements do not generate $G$).  

There are a number of ways of evaluating $p(n,k,t)$, one of which is
by inclusion-exclusion on the number of pairs contained.  This leads
to
\begin{equation} \label{eq:p} p(n,k,t) = \sum_{i=0}^t (-1)^i
\binom{t}{i} \binom{n - 2i}{k - 2i} {\binom{n}{k}}^{-1} = \sum_{i=0}^t
(-1)^i \binom{t}{i} \binom{n - 2i}{n - k}{\binom{n}{k}}^{-1}\, .
\end{equation}
Here, we define the binomial coefficient $\binom{p}{q}$ to be zero 
unless $0 \leq q \leq p$.

The numbers $a(n,k,t) := \binom{n}{k} p(n,k,t)$ are simpler to analyse
via generating functions; since for most purposes only a first-order
asymptotic is sought, we lose nothing in considering $a(n, k, t)$.
Note that from the description above, $a(n, k, t) = 0$ if $k + t > n$,
by the pigeonhole principle (since the complement of $S$ has size less
than $t$, $S$ must contain at least $t+1$ of the $2t$ chosen
elements). From the statement of the problem, $a(n, k, t)$ is not
defined if $2t > n$; however, formula~\eqref{eq:p} still makes sense
in that case, even though it does not define the probability of any
event. In fact, $a(n, k, t)$ can be negative for large $t$.

The most direct approach to finding a trivariate generating function
of $a(n, k, t)$ is to use some well-known bivariate generating
functions $\sum_{i,j} a_{ij} x^i y^j$. If $a_{ij} =  \binom{i+j}{j}$
then the GF is $(1- x - y)^{-1}$, while that for $a_{ij} =
\binom{i}{j}$ is $(1- x(1+y))^{-1}$.  We now compute 
\begin{align*}
\sum_{n, k, t, i}  x^n y^k z^t w^i \binom{t}{i} \binom{n - 2i}{k - 2i} 
& = \sum_{N, K, i, j} x^{N+2i} y^{K+2i} z^{i+j} w^i \binom{i+j}{i}  
   \binom{N}{K} \\
& =  \left(\sum_{i, j} \binom{i+j}{i}  (zwx^2y^2)^i z^j \right) 
   \left(\sum_{N, K} \binom{N}{K} x^N y^K \right) \\
& = \frac{1}{1 - z(1 + wx^2y^2)} \frac{1}{1 - x(1 + y)} ,
\end{align*}
which yields the trivariate generating function
\begin{equation} \label{eq:bin sum GF}
F(x,y,z) = \sum_{n,k,t} a(n,k,t) x^n y^k z^t 
   = \frac{1}{1 - z(1 - x^2 y^2)} \frac{1}{1 - x(1+y)} \, .
\end{equation}

Note that if we impose the restriction $2t \leq n$, then the sum over
$N, K$ is restricted to $N \geq 2j$. Now summing over $N, K, t, i$ as
above we obtain the restricted trivariate generating function
$$
F_2(x, y, z) = \sum_{\{n,k,t \colon 2t \leq n\}} a(n,k,t) x^n y^k z^t
=  \frac{1}{1 - x(1+y)}\frac{1}{1 - zx^2(1+2y)}. 
$$
The advantage of $F_2$ is that all its coefficients are nonnegative.

We show how to give asymptotics for $a(n,k,t)$ in an arbitrary
direction for which $n > k + t$.  For the specific case of $k =
\lfloor cn \rfloor, t = (n - 4)/12$ with $0 < c < 1$, such an analysis
was posed as an open question at the end of an early draft
of~\cite{LPSS}.

It turns out that the computations can be carried out equally well
with $F$ or $F_2$. However, with $F_2$ we can use
Theorem~\ref{th:universal} if we are only interested in the behaviour
of $a(n, k, t)$ for $2t\leq n$.

The denominator of $F$ factors into two smooth pieces, call them
$\denom_1 := 1 - x(1+y)$ and $\denom_2 := 1 - z(1-x^2y^2)$.  There is
a corresponding stratification of $\sing$ into two surfaces, $\sing_1
:= \{ \denom_1 = 0 \neq \denom_2 \}$ and $\sing_2 = \{ \denom_2 = 0
\neq \denom_1 \}$, and a curve, $\sing_0 := \{ \denom_1 = \denom_2 = 0
\}$.  For $\zz \in \sing_1$, lack of dependence on $z$ means that
$\dir (\zz) \perp (0,0,1)$, so $\dir$ cannot be in the strictly
positive orthant; hence for $\rbar \in \orthant$, there are no points
of $\crit_{\rbar}$ in $\sing_1$. It turns out there are no points of
$\crit_{\rbar}$ in $\sing_2$ either. This is discovered by
computation.  If $\zz \in \sing_2 \cap \crit_{\rbar}$ and $\rbar =
\overline{(r,s,1)}$, then $\zz$ satisfies the equations
\begin{align*}
\denom_2 (\zz) & =  0 \, ;  \\[2ex]
x \frac{\partial \denom_2}{\partial x} 
   - r z \frac{\partial \denom_2}{\partial z} & =  0 \, ; \\[2ex]
y \frac{\partial \denom_2}{\partial y} 
   - s z \frac{\partial \denom_2}{\partial z} & =  0 \, .
\end{align*}
These equations turn out to have no solutions: Maple tells us in an
instant that the ideal generated by the left-hand sides of the three
above equations in $\C[x,y,z,r,s]$ contains $r - s$; thus $\sing_2$
may contain points of $\crit_{\rbar}$ only when $r = s$, that is, only
governing asymptotics of $a(n,k,t)$ for which $n = k$, which are not
interesting.

Evidently, $\crit_{\rbar} \subseteq \sing_0$.  Two equations are 
$\denom_1 = \denom_2 = 0$, that is, 
$$
(x,y,z) = \left ( \frac{1}{1+y} \, , \, y \, , \frac{(1+y)^2}{1+2y} 
   \right ) \, .
$$
The last equation is that $\rbar$ is in the linear space spanned by
$\loggrad \denom_1$ and $\loggrad \denom_2$.  Setting the determinant
of $(\rbar, \loggrad \denom_1, \loggrad \denom_2)$ equal to zero gives
the equation
\begin{equation}
\label{eq:cayley-crit}
2(r - s - 1) y^2 + (r - 3s) y - s = 0.
\end{equation}
Note that exactly the same equations are obtained for $x, y, z$ when
using $F_2$.

If $r > s + 1$, the discriminant $(r+s)^2-8s$ of the quadratic in
\eqref{eq:cayley-crit} is positive, and  \eqref{eq:cayley-crit} has
one positive and one negative root, the positive root being
$$
y_+ := \frac{\sqrt{(r + s)^2 - 8s} - (r - 3s)}{4(r - s + 1)} \, .
$$


Plugging into for $x_+$ and $z_+$ yields expressions for these which
are also quadratic over $\Z[r,s]$.  The expression for $(x_-,y_-,z_-)$
is just the algebraic conjugate of $(x_+,y_+,z_+)$, but is negative
in the second and third coordinates. 

Having identified $\crit_{\rbar}$ as these two conjugate points, which
we will call $\zz_+$ and $\zz_-$, it remains to find
$\contrib_{\rbar}$. Since the two elements lie on different tori, we
may conclude from Theorem~\ref{th:universal} that $\contrib_{\rbar} =
\{ \zz_+ \}$. 

The form of the leading term asymptotic is then given by
Theorem~\ref{th:m<d}:  
$$
a(n,k,t) \sim C \left (\frac{k}{n} , \frac{t}{n} \right )
n^{-1/2} x_+^{-n} y_+^{-k} z_+^{-t} \, .
$$
For example, in the direction $\overline{(12, 12c, 1)}$ with $c <
11/12$, we can compute the exponential growth rate and compare with
that of $\binom{n}{\lfloor cn \rfloor}$, and thereby show that certain
random Cayley graphs have diameter 2 with very high probability. For
details, see~\cite{LPSS}.

Analysis in the unrestricted case where we may have $2t > n$ is more
difficult. In fact $[x^ny^kz^t] F$ need not be
zero when $k + t > n$ and $2t > n$, and may be negative. 

\subsection{Integer solutions to linear equations}
\label{ss:integer}

Let $a_\rr$ be the number of nonnegative integer solutions to $A \xx =
\rr$ where $A$ is a $d \times m$ integer matrix. Denote by $\vb^{(k)}$
the $k^{th}$ column of $A$.  Then $$F(\zz) = \sum_\rr a_\rr \zz^\rr =
\prod_{k=1}^m \frac{1} {1 - \zz^{\vb^{(k)}} }$$ This enumeration
problem has a long history.  We learned of it
from~\cite{DeSt2003} (see
also~\cite[Section~4.6]{Stan1997}). One special case is when $A =
A_{m,n}$ is the incidence matrix for a complete bipartite graph
between $m$ vertices and $n$ vertices. Solutions to $A \xx = \rr$
count nonnegative integer $m \times n$ matrices with row and column
sums prescribed by $\rr$; enumerating these is important when
constructing statistical tests for contingency tables.  Another
special case is when $A = A_n$ is the incidence matrix of a complete
directed graph on $n$ vertices, directed via a linear order on the
vertices.  Here solutions are counted in various cones other than the
positive orthant and the function enumerating them is known as
Kostant's partition function. 

It is known that $a_\rr$ is piecewise polynomial, and it has been a
benchmark problem in computation to determine these polynomials, and
the regions or \Em{chambers} of polynomiality, explicitly. Indeed,
several subproblems merit benchmark status.  Counting the chambers for
Kostant's partition function is one such problem.  Another is to
evaluate the leading term for the diagonal polynomial $a_{n \one}$ in
the case where $A = A_{k,k}$, the so-called \Em{Ehrhart polynomial},
which counts $k \times k$ nonnegative integer matrices with all rows
and columns summing to $n$.  Equivalently, this counts integer points
in the $n$-fold dilation of the \Em{Birkhoff polytope}, defined as the
set in $\R^{\binom{k}{2}}$ of all doubly stochastic $k \times k$
matrices; the leading term of $a_{n \one}$ is the volume of the
Birkhoff polytope.

Our methods do not improve on the computational efficiency of previous
researchers: our ending point is a well known representation, from
which other researchers have attempted to find efficient means of
computing.  Our methods do, however, give an effective answer, which
illustrates that this class of problems can be put in the framework
for which our methods give an automatic solution.  The remainder of
this section is devoted to such an analysis. 

Our examination of the pole variety $\sing$ begins with an answer
rather than a question: we see that $\one \in \sing$ and realize that 
life would be easy if $\contrib_{\rbar} = \{ \one \}$ for all $\rbar$.  
In order for $\one$ to be a singleton stratum it is necessary and
sufficient that the columns $\vb^{(k)}$ span all of $\R^d$.  For the
remainder of the section we assume this to be the case.  This satisfies
our standing assumption~\ref{ass:hole}.

The variety $\sing$ is the union of the smooth sheets $\sing_k$, where
for $1 \leq k \leq m$, $\sing_k$ is the binomial variety $\{ \zz :
\denom_k (\zz) := \zz^{\vb^{(k)}} - 1 = 0 \}$. On each of these
varieties, $\dir_k (\zz) = \loggrad (\denom_k) (\zz)$ is constant and
equal to $\vb^{(k)}$.  If a stratum $S$ is bigger than a singleton,
then at any $\zz \in S$, the vectors $\loggrad \denom_k$ span a proper
subspace of $\R^d$; but these vectors do not vary as $\zz$ varies in
$S$, so the union over $\zz \in S$ of $\linear (\zz)$ is this same
proper subspace of $\R^d$.  Consequently, the union over all points of
all non-singleton strata of $\linear (\zz)$ is a union of proper
subspaces, which we denote $\dirset'$.  For $\rbar \notin \dirset'$,
then, $\contrib_{\rbar}$ consists of one or more singleton strata.  

Taking logs, we see that $\log \sing_k$ is a hyperplane normal to
$\vb^{(k)}$ and is central (passes through the origin).  We see that
$\cone (\zero)$ is the positive hull of the vectors $\vb^{(k)}$; this
hull is $\dirset$ and outside of the closure of this, $a_\rr$
vanishes.  Forgetting about the logs, we see that $\one \in
\contrib_{\rbar}$ for all $\rbar \in \dirset_0 := \dirset \setminus
\dirset'$: this follows from Theorem~\ref{th:m>=d} since $\rbar \in
\cone (\one)$.

If the intersection of all the surfaces $\sing_k$ contains any other
points on the unit torus $\torus (\one)$ then these too are in
$\contrib_{\rbar}$ for all $\rbar \in \dirset_0$.  This is easy to
check, since it is equivalent to the integer combinations of the
columns of $A$ spanning a proper sublattice of $\Z^d$.  For instance,
in the example of counting matrices with constrained row and column
sums, the columns of $A$ span the alternating sublattice of $\Z^d$,
corresponding to the fact that $\contrib_{\rbar} = \{ \one , -\one \}$
for all $\rbar \in \dirset_0$.  

There may  be singleton strata given by intersections of 
subfamilies of $\{ \sing_k : 1 \leq k \leq m \}$, that lie on
the unit torus, but in computing the leading term asymptotics these
may be ignored because they yield polynomials in $\rr$ of lower degree.

In summary, for any $\rbar \in \dirset_0$, the leading term
asymptotics are given by summing~(\ref{eq:M>=d}) over a set containing
$\one$, and isomorphic to the quotient of $\Z^d$ by the integer span
of the columns of $A$.  Our \emph{a priori} knowledge that $a_\rr$ are
integers, teamed with~(\ref{eq:highpole}) as in the last part of
Theorem~\ref{th:m>=d}, shows that in fact $a_\rr$ are piecewise
polynomial, at least away from $\dirset'$.  Thus, except on a set of
codimension 1, we recover the well known piecewise polynomiality of
$a_\rr$.

In their paper, de Loera and Sturmfels use as a running example
the matrix 
$$
A = \left [ 
\begin{array}{ccccc} 
1&0&0&1&1 \\ 0&1&0&1&0 \\ 0&0&1&0&1 
\end{array} \right ] \; .
$$
In this example we see that the columns of $A$ span $\Z^3$ over $\Z$,
so the only contributing point is $\one$.  The first three $\vb$
vectors, in the order given, are the standard basis vectors, so the
cone $\dirset$ is the whole positive orthant of $\RP^2$, which is a
2-simplex.  The other two $\vb$ vectors are two of the three face
diagonals of this cone, which are two midpoints of edges of the
2-simplex in $\RP^2$.  In addition to the boundary of $\dirset$, there
are three projective line segments in $\dirset'$, corresponding to the
one-dimensional strata $\sing_4 \cap \sing_5$, $\sing_4 \cap \sing_3$
and $\sing_5 \cap \sing_2$.  The complement of these three line
segments (two medians of the simplex $\RP_2$ and the line segment
connecting two midpoints of edges) divides $\RP^2$ into five chambers,
and on each of these chambers $a_\rr$ is a quadratic polynomial.
Algorithms for computing these polynomials are given, for example,
in~\cite{DeSt2003} or~\cite[Section~5]{BaPe}.

\subsection{Queueing theory} \label{ss:queueing}

Queuing theory describes the evolution of a collection of
\Em{jobs} that enter the system, exit the system and change
types stochastically as they are worked on by one or more
\Em{servers}.  There are many variations on this model.  The
example we consider here is a closed network with no class-hopping,
meaning that no jobs enter, leave or change type.  The model
is defined in terms of the parameters $J, K, \{ \lambda_j :
1 \leq j \leq J \}, \{ \mu_{jk} : 1 \leq j \leq J , 1 \leq k \leq K \}$, 
and $\{ p_{j;kl} : 1 \leq j \leq J , 0 \leq k,l \leq K \}$ as follows.  

Let $J$ be the number of types of jobs.  Let $K+1$ be the number
of servers.  In our example, server~$k$ $(1 \leq k \leq K)$ 
distributes its resources evenly among all jobs in its queue,
serving jobs of type $j$ at rate $\mu_{jk} / n$ when there
are $n$ jobs in its queue, while server~0
is an \Em{infinite} server, which serves all waiting jobs 
simultaneously, serving a job of type $j$ at rate $\lambda_j$
independent of how many jobs it is serving.  A job of type of
type $j$ at a station $k \in \{ 0 , \ldots , K \}$ waits to 
be at the front of the queue (if $k \geq 1$), then waits to be
served, then moves to server~$l$ with probability $p_{j;kl}$.  
All service times are independent exponential random variables
and the movement between servers is Markovian and independent
of the service times.

The entire model is Markovian if one takes the state vector to
be the $J \times (K+1)$ matrix $Q = Q(t)$ giving the number of jobs of 
each type waiting at each server.  Let $n_j = n_j (t) = \sum_{k=0}^K
Q_{jk} (t)$ be the number of jobs of type $j$ present in the
system at time $t$.  The assumptions of our model (closed, with 
no class-hopping) imply that the vector $\nn = (n_1 , \ldots , n_J)$
is constant.  Assuming also that enough of the $p_{j;kl}$ are nonzero
that any job at any server can get to any other server, the process
started from any $\nn = \nn(0)$ will be an ergodic Markov chain on 
the set of states $S(\nn)$ of $J \times (K+1)$ matrices with row 
sums $\nn$.  There will be a unique stationary distribution $\pi_\nn$.
Explicit product formulae exist for $\pi_\nn (Q)$ of the form
$\pi_\nn (Q) = \frac{1}{G (\nn)} P_\nn (Q)$, where $P$ is a large product
and $G (\nn)$ is the normalizing constant or \Em{partition function}.
Typically, $P_\nn$ is easier to estimate than $G (\nn)$, so much of
the work in analyzing this sort of queueing network is in estimating
$G (\nn)$.  

It turns out that there is a relatively simple generating function
for the quantities $G (\nn)$.  One must first derive quantities 
$\rho_{ji}$, which are stationary probabilities for a single job
migrating through the network.  The \Em{partition generating function}
is then given by
\begin{equation}
F(\zz) := \sum_\nn G(\nn) \zz^\nn = \frac{\exp (z_1 + \cdots
   + z_j)}{\prod_{i=1}^K \left ( 1 - \sum_{j=1}^J \rho_{ji} 
   z_j \right ) } \, .
\end{equation}
This equation is given in~\cite[Equation~(44)]{Koga2002}
(though the definition of $\rho_{ji}$ here is sketchy) or 
in~\cite[Equation~(2.26)]{BeMc1993} for
a slightly more general model (more than one infinite server).

Evidently, the denominator of $F$ decomposes into linear factors.
A general analysis of this case, including an algorithm for
identifying $\contrib$, is provided in~\cite{BaPe}.  Here, we
work the simplest nontrivial case, $J=K=2$.  Thus
$$F(x,y) = \frac{\exp(x+y)}{(1 - \rho_{11}\, x - \rho_{21}\, y) 
   (1 - \rho_{12}\, x - \rho_{22}\, y)}$$
for positive constants $\rho_{ji}$.  This example is worked
in~\cite{BeMc1993} but without the infinite
server.  

Assume that $\rho_{11} > \rho_{12}$, which loses no generality 
except for allowing $\rho_{11} = \rho_{12}$.  The most interesting 
case is then $\rho_{22} > \rho_{21}$.  The singular variety $\sing$ 
consists of two lines with negative slopes, and we have assumed 
order relations implying that the lines intersect in the positive
quadrant.  This is exactly what is shown in Figure~\ref{fig:lines}.

\begin{figure}[ht]
\includegraphics[scale=0.2]{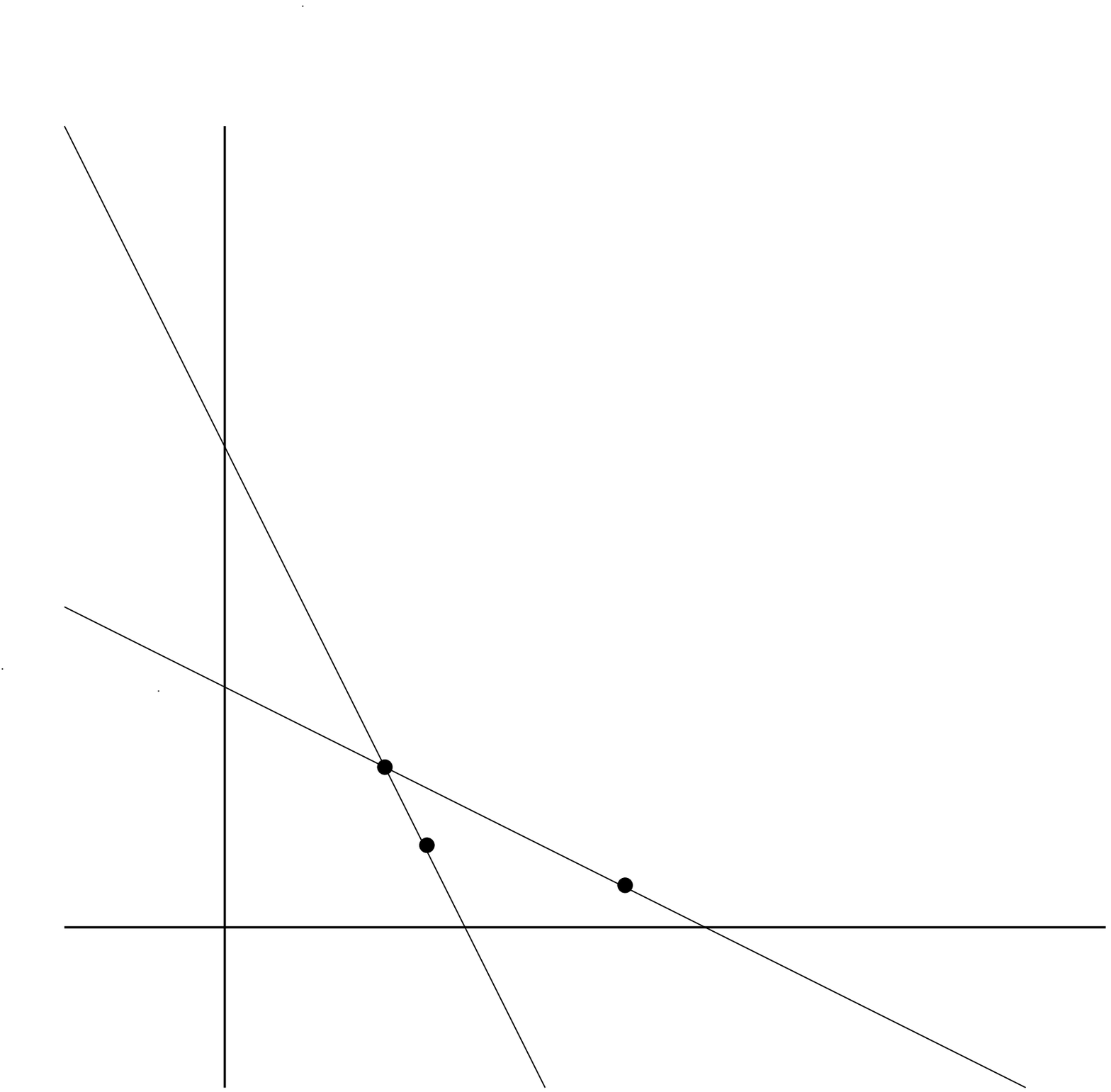}
\hspace{0.5in}
\includegraphics[scale=0.2]{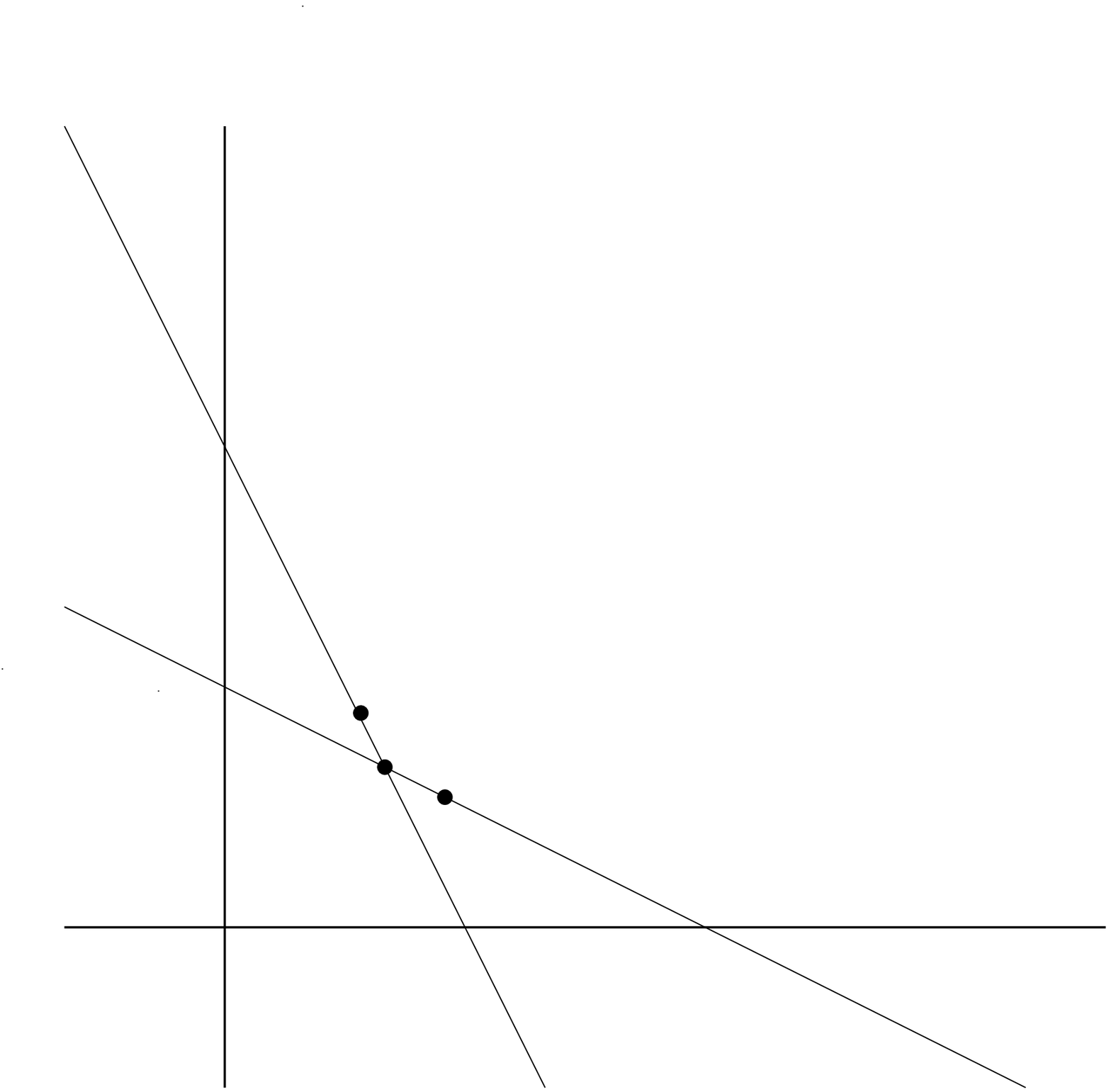}
\hspace{0.5in}
\includegraphics[scale=0.2]{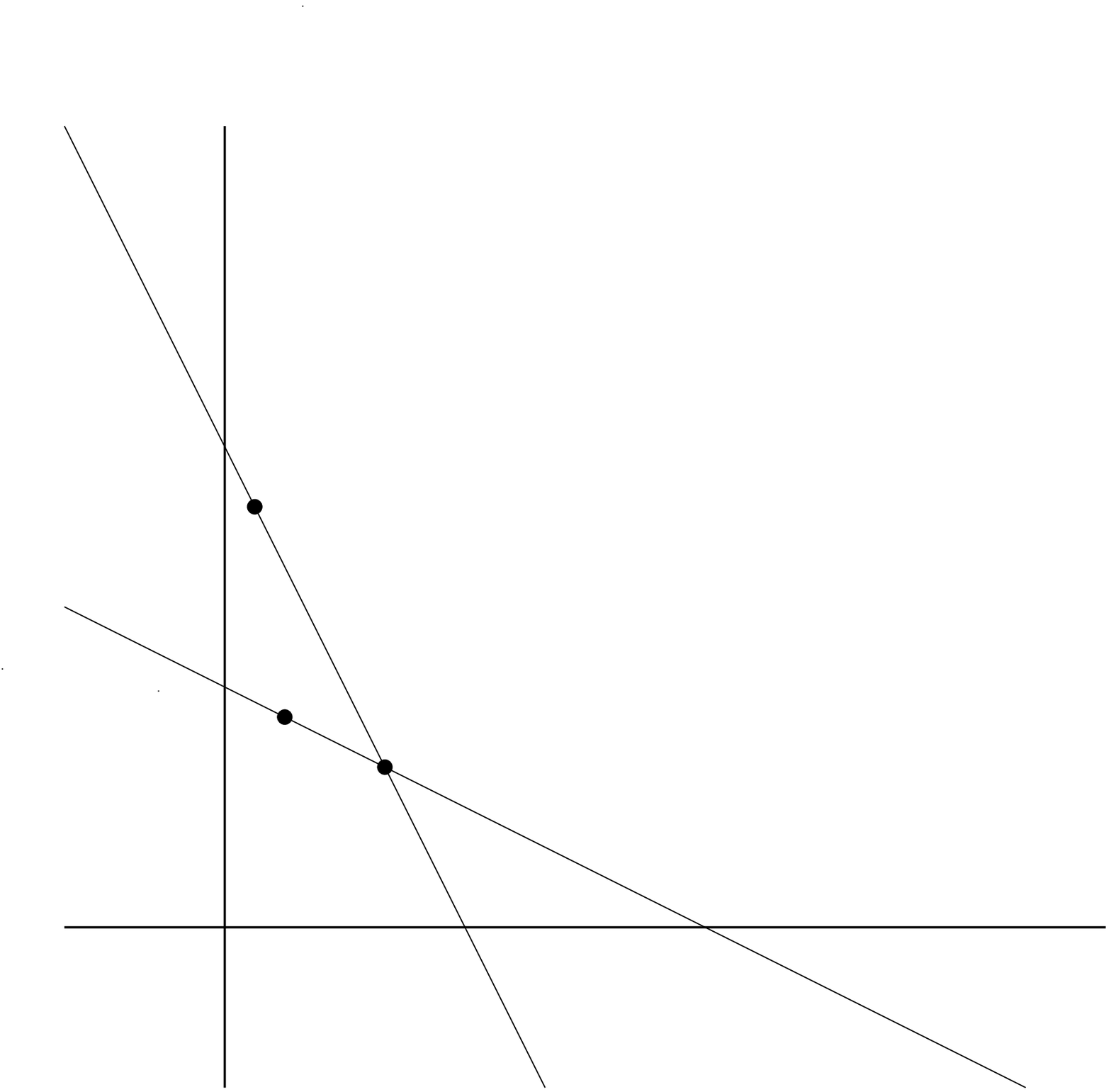}
\caption{The three possible configurations of critical points}
\label{fig:three cases}
\end{figure}

To recall the information from that figure, let
$$D := \rho_{11} \rho_{22} - \rho_{12} \rho_{21}$$ 
denote the determinant of the matrix $[ \rho_{ji} ]$. 
The set $\sing$ has three strata, one being the point 
$$(x_0 , y_0) := \left ( \frac{\rho_{22} - \rho_{21}}{D} ,
   \frac{\rho_{11} - \rho_{12}}{D} \right )$$
where the lines meet and the others being the two lines, 
each with the point $(x_0,y_0)$ removed.  Whenever a stratum is linear, 
there is a unique solution to the critical point equations~\eqref{eq:zero-dim}
on the stratum (use concavity of the logarithm to see this), hence
$\crit_{\rbar}$ has cardinality~3 unless two of the critical points 
coincide.  One of the critical points is always the point
$(x_0 , y_0)$.  Figure~\ref{fig:three cases} shows the location of
all three critical points as $\rbar$ varies from near horizontal
to near diagonal to near vertical.  

The three configurations correspond to a division of 
the positive orthant of $\RP^1$ into three regions.  
Transitions between the regions occur when $\rbar$ is equal to
$$\lambda_1 := \frac{\rho_{11} (\rho_{22} - \rho_{21})}
   {\rho_{21} (\rho_{11} - \rho_{12})}$$
or
$$\lambda_2 := \frac{\rho_{12} (\rho_{22} - \rho_{21})}
   {\rho_{22} (\rho_{11} - \rho_{12})} \, .$$
The middle region is the $\cone (x_0 , y_0)$ of directions
between $\lambda_1$ and $\lambda_2$.  
In the interior of this cone, asymptotics are given by
Corollary~\ref{cor:m=d}.  The determinant of the Hessian 
of the denominator of $F$,
evaluated at $(x_0 , y_0)$, is equal to $-D^2$.  The corollary
therefore yields
$$a_{rs} \sim x_0^{-r} y_0^{-s} \frac{\exp (x_0 + y_0)}
   {x_0 \, y_0 \, D} \, .$$

To be completely explicit, consider the denominator from 
Figure~\ref{fig:lines}.  Here $\rho_{11} = \rho_{22} = 2/3$ 
and $\rho_{12} = \rho_{21} = 1/3$.  The boundaries of the cone 
have slopes $1/2$ and~2.  Both $x_0$ and $y_0$ are equal to~1, 
so for $(r,s) \to \infty$ in a compact subcone of $\cone (x_0 , y_0)$,
$$a_{rs} = \frac{e^2}{1/3} = 3 e^2 \approx 22.1671 \, .$$
The actual value for $r=s=80$ is $a_{80,80} \approx 22.1668$ and the
the relative error is $1.7 \times 10^{-5}$.  Even at $r=s=40$, 
the error is only one quarter of one percent.  
As in Section~\ref{ss:distinct subseq}, there is
a Gaussian scaling window at the boundary of $\cone (x_0 , y_0)$.
Thus, for example, $a_{n,2n} \sim \frac{3}{2} e^2 \approx 11.083 \,$.  
The error decays more slowly on the boundary (order $n^{-1/2}$, 
though our treatment here does not extend to this depth), 
so for instance $a_{40,80} \approx 10.893$ differs from the limiting value 
of $\frac{3}{2} e^2$ by $1.7\%$.  The standard deviation in the 
scaling window is a little under $\sqrt{n}$: for example, when 
$n = 80$ the standard deviation is around~8.  Checking $a_{48,80}$ 
we see a value of roughly $18.45$, which is $83\%$ of the 
interior value of $3 e^2$ while at two standard deviations above 
we get $a_{56,80} \approx 21.49$ or $97\%$ of the interior value; 
at one and two standard deviations below we have 
$a_{32,80} \approx 3.018$ and $a_{24,80} \approx 0.241$ which are 
respectively $13.6\%$ and $1\%$ of the interior value.

\setcounter{equation}{0}
\section{Transfer matrices}
\label{sec:transfer}

Suppose we have a class $\CC$ of combinatorial objects that may be put
in bijective correspondence with the collection of paths in a 
finite directed graph.  For example, it is shown 
in~\cite[Example~4.7.7]{Stan1997} that the class $\CC$ of 
permutations $\pi \in \bigcup_{n=0}^\infty \mathcal{S}_n$ 
such that $|\pi (j) - j| \leq 1$ for all $j$ in the domain of $\pi$
is of this type.  Specifically, the digraph has vertices
$V := \{ -1 , 0 , 1 \}^2 \setminus \{ (0,-1) , (1,0) \}$, with
a visit to the vertex $(a,b)$ corresponding to $j$ for which
$\pi (j) = a$ and $\pi (j+1) = b$; the vertices $(0,-1)$ and $(1,0)$
do not occur because this would require $\pi (j) = \pi (j+1)$,
which cannot happen; the edge set $E$ is all elements of $V^2$ except
those connecting $(1 , x)$ to $(y, -1)$ for some $x$ and $y$, 
since this would require $\pi(j) = \pi(j+2)$; the correspondence 
is complete because once $\pi (j), \pi (j+1)$ and $\pi(j+2)$ 
are distinct, no more collisions are possible. 

The \Em{transfer-matrix method} is a general device for producing 
a generating function which counts paths from a vertex $i$ to a 
vertex $j$ according to a vector weight function $W$ such that
$W(v_1 , \ldots , v_n) = \sum_{j=1}^n w(v_i , v_j)$ for some function 
$w: V \times V \to \N^d$.  
\begin{pr}[Theorem~4.7.2~of~\protect{\cite{Stan1997}}] 
\label{pr:transfer matrix}
Let $A$ be the weight matrix, that is, the 
matrix defined by $A_{ij} = \delta_{i,j} \zz^{w(i,j)}$,
where $\zz = (z_1 , \ldots , z_d)$ and $\delta_{i,j}$ denotes~1
if $(i,j) \in E$ and~0 otherwise.  Let $\CC_{ij}$ be the subclass of
paths starting at $i$ and ending at $j$ and define $F(\zz) 
= \sum_{\gamma \in \CC_{ij}} \zz^{W(\gamma)}$.  
Then 
$$F (\zz) = \left ( (I-A)^{-1} \right )_{ij} 
   = \frac{(-1)^{i+j} \det (I-A:i,j)}{\det(I-A)}$$
where $(M:i,j)$ denotes $(i,j)$-cofactor of $M$, that is,
$M$ with row $i$ and column $j$ removed.
$\Cox$
\end{pr}

Thus any class to which the transfer-matrix method applies will have a
multivariate rational generating function.  It is obvious, for
example, that the transfer-matrix method applies to the class of words
in a finite alphabet with certain transitions forbidden. The examples
discussed in~\cite[pages~241--260]{Stan1997} range from restricted
permutations to forbidden transition problems to the derivation of
polyomino identities including~(\ref{eq:HCP by n}), and finally to a
very general result about counting classes that factor when viewed as
monoids.  Probably because techniques for extracting asymptotics were
not widely known in the multivariate case, the discussion
in~\cite{Stan1997} centers around counting by a single variable:
in the end all weights are set equal to a single variable, $x$,
reducing multivariate formulae such as~(\ref{eq:HCP by n,k}) to
univariate ones such as~(\ref{eq:HCP by n}).  But the methods in the
present paper allow us to handle multivariate rational functions
almost as easily as univariate functions, so we are able to derive
joint asymptotics for several statistics at once, which is useful to
the degree that we care about joint statistics.

\subsection{Restricted switching}
\label{ss:switching}

Some of the examples we have already seen, such as enumeration
of Smirnov words, may be cast as transfer-matrix computations.  
A simpler example of the method is the following path counting problem.

Let $G$ be the graph on $K+L+2$ vertices which is the 
union of two complete graphs of sizes $K+1$ and $L+1$ with
a loop at every vertex and one edge $\overline{xy}$ between them.  
Paths on this graph correspond
perhaps to a message or task being passed around two workgroups,
with communication between the workgroups not allowed except
between the bosses.  If we sample uniformly among paths of
length $n$, how much time does the message spend, say, among the common
members of group~1 (excluding the boss)?

We can model this efficiently by a four-vertex graph, with vertices
$C_1,B_1,B_2$ and $C_2$, where $B$ stands for boss (so $x=B_1$ and $y=B_2$)
and $C$ refers to the workgroup.  We have collapsed $K$
vertices to form $C_1$ and $L$ to form $C_2$, so a path whose 
number of visits to $C_1$ is $r$ and to $C_2$ is $s$ will count
for $K^r L^s$ actual paths.

Let $F(u,v,z) = \sum_\omega u^{N(\omega , C_1)} v^{N(\omega , C_2)} 
z^{|\omega|}$
where $N(\omega , C_j)$ denotes the number of visits by $\omega$ to $C_j$
and $|\omega|$ is the length of the path $\omega$.  
By Proposition~\ref{pr:transfer matrix}, the function $F$ is rational
with denominator $\denom = \det (I-A)$, where 
$$
A = \left [ \begin{array}{cccc}
   uz & uz & 0 & 0 \\
   z & z & z & 0 \\
   0 & z & z & z \\
   0 & 0 & vz & vz
   \end{array} \right ]
$$
and the states are ordered $C_1 , B_1 , B_2 , C_2$.
Subtracting from the identity and taking the determinant yields 
$$
\denom =
u{z}^{2}+u{z}^{2}v-uz-u{z}^{4}v+{z}^{2}v-2\,z-zv+1+{z}^{3}v+u{z}^{3}
\, .
$$
Giving each visit to $C_1$ weight $K$ corresponds to the substitution
$u=K u'$ and similarly, $v = L v'$.  Applying Theorem~\ref{th:WLLN}
and reframing the results in terms of the original $u$ and $v$ shows
that the proportion of time a long path spends among $C_1$ tends
as $n \to \infty$ to 
$$
\left. \frac{u \partial \denom / \partial u}{z \partial \denom / 
   \partial z} \right |_{(K,L,z_0)}
$$
where $z_0$ is the minimum modulus root of $\denom (K,L,z)$.  In the  
present example, this yields
$$
\left. {\frac {K \left( -z-zL+1+{z}^{3}L-{z}^{2} \right)
}{-2\,Kz-2\,KzL+K+4 \,K{z}^{3}L-2\,zL+2+L-3\,{z}^{2}L-3\,K{z}^{2}}}
\right |_{z=z_0}
$$
for the proportion of time spent in $C_1$, where $z_0$ is the minimal
modulus root of $\denom (K,L,z)$.  Trying $K=1$, $L=1$ as an example,
we find that $H(1,1,z) = 1 - 4z + 2 z^2 + 3 z^3 - z^4$, leading to
$z_0 \approx 0.381966$ and a proportion of just over $1/8$ for the
time spent in $C_1$.  There are 4 vertices, so the portion of time the
task spends at the isolated employee is just over half what it would
have been had bosses and employees had equal access to communication.
This effect is more marked when the workgroups have different sizes.
Increasing the size of the second group to~2, we plug in $K=1,L=2$ and
find that $z_0 \approx 0.311108$ and that the fraction of time spent
in $C_1$ has plummeted to just under $0.024$.

\subsection{Connector matrices}
\label{ss:connector}

For a more substantial example, we have chosen a further refinement of
the transfer-matrix method, called the \Em{connector matrix method}.
The specific problem we will look at is the enumeration of sequences
with forbidden substrings, enumerated by the composition of the
sequence.  This is a problem in which we are guaranteed (by the
general transfer matrix methodology) to get a multivariate rational
generating function, but for which a more clever analysis greatly
reduces the complexity of the computations.  Our discussion of this
method is distilled from Goulden and Jackson's
exposition~\cite{GoJa2004}.

Let $S$ be a finite alphabet and let $T$ be a finite set of words
(elements of $\bigcup_n S^n$). Let $\CC$ be the class of words with no
substrings in $T$, that is, those $(y_0 , \ldots , y_n) \in \CC$ such
that $(y_j , \ldots , y_{j+k}) \notin T$ for all $j,k \geq 0$. We may
reduce this to a forbidden transition problem with vertex set $V :=
S^k$ with $k$ the maximum length of a word in $T$. This proves the
generating function will be rational, but is a very unwieldy way to
compute it, involving a $|S|^k$ by $|S|^k$ determinant.

The connector matrix method of~\cite{GuOd1981b}, as presented
in~\cite[pages~135--136]{GoJa2004}, finds a much more efficient
way to solve the forbidden substring problem (see also an
 elementary proof of these results via
martingales~\cite{Li1980}).  They discover that it is easy to
enumerate sequences containing at least $v_i$ occurrences of the
forbidden substring $\omega_i$.  By inclusion-exclusion, one can then
determine the number containing precisely $y_i$ occurrences of
$\omega_i$, and then set $\yy = \zero$ to count sequences entirely
avoiding the forbidden substrings.

To state their result, let $T = \{ \omega_1 , \ldots , \omega_m \}$ be
a finite set of forbidden substrings in the alphabet $[d] := \{ 1 ,
\ldots , d \}$.  We wish to count words according to how many
occurrence of each symbol $1 , \ldots , d$ they contain and how many
occurrences of each forbidden word they contain.  Thus for a word
$\eta$, we define the weight $\tau (\eta)$ to be the $d$-vector
counting how many occurrences of each letter and we let $\sigma
(\eta)$ be the $m$-vector counting how many occurrences of each
forbidden substring (possibly overlapping) occur in $\eta$.  Let $F$
be the $(d+m)$-variate generating function with variables $x_1 ,
\ldots , x_d , y_1 , \ldots , y_m$ defined by
$$
F(\xx , \yy) := \sum_{\eta} \xx^{\tau (\eta)} \yy^{\sigma (\eta)} 
   \, .
$$
\begin{pr}[Theorem~2.8.6 and Lemma~2.8.10
of~\protect{\cite{GoJa2004}}]
\label{pr:connector matrix}
Given a pair of finite words, $\omega$ and $\omega'$, let
$\connect (\omega,\omega')$ denote the sum of the weight of 
$\alpha$ over all words $\alpha$ such that some initial segment 
$\beta$ of $\omega'$ is equal to a final segment of $\omega$ and 
$\alpha$ is the initial unused segment of $\omega$.  Formally,
$$
\connect (\omega,\omega') = \sum_{\alpha : (\exists \beta,\gamma ) \,
\omega = \alpha \beta , \omega' = \beta \gamma} \tau (\alpha) \, .
$$
Let $\V$ be the square matrix of size $m$ defined by $\V_{ij} =
\connect (\omega_i , \omega_j)$, denote the diagonal matrices $\Y :=
\diag (y_1 , \ldots , y_m)$, $\L := \diag (\xx^{\tau (\omega_1)} ,
\ldots , \xx^{\tau (\omega_m)})$, and let $\J$ to be the $m$ by $m$
matrix of ones.  Then
$$
F(\xx,\yy) = \left [ 1 - (x_1 + \cdots + x_d ) 
   - C(\xx , \yy - \one) \right ]^{-1}
$$
where
$$
C(\xx , \yy) = \trace \left ( (\I - \Y \V)^{-1} \Y \L \J \right ) \, .
$$
In particular, setting $\yy = \zero$, the generating
function for the words with no occurrences of any forbidden
substring is
\begin{equation} 
\label{eq:connect}
F(\xx , \zero) = \left [ 1 - (x_1 + \cdots + x_d) 
   + \trace \left ( (\I + \V)^{-1} \L \J \right ) \right ]^{-1} \, .
\end{equation}
$\Cox$
\end{pr}

\begin{unremark} 
The function $C$ is the so-called cluster generating function, whose $
(\rr , \sss)$-coefficient counts strings of composition $\rr$ once for
each way that the collection of $s_j$ occurrence of the substring
$\omega_j$ may be found in the string.
\end{unremark}

\subsection{Forbidden substring example}
\label{ss:forbidden}

Goulden and Jackson apply the results of the previous subsection to an
example.  Let $S = \{ 0 , 1 \}$ be the binary alphabet, and let $T =
\{ \omega_1 , \omega_2 \} = \{ 10101101 , 1110101 \}$.  The final
substrings of length~1 and~3 of $\omega_1$ are initial substrings of
$\omega_1$ with corresponding leftover pieces 1010110 and 10101.  Thus
$V_{11} = x_1^3 x_2^4 + x_1^2 x_2^3$.  Computing the other three
entries similarly, we get
$$
\V = \left [ \begin{array}{cc} x_1^3 x_2^4 + x_1^2 x_2^3 & x_1^3 x_2^4 
   \\ x_1^2 x_2^4 + x_1 x_2^3 + x_2^2 & x_1^2 x_2^4 
   \end{array} \right ] \, .
$$
We also obtain 
$$
\L =  \left [ \begin{array}{cc} x_1^3 x_2^5 & 0 \\
   0 & x_1^2 x_2^5 \end{array} \right ] \, .
$$
Plugging into~(\ref{eq:connect}), denoting $(x_1 , x_2)$ by $(x,y)$,
finally yields
\begin{equation} 
\label{eq:connect eg}
F(x,y) = \frac{1 + x^2 y^3 + x^2 y^4 + x^3 y^4 - x^3 y^6}
{1 - x - y + x^2 y^3 - x^3 y^3 - x^4 y^4 - x^3 y^6 + x^4 y^6} \, .
\end{equation}

We make the usual preliminary computations on $F$.  Write $\numer$ and
$\denom$ for the numerator and denominator in~(\ref{eq:connect eg}). A
Gr\"{o}bner basis computation in Maple quickly tells us $\denom$
has no singularities and that $\denom = \numer = 0$ does not occur in
the positive quadrant of $\R^2$ (the last command finds real roots of
the relevant polynomial):
\begin{verbatim}
> Basis([H , diff(H,x) , diff(H,y)] , plex(x,y));

                              [1]
> gb:=Basis([G,H],plex(x,y));

               4    3      2       3            2
             [y  + y  - 2 y  + 1, y  + x - y + y  + 1] 
             
> fsolve(gb[1], y);
                          -1.905166168, -0.6710436067

\end{verbatim}
Similarly, we see that $Q$ and $\denom$ vanish simultaneously
only at $(1,0)$ and at $(0,1)$.  Using our \emph{a priori}
knowledge that the coefficients of $F$ are nonnegative,
we apply Corollary~\ref{cor:procedure} and look for minimal 
points on the lowest arc of the graph of $\sing$ in 
the first quadrant.  The plot of this is visually indistinguishable
from the line segment $x + y=1$, which is not surprising because
the forbidden substrings only affect the terms of the generating
function of order~7 and higher.  More computer algebra shows the
point $(x,y) \in \contrib_{\rbar}$ to be algebraic of degree~21.

One question we might ask next is to what degree the forbidden
substrings bias the typical word to contain an unequal number of 
zeros and ones.  One might imagine there will be a slight preference
for zeros since the forbidden substrings contain mostly ones.  
To find out the composition of the typical word, we apply 
Theorem~\ref{th:WLLN2}.  Setting $y = x$ gives the 
univariate generating function 
$$
f(x) = F(x,x) = \frac{1+x^5+x^6+x^7-x^9}{1-2x+x^5-x^6-x^8-x^9+x^{10}}
   = \sum_n N(n) x^n
$$
for the number $N(n)$ of words of length $n$ with no occurrences of
forbidden substrings.  The root of the denominator of minimum modulus
is $x_0 = 0.505496 \ldots$, whence $N(n) \sim C (1/x_0)^n$.

The direction associated to the point $(x_0 , x_0)$ is given by the
projective point
\begin{equation} 
\label{eq:max rbar}
\rbar = \left ( x \frac{\partial \denom}{\partial x} , 
   y \frac{\partial \denom}{\partial y} \right ) (x_0 , x_0) \; .
\end{equation}
The ratio $\lambda_0 = r/s$ of zeros to ones for this $\rbar$ is a
rational function of $x_0$.  We may evaluate $x_0$ numerically and
plug it into this rational function, but the numerics will be more
accurate if we reduce algebraically first in $\Q [x_0]$. Specifically,
we may first reduce the rational function to a polynomial by inverting
the denominator modulo the minimal polynomial for $x_0$ (Maple's {\tt
gcdex} function) then multiplying by the numerator and reducing again.
Then, from the representation $\lambda_0 = P(x_0)$, we may produce the
minimal polynomial for $\lambda_0$ by writing the powers of
$\lambda_0$ all as polynomials of degree~9 in $x_0$ and determining
the linear relation holding among the powers $\lambda_0^0 , \ldots ,
\lambda_0^{10}$.  We find in the end that $\lambda_0 = 1.059834 \ldots
$. Thus indeed, there is a slight bias toward zeros.

Suppose we wish to know how long a string must be before the count of
zeros and ones tells us whether the string was sampled from the
measure avoiding 10101101 and 1110101 versus the uniform measure. We
may answer this by means of the local central limit behavior described
in Theorem~\ref{th:WLLN2}.  We may verify that the proportion of zeros
is distributed as $$\frac{\lambda_0}{\lambda_0 + 1} + c n^{-1/2}
N(0,1)$$ where $N(0,1)$ denotes a standard normal.  Once we have
computed the constant $c$, we will know how big $n$ must be before
$\displaystyle{\frac{\lambda_0}{\lambda_0 + 1} - \frac{1}{2} \gg c
n^{-1/2}}$ and the count of zeros and ones will tip us off as to which
of the two measures we are seeing.

\setcounter{equation}{0}
\section{Lagrange inversion}
\label{sec:lagrange}

Suppose that a univariate generating function $f(z)$ satisfies
the functional equation $f(z) = z \phi (f(z))$ for some 
function $\phi$ analytic at the origin and not vanishing there.
Such functions often arise, among other places, in graph and tree
enumeration problems.  If $\phi$ is a polynomial, then $f$
is algebraic, but even in this case it may not be possible
to solve for $f$ explicitly.  A better way to analyse $f$
is via the Lagrange inversion formula.  One common
formulation states~\cite[Thm 1.2.4]{GoJa2004} that
\begin{equation} 
\label{eq:lagrange 1}
[z^n] f(z) = \frac{1}{n} \left [ y^{n-1} \right ] \phi(y)^n
\end{equation}
where $[y^n]$ denotes the coefficient of $y^n$.

To evaluate the right side of~(\ref{eq:lagrange 1}), we look at 
the generating function
$$
\frac{1}{1 - x \phi (y)} = \sum_{n=0}^\infty x^n \phi (y)^n
$$
which generates the powers of $\phi$.  The $x^n y^{n-1}$ term
of this is the same as the $y^{n-1}$ term of of $\phi (y)^n$. 
In other words,
\begin{equation} 
\label{eq:lagrange 2}
[z^n] f(z) = \frac{1}{n} \left [ x^n y^n \right ] 
   \frac{y}{1 - x \phi (y)} \, .
\end{equation}
This is a special case of the more general formula for 
$\psi (f(z))$:
\begin{equation} 
\label{eq:lagrange 4}
[z^n] \psi (f(z)) = \frac{1}{n} \left [ x^n y^{n} \right ] 
   \frac{y \psi' (y)}{1 - x \phi (y)} \, .
\end{equation}
These formulae hold at the level of formal power series, and, if $\psi$
and $\phi$ have nonzero radius of convergence, at the level of analytic
functions.

We can now apply the analysis leading to~(\ref{eq:riordan-asymp})
(taking note that the series $v(x)$ in that equation is presently 
called $\phi(y)$, while the series $\phi(x)$ there is here called 
$y \psi' (y)$).  Assume $\phi$ has degree at least~2 -- the
easy case $\phi (z) = az + b$ may be handled directly.  

Recall the definitions of $\mu$ and $\sigma^2$ from 
Subsection~\ref{ss:riordan}.  We are interested only in 
coefficients of $x^n y^n$, that is, the diagonal direction. 
Set $\mu(\phi; y) := y \phi' / \phi$ equal to~$1$; geometrically, we
graph $\phi (y)$ against $y$ and ask that the secant line from the
origin to the point $(y,\phi(y))$ be tangent to the graph there. 
\begin{figure}[ht] \hspace{1.5in}
\includegraphics[scale=0.6]{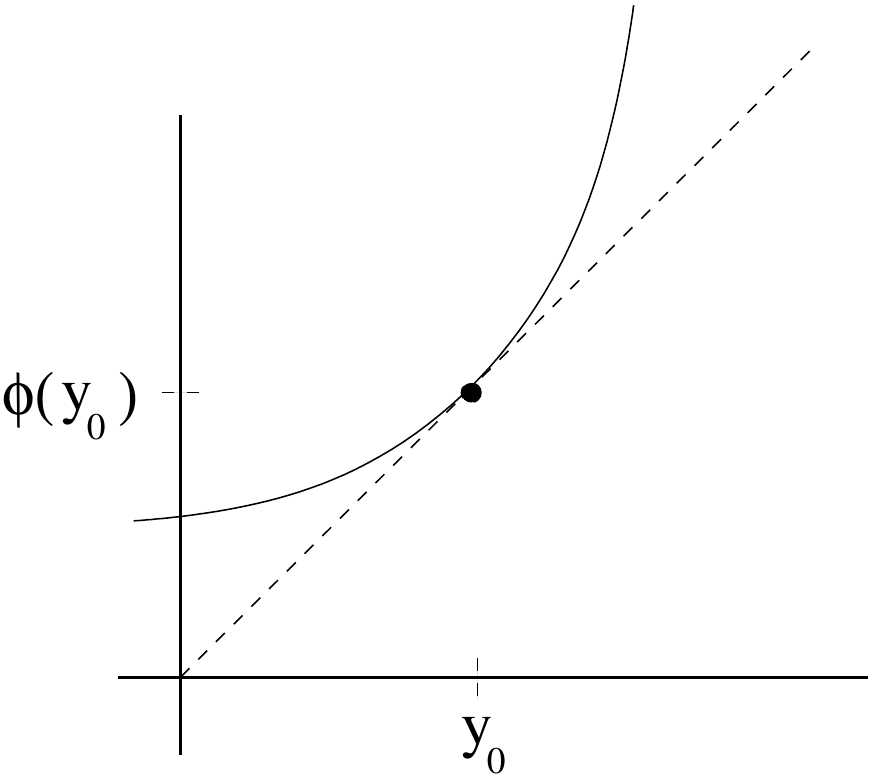}
\label{fig:phi}
\end{figure}
Letting $y_0$ denote a solution to this, we then have a point $(1 /
\phi (y_0) , y_0)$ in the set $\crit_{\one}$ and at this point the
quantity $x^{-n} y^{-n}$ is equal to $(\phi (y_0) / y_0)^n =
\phi'(y_0)^n$.  In equations~(\ref{eq:mu}) and~(\ref{eq:sigma}) we
have $\mu(\phi; y_0) = 1$ and consequently $\sigma^2(\phi; y_0)$
simplifies to $y_0^2 \phi''(y_0) / \phi (y_0)$.  Putting this together
with the asymptotic formula~(\ref{eq:riordan-asymp}), setting $r=s=n$
and simplifying  leads to the following proposition.  Note that it is
easily shown that $f$ is aperiodic if and only if $\phi$ is.

\begin{pr}
\label{pr:lagrange}
Let $\phi$ be analytic and nonvanishing at the origin, aperiodic 
with nonnegative coefficients, with degree at least $2$ at infinity. 
Let $f$ be the nonnegative series satisfying $f(z) = z \phi(f(z))$. 
Let $y = y_0$ be the positive solution of $y \phi' (y) = \phi (y)$. 
Then if $\psi$ has radius of convergence strictly greater than $y_0$, 
we have
\begin{equation}
\label{eq:riordan-univariate}
[z^n] \psi(f(z)) \sim \phi'(y_0)^n n^{-3/2} \sum_{k\geq 0} b_k n^{- k}
\end{equation}
where $\displaystyle{b_0 = \frac{y_0 \psi'(y_0)} 
{\sqrt{2\pi \phi''(y_0)/\phi(y_0)}}\,}$.
$\Cox$
\end{pr}

A variant of Proposition~\ref{pr:lagrange}, proved by other means, is
found in \cite[Thm VI.6]{FlSe}.  That result,
proved by univariate methods, is stronger than
Proposition~\ref{pr:lagrange} in some ways.  For example, it can
handle the estimation of $[z^n] T(z) (1 - T(z))^{-2}$ where $T(z) = z
\exp(T(z))$, which occurs in the study of random mappings, whereas an
attempt to use Proposition~\ref{pr:lagrange} directly runs into the
problem that $y_0 = 1$ and $\psi (z) = (1-z)^{-2}$ has radius of
convergence $1$.

\subsection{Bivariate asymptotics}
\label{ss:lagrange-bivariate}

We now discuss bivariate asymptotics. In the special case $\psi(y) =
y^k$, for fixed $k$, the above results on Lagrange inversion 
yield the first order asymptotic
$$
[z^n] f(z)^k \sim \frac{k}{n} 
   \phi'(y_0)^n \frac{y_0^{k}}{\sqrt{2\pi n\phi''(y_0)/\phi(y_0)}}
$$
where $y_0 \phi'(y_0) = \phi'(y_0)$.
We can also derive an asymptotic as both $n$ and $k$ approach
$\infty$. We sketch an argument here (see \cite{Wils2005} for
details). The formula~(\ref{eq:riordan-asymp}) supplies asymptotics
whenever $n/k$ belongs to a compact set of the interior of $(A, B)$,
where $A$ is the order of $f$ at $0$ and $B$ the order at $\infty$. In
the present case we have $A=1$ and $B=\infty$. We may then use the
defining relation for $f(z)$ to express the formula in terms of $\phi$
only, leading to the following result.

\begin{pr}
\label{pr:lagrange-bivariate}
Let $\phi$ be analytic and nonvanishing at the origin, with
nonnegative coefficients, aperiodic and of order at least $2$ at
infinity. Let $f$ be the positive series satisfying $f(z) = z
\phi(f(z))$. For each $n, k$, set  $\lambda = k/n$ and let $y =
y_\lambda$ be the positive real solution of the equation $\mu(\phi; y)
=  (1 - \lambda)$. Then
\begin{equation}
\label{eq:riordan-bivariate}
[z^n] f(z)^k \sim  \lambda \phi(y_\lambda)^n y_\lambda^{k-n} 
   \frac{1}{\sqrt{2\pi n \sigma^2(\phi; y_\lambda)}}
   = \lambda (1 - \lambda)^{-n} \phi'(y_\lambda)^n  
   \frac{y_\lambda^k }{\sqrt{2\pi n \sigma^2(\phi; y_\lambda)}}.
\end{equation}
Here $\mu$ and $\sigma^2$ are given by equations~(\ref{eq:mu}) 
and~(\ref{eq:sigma}) respectively.  The asymptotic approximation 
holds uniformly provided that $\lambda$ lies in a compact subset of $(0,1)$.
$\Cox$
\end{pr}

The similarity between~(\ref{eq:riordan-univariate})
and~(\ref{eq:riordan-bivariate}) seems to indicate that a version of
Proposition~\ref{pr:lagrange-bivariate} that holds uniformly for $0
\leq k/n \leq 1 - \varepsilon$ should apply.  This is consistent with
a result of Drmota \cite{Drmo1994} but is seemingly
more general, and warrants further study. Such a result can very
likely be obtained using the extension by
Lladser~\cite{Llad2003} of the methods of~\cite{PeWi2002}, as can
results of Meir and Moon \cite{MeMo1990} and
Gardy~\cite{Gard1995} related to
Proposition~\ref{pr:lagrange-bivariate}.

\subsection{Trees in a simple variety}
\label{ss:simple variety}

We now discuss a well-known 
situation~\cite[VII.2]{FlSe} 
in which the foregoing results can be applied.  

Consider the class of ordered (plane) unlabelled trees belonging to 
a so-called \Em{simple variety}, namely a class $\var$ defined by 
the restriction that each node may have a number of children belonging 
to a fixed subset $\Omega$ of $\mathbb{N}$. Some commonly used simple 
varieties are: regular $d$-ary trees, $\Omega = \{0, d\}$; 
unary-binary trees, $\Omega = \{0, 1, 2\}$; general plane trees, 
$\Omega = \mathbb{N}$. The generating function $f(z)$ counting trees 
by nodes satisfies $f(z) = z \omega(f(z))$ where $\omega(y) 
= \sum_{t\in \Omega} y^t$. 

Provided that $\omega$ is aperiodic, Proposition~\ref{pr:lagrange}
applies. The form of $\omega$ shows that the equation $y \omega'(y) =
\omega(y)$ always has a unique solution strictly between $0$ and $1$.
As a simple example, we compute the asymptotics for the number of
general plane trees with $n$ nodes (the exact answer being the Catalan
number $C_{n-1}$).  The equation $y \phi'(y) = \phi(y)$ has solution
$y = 1/2$, corresponding to $\phi'(y) = 1/4$. Thus we obtain 
$$
[z^n] f(z) \sim 4^{n-1} \frac{1}{\sqrt{\pi n^3}}
$$ 
in accordance with Stirling's approximation applied to the 
expression of $C_{n-1}$ in terms of factorials.

Proposition~\ref{pr:lagrange} does not directly apply to the case of
regular $d$-ary trees.  Rather, we use a version where the
contribution from several minimal points on the same torus must be
added. The details are as follows.  The equation $\mu(\phi; y) = 1$
has solutions $\omega \rho$ where $\rho = (d - 1)^{-1/d}$ and
$\omega^d = 1$.  Each of these is a contributing critical point and
the corresponding leading term asymptotic contribution is $C
\omega^{(n - 1)} \alpha^n n^{-3/2} $ where $C = \rho^2 (2\pi(d -
1))^{-1/2}$ and $\alpha^d = d^d/((d-1)^{d-1})$. Thus the asymptotic
leading term is $0$ unless $d$ divides $n - 1$, in which case it is $d
C \alpha^n n^{-3/2}$.

We note in passing that one can avoid the last computation  as
follows~\cite[I.5, Example 13]{FlSe}. There is a bijection
between the class of $d$-ary trees and the class $\mathcal{C}$ of
trees with vertices of degree at most $d$ but allowing for
$\binom{d}{j}$ types of nodes of degree $j$. The pruning map removes
all external nodes (nodes of degree 0) from a $d$-ary tree, then
labelling each node of the pruned tree with the set of children that
were removed.  Pruning always maps a tree of $1 + dm$ nodes to a tree
of $m$ nodes.  The extension of the above asymptotics to a multiset of
degrees (allowing for different types of children) is straightforward.
Thus we can compute using the degree enumerator $g(z) = (1 + z)^d$ of
$\mathcal{C}$ in Proposition~\ref{pr:lagrange}, and $[z^n] f(z) =
[z^m] g(z)$.

Now we consider the \Em{mean degree profile} of trees in $\var$. 
Let $\xi_k(t)$ be the number of nodes of degree $k$ in the tree $t$, 
and $|t|$ the total number of nodes in $t$.  A standard
calculation~\cite[VII.2.2, Example 5]{FlSe} shows that 
the cumulative generating function is 
$$
\sum_{t \in \var} z^{|t|} \xi_k(t) = z^2 \phi_k f(z)^{k - 1} \phi'(f(z))
$$
where $\phi_k = [y^k] \phi(y)$.
Thus we have
$$
F(z, u):= \sum_{k\geq 0} \sum_{t \in \var} \xi_k(t) z^{|t|} u^k 
   = \mu(f; z) z \phi(u f(z)) \, .  
$$
The mean number of nodes of degree $k$ in a uniformly randomly chosen 
tree of size $n$ from $\var$ is then given by
$$
M_{nk} = \frac{[z^n u^k] F(z, u)}{[z^n] f(z)}.
$$
Consider again the simplest case, general plane trees, with $\phi(y) =
(1 - y)^{-1}$.  Then $F(z, u)$ corresponds to a Riordan array.  
A simple variant of Proposition~\ref{pr:lagrange-bivariate} shows that 
$$
[z^n u^k] F(z, u) \sim y^{k} \phi'(y)^{n - 1}  
   \frac{1}{\sqrt{2 \pi n \sigma^2(\phi; y)}}
$$
where $\mu(\phi; y) = 1 - k/n$. This is easily solved to obtain 
$y = (n - k)/(2n - k)$ and hence we obtain routinely:
\begin{equation} 
\label{eq:tree-level}
M_{nk} \sim n \left(\frac{2n-k}{2n}\right)^{2n-2} 
   \left(\frac{n-k}{2n-k}\right)^k \sqrt{\frac{n^2}{2 (n - k) (2n - k)}}
\end{equation}
uniformly as long as $k/n$ is in a compact subset of $(0, 1)$.
The mean number of leaves (nodes of degree $0$) in such a tree is 
well known~\cite[III.5, Example 12]{FlSe} to be $n/2$, 
which is obtained by substituting $k = 0$ in the right side 
of~(\ref{eq:tree-level}).  Thus it again seems likely that the approximation  
is in fact uniform on $[0, 1 - \varepsilon]$.

\subsection{Cores of planar graphs}
\label{ss:core}

A \Em{rooted planar map} 
is a graph with a distinguished edge (the root) that can be 
embedded in the plane.  A rooted planar map is completely 
specified by a planar graph, a distinguished edge, and a cyclic 
ordering of edges around each vertex.  
The \Em{core} of a map (henceforth always a rooted planar map) is the 
largest $2$-connected subgraph containing the root edge.  The problem
of the typical core size (cardinality of the core) of a map is
considered in~\cite{GaWo1999,BFSS2001}.  They obtain a 
functional equation for the generating function 
$$
M(u,z) = \sum_{n,k} a_{n,k} z^n u^k
$$
where $a_{n,k}$ is the number of maps with $n$ edges and core size $k$
(that is, the core has cardinality $k$).  They are interested in
computing the probability distribution of the core size of a map
sampled uniformly from among all maps of with $n$ edges.  Applying
the general inversion formula~(\ref{eq:lagrange 4}) to their functional
equations, they are able to compute
\begin{equation} \label{eq:BFSS2001}
p(n,k) := \frac{a_{n,k}}{\sum_j a_{n,j}} = \frac{k}{n}
   [z^{n-1}] \psi' (z) \psi(z)^{k-1} \phi(z)^n
\end{equation}
where $\psi (z) = (z/3) (1 - z/3)^2$ and $\phi (z) = 3 (1+z)^2$
are the operators that arise in the functional equation for $M(u,z)$.

The mean size of the core was computed in~\cite{GaWo1999}.
For this it sufficed to prove a specialized result about coefficients
of powers of functions asymptotic to $(1-z)^{-3/2}$ as $z \to 1$.  
The problem of finding a limit law for the distribution about the mean 
was taken up in~\cite{BFSS2001}.  
One may arrive directly at an asymptotic formula for $p(n,k)$ 
if one rewrites~\eqref{eq:BFSS2001} as
$$
p(n,k) = \frac{k}{n} [x^k y^n z^n] \frac{x z \psi'(z)}
   {(1 - x \psi (z)) (1 - y \phi (z))} \, .
$$
The analysis of the resulting trivariate generating function is,
however, quite challenging.  In particular, there is a point 
where $Q$ vanishes; a generalization of 
Theorem~\ref{th:smooth asym}~\cite[Theorem~3.3]{PeWi2002}
tells us that the asymptotics in this precise direction, but does not
answer the more interesting question of asymptotics in a scaling
window near that direction, which turns out to be $k = n/3 \pm
n^{-2/3}$.  An answer, in the framework of~\cite{PeWi2002}, is given 
in the  doctoral work of 
M.\ Lladser~\cite[Section~5.5]{Llad2003}. 
A complete answer is given in~\cite{BFSS2001} by
reductions to a one-variable generating function, 
on which a coalescing-saddle approximation is used. 
We will not go into details here.

\setcounter{equation}{0}
\section{The kernel method}
\label{sec:kernel}

The kernel method is a means of producing a generating function for
an array $\{ a_\rr : \rr \in \N^d \}$ of numbers satisfying a
linear recurrence
\begin{equation} \label{eq:recursion}
a_\rr = \sum_{\sss \in E} c_\sss a_{\rr - \sss} \, .
\end{equation}
Here $E$ is a finite subset of $\Z^d \setminus \{ 0 \}$ which is 
not necessarily a subset of $\N^d$ but whose convex hull must 
not intersect the closed negative orthant. 
The numbers $\{ c_\sss : s \in E \}$ are constants
and the relation~(\ref{eq:recursion}) holds for all $\rr$ except those
in the \Em{boundary condition}, which will be made precise below.
As usual, let $F(\zz) = \sum a_\rr \zz^\rr$.
In one variable $F$ is always a rational function, 
but in more than one variable the generating function
can be rational, algebraic, D-finite, or differentially transcendental
(not D-finite).  A classification along these lines, determined more
or less by the number of coordinates in which points of $E$ can be
negative, is given in~\cite{BoPe1998}, which is a
very nice exposition of the kernel method at an elementary level.

We are interested in the kernel method because it often produces
generating functions to which Theorem~\ref{th:smooth asym} may be
applied.  Because the method is not all that well known, we include a
detailed description, drawing heavily
on~\cite{BoPe1998}.  We begin, though, with an example.

\subsection{A random walk problem}
\label{ss:random walk}
Two players move their tokens toward the finish square, 
flipping a fair coin each time to see who moves forward one square.  
At present the distances to the finish are $1+r$ and $1+r+s$.  If the 
second player passes the first player, the second player wins; if the 
first player reaches the finish square, the first player wins; 
if both players are on the square before the finish square, it is 
a draw.  What is the probability of a draw?  

Let $a_{rs}$ be the probability of a draw, starting with initial
positions $1+r$ and $1+r+s$.  Conditioning on which player
moves first, one find the recursion 
$$
a_{rs} = \frac{a_{r,s-1} + a_{r-1,s+1}}{2}
$$ 
which is valid for all $(r,s) \geq (0,0)$ except for $(0,0)$, provided
that we define $a_{rs}$ to be zero when one or more coordinate is
negative. The relation $a_{rs} - (1/2) a_{r,s-1} - (1/2) a_{r-1,s+1} =
0$ suggests we multiply the generating function $F(x,y) := \sum a_{rs}
x^r y^s$ by $1 - (1/2) y - (1/2) (x/y)$.  To clear denominators, we
multiply by $2y$: define $Q(x,y) = 2y - y^2 - x$ and compute $Q \cdot
F$.  We see that the coefficients of this vanish with two exceptions:
the $x^0 y^1$ coefficient corresponds is $2 a_{0,0} - a_{0,-1} -
a_{-1,1}$ which is equal to~2, not~0, because the recursion does not
hold at $(0,0)$ ($a_{00}$ is set equal to~1); the $y^0 x^j$
coefficients do not vanish for $j \geq 1$ because, due to clearing the
denominator, these correspond to $2 a_{j,-1} - a_{j,-2} -  a_{j-1,0}$.
This expression is nonzero since, by definition, only the third term
is nonzero, but the value of the expression is not given by prescribed
boundary conditions. That is, we have
\begin{equation} \label{eq:kernel 1}
Q(x,y) F(x,y) = 2 y - h(x)
\end{equation}
where $h(x) = \sum_{j \geq 1} a_{j-1,0} x^j = x F(x,0)$ will not be 
known until we solve for $F$.  

This generating function is in fact a simpler variant of the one
derived in~\cite{LaLy1999} for the waiting time 
until the two players collide, which is needed in the analysis of
a sorting algorithm.  Their solution is to observe that there is
an analytic curve in a neighborhood of the origin on which $Q$ 
vanishes.  Solving $Q=0$ for $y$ in fact yields two solutions,
one of which, $y = \xi (x) := 1 - \sqrt{1-x}$, vanishes at the
origin.  Since $\xi$ has a positive radius of convergence, we
have, at the level of formal power series, that $Q(x , \xi (x)) = 0$,
and substituting $\xi (x)$ for $y$ in~(\ref{eq:kernel 1}) gives
$$
0 = Q(x,\xi(x)) F(x,\xi (x)) = 2 \xi (x) - h(x) \, .
$$
Thus $h(x) = 2 \xi (x)$ and 
$$
F(x,y) = 2 \frac{y - \xi (x)}{Q(x,y)} = \frac{1}{1 + \sqrt{1-x} - y}
\, .
$$

As is typical of the kernel method, the generating function $F$ has a
pole along the branch of the kernel variety $\{ Q(x,y) = 0 \}$ that
does not pass through the origin.  The function $F$ is not meromorphic
everywhere, having a branch singularity on the line $x=1$, but it is
meromorphic in neighborhoods of polydisks $\disk (x,y)$ for minimal
points $(x,1 + \sqrt{1-x})$ on the pole variety for $0 < x < 1$. For
$0 < x < 1$, $\dir (x , 1 + \sqrt{1-x}) = (2 \sqrt{1-x} + 2 - 2x)/x$.
If we set this equal to $\lambda$ and solve for $x$ we find $x = 4
(1+\lambda)/(2+\lambda)^2$ and $1+\sqrt{1-x} = (2+2\lambda) /
(2+\lambda)$.  In other words, 
$$
\contrib_{\rbar} = \left ( \frac{4(1+\lambda)}{(2+\lambda)^2} , 
\frac{2 (1+\lambda)}{2+\lambda} \right )
$$
where $\lambda = s/r$.  Plugging this into~(\ref{eq:smooth asym})
gives
\begin{eqnarray*}
a_{rs} & \sim & C (r+s)^{-1/2} \left ( \frac{4 r (r+s)}{(2r + s)^2} 
   \right )^{-r} \left ( \frac{2(r+s)}{2r+s} \right )^{-s} \\[2ex]
& = & \frac{C}{2^{2r+s}} \frac{(2r+s)^{2r+s}}{r^r (r+s)^{(r+s)}} \, .
\end{eqnarray*}
One recognizes in this formula the asymptotics of the binomial
coefficient $\binom{2r+s}{r}$ and indeed the binomial coefficient may
be obtained via a combinatorial analysis of the random walk paths.

\subsection{Explanation of the kernel method}
Because of the applicability of Theorem~\ref{th:smooth asym} to 
generating functions derived from the kernel method, we now give 
a short explanation of this method.  We adopt the notation from 
the first paragraph of this section.

Let $\pp$ be the coordinatewise infimum of points in $E \cup \{ \zero
\}$, that is the greatest element of $\Z^d$ such that $\pp \leq \sss$
for every $\sss \in E \cup \{ \zero \}$.  Let $$Q(\zz) :=
\zz^{-\pp}\left(1 - \sum_{\sss \in E} c_\sss \zz^{\sss}\right),$$
where the normalization by $\zz^{-\pp}$ guarantees that $Q$ is a
polynomial but not divisible by any $z_j$.  Partition $\N^d$ into two
sets, $Z$ and $B$, and assume that the relation~(\ref{eq:recursion})
holds for all $\rr \in Z$; the set $B$ must be closed under $\leq$ and
the values $\{ a_\rr : \rr \in B \}$ are specified explicitly by
constants $\{ b_\rr : \rr \in B \}$ rather than
by~(\ref{eq:recursion}).

If $E \subseteq \N^d$ then $\pp = \zero$, $B$ can be arbitrary, $F$ is
a function of the form $\numer / Q$, and $\numer$ is rational if the
boundary conditions are rational.  The analysis in this case is
straightforward and the kernel method yields only what may be derived
directly from the recursion for $a_\rr$ in terms of $\{ a_\sss : \sss
\leq \rr \}$.  For examples of this, see Sections~\ref{ss:binomial}
and~\ref{ss:delannoy}. We concentrate instead on the case where $d=2$
and the second coordinate of points in $E$ may be negative.  It is
known in this case~\cite[Theorem~13]{BoPe1998} that if
the generating function for the boundary conditions is algebraic then
$F$ is algebraic.  On the other hand, we shall see that an outcome of
the kernel method is that $F$ will have a pole variety, and will
usually satisfy the meromorphicity condition in the remark after
Theorem~\ref{th:universal}.

To apply the kernel method, one examines the product $Q F_Z$, 
where for convenience we have let $F_Z := \sum_{\rr \in Z} 
a_\rr \zz^\rr$ be the generating function for those
values for which the recursion~(\ref{eq:recursion}) holds.
If we assume that the generating function 
$F_B := \sum_{\rr \in B} b_\rr \zz^\rr$ for the prescribed
boundary conditions is known, then since $F = F_Z + F_B$,
finding $F_Z$ is equivalent to finding $F$.

There are two kinds of contribution to $Q F_Z$.  Firstly, for every 
pair $(\rr , \sss)$ with $\sss \in E$, $\rr \in Z$ and $\rr - \sss 
\in B$, there is a term $c_\sss b_{\rr - \sss} \zz^{\rr -\pp}$:
the corresponding coefficient of $Q F$ vanishes, but $F_Z$ has been 
stripped of the boundary terms that cause the cancellation.  Let 
$$
K(\zz) := \sum_{\rr \in Z, \sss \in E, \rr - \sss \in B} c_\sss 
b_{\rr - \sss} \zz^{\rr - \pp}
$$ 
denote the sum of these terms.
The ``$K$'' stands for ``known'', because the coefficients of
$K$ are determined by the boundary conditions, which are known.
Secondly, for every pair $(\rr , \sss)$ with $\sss \in E, 
\rr - \sss \in B$ and $\rr \notin Z$ there is a term 
$-c_\sss a_{\rr - \sss} \zz^{\rr -\pp}$ due to the fact that
the recursion does not hold at $\rr$.
$$
U(\zz) = \sum_{\rr - \sss \in Z , \sss \in E , \rr \notin Z}
c_\sss a_{\rr - \sss} \zz^{\rr - \pp}
$$ 
denote these terms.  The ``$U$'' stands
for ``unknown'', because these coefficients are not
explicitly determined from the boundary conditions.
It is not hard to show that for any dimension $d$ and any
$E$ whose convex hull does not contain a neighborhood of the origin, 
the following result holds.
\begin{pr}[\protect{\cite[Theorem~5]{BoPe1998}}] 
\label{pr:BMP}
Let $E$, $\{ c_\sss \}$, $\pp$, $B$, $Z$, $\{ b_\sss : \sss \in B \}$,
$F_Z$ and $K$ be as above.  Then there is a unique set of values 
$\{a_\rr : \rr \in Z \}$ for such that~(\ref{eq:recursion}) holds 
for all $\rr \in Z$.  Consequently, there is a unique pair $(F, U)$ of
formal power series  such that
$$
Q F_Z = K - U \, .
$$
$\Cox$
\end{pr}

In other words, the unknown power series $U$ is determined from
the data along with $F_Z$.  Another way of thinking about this is
that $F_Z$ is trying to be the power series $K/Q$ but since $Q$ 
vanishes at the origin, one must subtract some terms from $K$ to
cancel the function $Q_0$ where $Q = Q_0 Q_1$ and $Q_0$ consists
of the branches of $Q$ passing through the origin.  The kernel
method, as presented in~\cite{BoPe1998} turns
this intuition into a precise statement.
\begin{pr}[\protect{\cite[equation~(24)]{BoPe1998}}] 
\label{pr:branches}
Suppose $d=2$, $\pp = (0 , -p)$ and $F_B = 1$.  There will be exactly
$p$ formal power series $\xi_1 , \ldots , \xi_p$ such that $\xi_j (0) =
0$ and $Q(x,\xi_j (x)) = 0$, and we may write $Q(x,y) = -C(x)
\prod_{j=1}^p (y - \xi_j (x)) \prod_{j=1}^{r} (y - \rho_j (x))$ for some
$r$ and $\rho_1 , \ldots , \rho_r$.  The generating function $F_Z$ will
then be given by
$$
F_Z (x,y) = \frac{K(x,y) - U(x,y)}{Q(x,y)} 
   = \frac{\prod_{j=1}^p (y - \xi_j (x))}{Q(x,y)} 
   = \frac{1}{-C(x) \prod_{j=1}^r (y - \rho_j(x))} \, .
$$
$\Cox$
\end{pr}
We turn to some examples.

\subsection{Dyck, Motzkin, Schr\"{o}der, and generalized Dyck paths}
\label{ss:paths}

Let $E$ be a set $\{ (r_1 , s_1) , \ldots , (r_k , s_k) \}$ of integer
vectors with $r_j > 0$ for all $j$ and $\min_j s_j = -p < 0 < \max_j
s_j = P$. The generalized Dyck paths with increments in $E$ to the
point $(r,s)$ in the first quadrant are the paths from $(0,0)$ to
$(r,s)$, with increments in $E$, which never go below the horizontal
axis.

\begin{figure}[ht] \hspace{1.5in}
\includegraphics[scale=0.6]{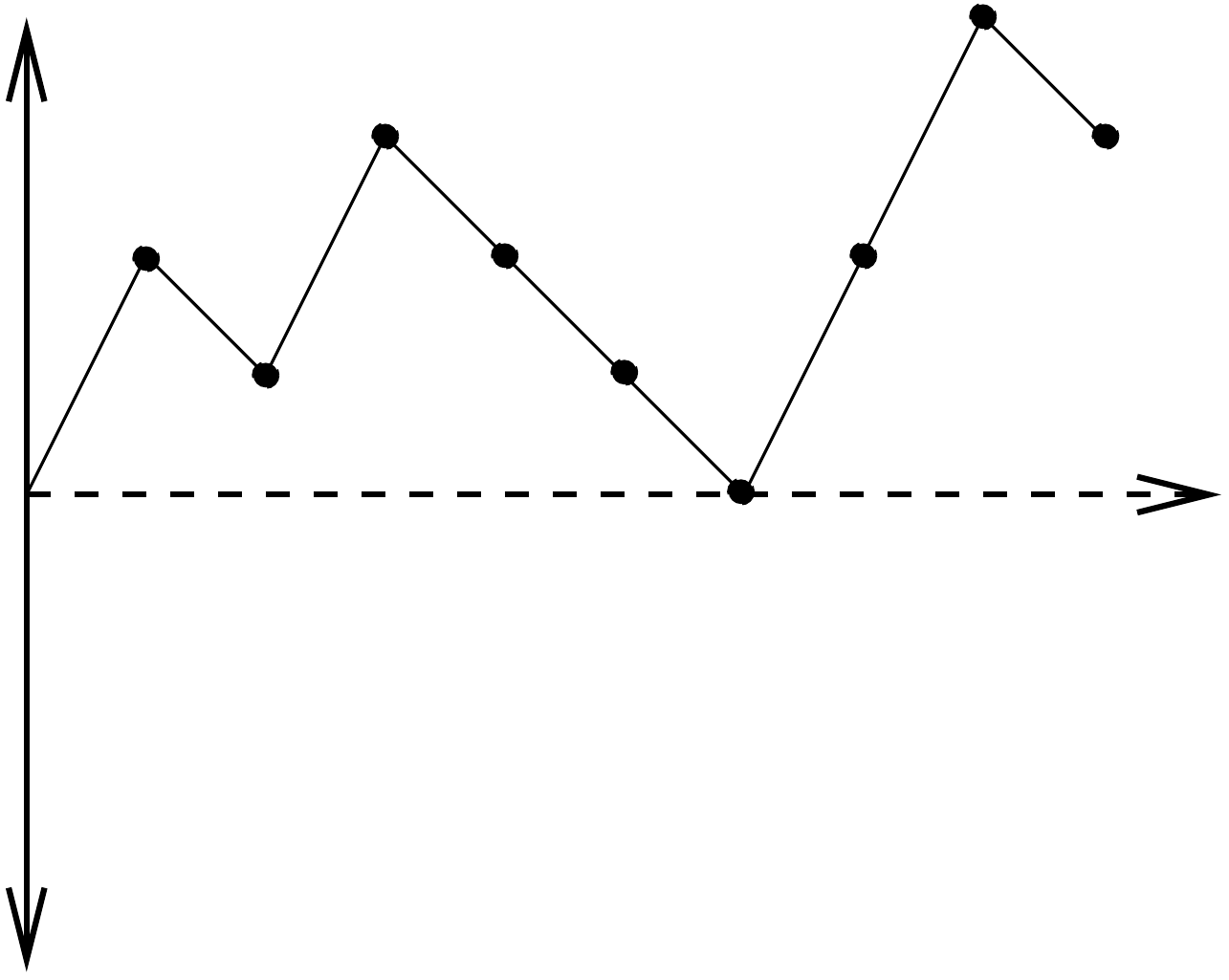}
\caption{a generalized Dyck path of length nine with $E = \{
(1,2),(1,-1) \}$}
\label{fig:gen-dyck}
\end{figure}

Let $F(x,y) = \sum_{r,s} a_{rs} x^r y^s$ where $a_{rs}$ is the number
of generalized Dyck paths to the point $(r,s)$.  We will have
$B = \{ (0, 0) \}$, $F_B = 1$, $Q(x, y) = y^p ( 1 - \sum_i 
x^{r_i} y^{s_i})$, and $C(x) =  \sum_{i: s_i = P} x^{r_i}$.  
We now discuss three well known instances of such paths 
taken from~\cite{BoPe1998}.

\setlength{\unitlength}{1.5pt}
\begin{picture}(170,80)(0,5)
\put(30,60){Dyck paths}
\put(30,30){\circle*{2}}
\put(30,30){\vector(1,1){20}}
\put(30,30){\vector(1,-1){20}}
\put(130,60){Motzkin paths}
\put(130,30){\circle*{2}}
\put(130,30){\vector(1,1){20}}
\put(130,30){\vector(1,-1){20}}
\put(130,30){\vector(1,0){20}}
\put(230,60){Schr\"{o}der paths}
\put(230,30){\circle*{2}}
\put(230,30){\vector(1,1){20}}
\put(230,30){\vector(1,0){40}}
\put(230,30){\vector(1,-1){20}}
\end{picture}
\begin{figure}[ht] \hspace{1.5in}
\caption{legal steps for three types of paths}
\label{fig:Dyck}
\end{figure}

\noindent{\bf Dyck paths:}
When $E = \{ (1,1) , (1,-1) \}$ we have the original Dyck paths.
We have $Q(x,y) = y - xy^2 - x$. Here $C(x) = x$, and  
$Q(x,y) = -x (y - \xi (x)) (y - \rho (x))$ where 
$\xi (x) = (1 - \sqrt{1 - 4 x^2})/(2x)$ and 
$\rho (x) = (1 + \sqrt{1 - 4 x^2})/(2x)$ is the 
algebraic conjugate of $\xi$. Note that $\rho$ is a formal 
Laurent series and $\rho \xi = 1$.

Thus we have, following the discussion above, 
$$
F(x,y) = \frac{1}{-x (y - \rho(x))} = \frac{\xi(x)/x}{1 - y\xi(x)} .
$$ 
Setting $y = 0$ recovers the fact that the Dyck paths coming back to
the $x$-axis at $(2n,0)$ are counted by the Catalan number $C_n$.

Asymptotics are readily obtained either using the explicit or implicit
form of $\xi$ (noting the periodicity of $\xi$).  Let us use the
implicit form in this example, since we will be illustrating use of
the explicit form below in the case of Schr\"{o}der paths, where
Lagrange inversion does not apply.  The vanishing of $Q (x,y) = y -
xy^2 - x$ occurs when $y = x \phi (y)$ with $\phi (y) = 1 + y^2$.
Thus $\xi = x \phi (\xi)$.  We use a slight variant
of~(\ref{eq:riordan-bivariate}), namely
$$
[z^n] \xi(z)^{k+1}/z \sim \lambda \phi(y_\lambda)^{n+1}
y_\lambda^{k-n} \frac{1}{\sqrt{2\pi n \sigma^2(\phi; y_\lambda)}}
$$
where $\lambda = k/n$ and $\mu(\phi; y_\lambda) = 1 - \lambda$. 

The chain rule yields $\mu(\phi; y) = y \phi_y/\phi = y \phi_t t_y/\phi 
= (y t_y/t) (t\phi_t/\phi) = \mu(y^2; y) \mu(1+t; t)$ with
$t = y^2$. Thus we solve $1 - \lambda = 2t/(1 + t)$, or $y^2 = (1 -
\lambda)/(1 + \lambda) = (r - s)/(r + s)$. The two contributing points
$y_\lambda$ cancel out if $r-s$ is odd and reinforce if $r-s$ is even.
We compute $\phi(y_\lambda) = 2/(1+\lambda) = 2r/(r+s)$ and
$\sigma^2(\phi; y_\lambda) = 1 - \lambda^2 = (r^2 - s^2)/r^2$. We
obtain 
\begin{eqnarray*}
a_{rs} & \sim & 2\, \frac{s}{r} \left(\frac{2r}{r+s}\right)^{r+1} 
   \left(\frac{r-s}{r+s}\right)^{\frac{s-r}{2}} 
   \frac{\sqrt{r}}{\sqrt{2\pi (r - s) (r + s)}} \\
& = & \frac{2s}{(r+s)} \frac{r^r}{(\frac{r+s}{2})^{\frac{r+s}{2}} 
   (\frac{r-s}{2})^{\frac{r-s}{2}}} 
   \frac{\sqrt{r}}{\sqrt{2\pi (\frac{r - s}{2}) (\frac{r + s}{2})}} 
\end{eqnarray*}
provided $r-s$ is even, and $0$ otherwise. This is uniform for 
$0 < \delta \leq s/r \leq 1 - \varepsilon < 1$.

Of course, in this simple example Lagrange inversion also gives an
exact formula involving binomial coefficients, namely 
\begin{equation}
\label{eq:dyck-exact}
a_{rs} = [x^{r+1}] \xi(x)^{s+1} = \frac{s+1}{r+1} [y^{r-s}]  
   (1 + y^2)^{r+1} = \frac{s+1}{r+1} \binom{r+1}{(r-s)/2} 
   = \frac{2(s+1)}{r+s+2} \frac{r!}{(\frac{r-s}{2})! (\frac{r+s}{2})!} 
\end{equation}
when $r-s$ is even, and $0$ otherwise. If we had instead computed
$[z^{r+1}] \xi(z)^{s+1}$ using \eqref{eq:riordan-bivariate}, we would
have obtained a correct leading order asymptotic of a slightly
different form ($r, s$ replaced by $r+1, s+1$ in some places). In each
case these asymptotics are consistent with what is obtained by
applying Stirling's formula to the factorials in
\eqref{eq:dyck-exact}.

\noindent{\bf Motzkin paths:}
Let $E = \{ (1,1) , (1,0) , (1,-1) \}$.  In this case the generalized 
Dyck paths are known as \Em{Motzkin paths}.  Again we have
case $Q(x,y) = y - x y^2 - x - xy$.  Now $\rho$ and $\xi$ are given by 
$(1 - x \pm \sqrt{1 - 2x - 3x^2}) / (2x)$ and again
$$
F(x,y) =  \frac{\xi(x)/x}{1 - y\xi(x)} = \frac{2}
   {1 - x + \sqrt{1 - 2x - 3x^2} - 2xy} \, .
$$
This time $\xi$ is given implicitly by $\xi = x(1 + \xi + \xi^2)$ and
the coefficients are not binomial coefficients, but the asymptotics
are no harder to compute. Here, with $\lambda = s/r$, we have that
$\contrib_{1 - \lambda}$ is a singleton $\{y_\lambda\}$ by
aperiodicity. The critical point equation is 
$$
(1 + \lambda) y^2 + \lambda y - (1 - \lambda) = 0.
$$
The solution is
$$ 
y_\lambda = \frac{\sqrt{4 - 3 \lambda^2} - \lambda}{2(1 + \lambda)} \,
.
$$
The minimal polynomial for $\sigma^2(\phi; y_\lambda)$ is found as in
earlier sections to be
$$
3S^2 + (6 \lambda^2 + 12 \lambda - 2) S + 3 \lambda^4 - 24 \lambda^3 +
65 \lambda^2 - 68 \lambda + 24.
$$
This polynomial has two positive solutions for $S$ for each given
$\lambda$. The correct one is found by noting that $\sigma^2$
approaches $0$ as $\lambda$ approaches $0$.

\noindent{\bf Schr\"{o}der paths:} 
Here $E = \{ (1,1) , (2,0) , (1,-1) \}$. We have $C(x) = x, Q(x, y) =
y - xy^2 - x^2 y - x$, and $\rho$ and $\xi$ are given by $(1 - x^2 \pm
\sqrt{1 - 6x^2 + x^4})/(2x)$.  This time Lagrange inversion does not
obviously apply.  We perform an explicit computation, noting the
periodicity of $F$.  We have $\dir(x,y) = \rbar$ with 
$$
\lambda:=\frac{s}{r} = \frac{\sqrt{1 - 6x^2 + x^4}}{1+x^2}.
$$ 
Then $\lambda$ decreases from $1$ to $0$ as $x$ increases from $0$ to
the smaller positive root of $1 - 6x^2 + x^4$, namely $\sqrt{2} - 1$.
We also have 
$$
x^2 = \frac{3+\lambda^2 - 2\sqrt{(2 + 2\lambda^2)}}{1 - \lambda^2}.
$$
Choosing the positive value of $x$, we see that asymptotics are given
by 
$$
a_{rs} \sim 2 C x^{-r} y^{-s} s^{-1/2}
$$
where $y = 1/\xi(x)$, when $r+s$ is even, and $0$ otherwise. Any
particular diagonal (with a value of $\lambda$ between $0$ and $1$)
can be extracted easily. For example, with $\lambda = 1/3$, we obtain
$a_{3s,s} \sim 2C \gamma^s s^{-1/2}$, with $x = (3 - \sqrt{5})/2, y =
(1 + \sqrt{5})/2$, $\gamma = (11 + 5 \sqrt{5})/2 \approx 11.09016992$,
and $C \approx 0.1526195310$. For $s = 12$ this approximation
underestimates by about $2.9 \%$, and for $s = 24$ by about $1.5 \%$.

\subsection{Pebble configurations}
\label{ss:pebble}

Chung, Graham, Morrison and Odlyzlo~\cite{CGMO1995} consider the
following problem.  Pebbles are placed on the nonnegative integer
points of the plane.  The pebble at $(i,j)$ may be replaced by two
pebbles, one at $(i+1,j)$ and one at $(i,j+1)$, provided this does not
cause two pebbles to occupy the same point.  Starting from a single
pebble at the origin, it is known to be impossible to move all pebbles
to infinity; in fact it is impossible to clear the region $\{ (i,j) :
1 \leq i + j \leq 2 \}$~\cite[Lemma~2]{CGMO1995}.

They consider the problem of enumerating minimal unavoidable
configurations. More specifically, say that a set $T$ is a \Em{minimal
unavoidable configuration} with respect to some starting configuration
$S$ if it is impossible starting from $S$ to move all pebbles off of
$T$, but pebbles may be cleared from any proper subset of $T$. Let
$S_t$ denote the starting configuration where $(i,j)$ is occupied if
and only if $i+j = t$.  Let $f_t(k)$ denote the number of sets in the
region $\{ (i,j) : i , j \geq 0 \, ; i+j \geq t+1 \}$ that are minimal
unavoidable configurations for the starting configuration $S_t$.

They derive the recurrence
$$
f_t (k) = f_{t-1} (k) + 2 f_t (k-1) + f_{t+1} (k-2) 
$$
which holds whenever $t \geq 3$ and $k \geq 2$.  Let 
$$
F(x,y) = \sum_{t,k \geq 0} f_t(k) x^k y^t \, .
$$
According to the kernel method, we will have $F = \eta / Q$ for some
$\eta$ vanishing on the zero set of $Q$ near the origin, where $Q(x,y)
= x - (x+y)^2$.  Using some more identities, Chung {\it et al.} are
able to evaluate $F$ explicitly.  They state that they are primarily
interested in $f_0(k)$, so they specialize to $F(x,0)$ and compute the
univariate asymptotics.  It seems to us that the values $f_t (k)$ are
of comparable interest, and we pursue asymptotics of the full
generating function.

The formula for $F$ is cumbersome, but its principal features are (i)
a denominator of $P \cdot Q$ where $P$ is the univariate polynomial $1
- 7x + 14 x^2 - 9 x^3$, and (ii) $F$ is algebraic of degree~2 and is
in $\CC[x,y][\sqrt{1-4x}]$. The minimum modulus root of $P$ is $x_0
\approx 0.2410859 \ldots$.  The algebraic singularity of $F$ occurs
along the line $x = 1/4$, so conveniently, the branching is completely
outside the closure of the domain of convergence and the
meromorphicity assumption of the remark following
Theorem~\ref{th:universal} is satisfied.

\begin{figure}[ht] \hspace{1.5in}
\includegraphics[scale=0.6]{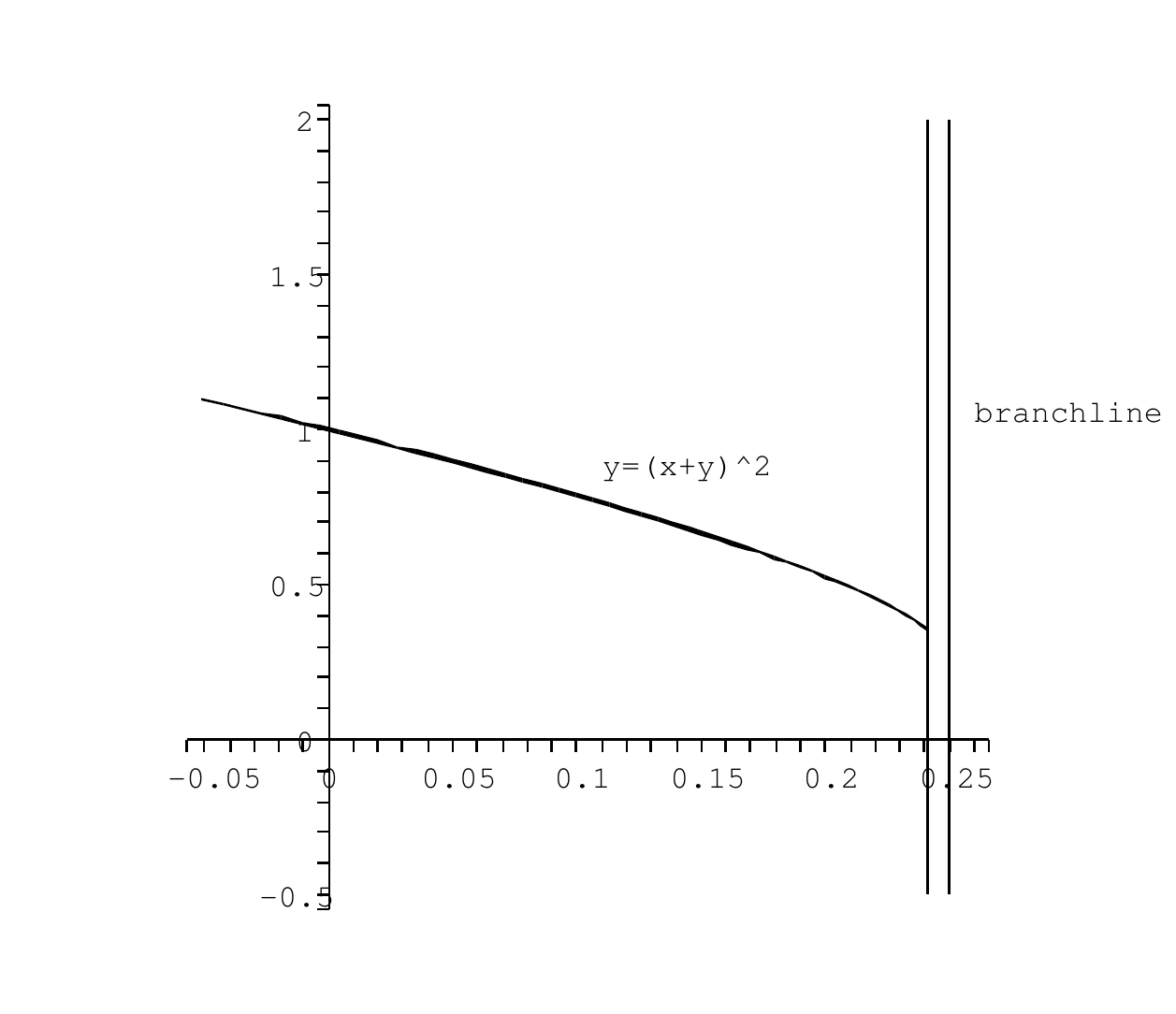}
\caption{the domain of convergence and the algebraic singularity}
\label{fig:two-curves} 
\end{figure}

The boundary of the domain of convergence in the first quadrant is
composed of pieces of two curves, namely $x = x_0$ and $y = (x+y)^2$.
These intersect at the point $(x_0 , y_0)$ where $y_0 = (1 - 2x +
\sqrt{1-4x})/2 \approx 0.3533286$.  We are in the combinatorial case,
so we know that minimal points will be found along these curves. Along
the curve $y = (x+y)^2$ the direction $\dir(x,y)$ is given by $\lambda
= r/s = (2 - 2 \sqrt{y}) / (2 \sqrt{y} - 1)$.  As $(x,y)$ travels from
$(0,1)$ down the curve to $(x_0 , y_0)$, $\lambda$ increases from 0 to
$\lambda_0 \approx 4.295798\ldots$. At the point $(x_0 , y_0)$, the
cone $\cone (x_0 , y_0)$ is the convex hull of the positive
$x$-direction the direction $\overline{(\lambda_0,1)}$.

We then have two sorts of asymptotics for $f_t (k)$.  When $t < k /
\lambda_0$, the asymptotics are given by Corollary~\ref{cor:m=d=2}. In
this case we may evaluate $\numer (x_0 , y_0) \approx 0.00154376$ and
$\sqrt{- x_0^2 y_0^2 \hess} \approx 0.02925688$ so that
$$
f_t(k) \sim C x_0^{-k} y_0^{-t} \;\; \mbox{ with } C = 0.05276 \ldots
$$
uniformly as $t/k$ varies over compact subsets of $(0 , 1 /
\lambda_0)$. It is interesting to compare to the asymptotics for $f_0
(k)$.  Setting $y=0$ gives the univariate generating
function~\cite{CGMO1995} 
$$
f(x) = \sum a_n x^n = x^2 \frac{(1-4x)^{1/2} (1-3x+x^2) - 1 + 5x - x^2
- 6 x^3}{P(x)} \, .
$$
We may compute 
$$\lim_{n \to \infty} a_n x_0^k = \lim_{x \to x_0} (x_0 - x) f(x)$$
and we find that
$$
f_0 (k) \sim c x_0^{-k} \;\; \mbox{ with } c = 0.016762 \ldots \, .
$$
In fact, one may calculate $\lim f_t (k) x_0^k$ for any fixed $t$
by computing 
\begin{equation} \label{eq:t}
f_t (x) := (t!)^{-1} \left ( \frac{\partial}{\partial y} \right )^t F(x,0)
\end{equation}
and again computing $c(t) := \lim_{x \to x_0} (x_0 - x) f_t (x)$. The
pole at $x = x_0$ in $F$ is removable: replacing the denominator $PQ$
of $F$ by $(P/(x-x_0)) Q$ we see that $g(y) := \lim_{x \to x_0} (x_0 -
x) F(x,y)$ has a simple pole at $y=y_0$ and that $\lim_{t \to \infty}
y_0^t c(t) = \lim_{y \to y_0} (y_0 - y) g(y)$ is equal to $C$.  In
other words, we see that the asymptotics known to hold uniformly as
$t/k$ varies over compact subsets of $(0 , 1/\lambda_0)$ actually hold
over $[0 , 1/\lambda_0)$ as long as $t \to \infty$, while for $k \to
\infty$ with $t$ fixed we use~(\ref{eq:t}).

On the other hand, when $t/k > 1 / \lambda_0$ we may solve for $x$ 
and $y$ to get 
$$
\contrib_{\rbar} = (x,y) := \left \{ \left ( \frac{k (2t + k)}
   {(2t+2r)^2} \, , \, \frac{(2t+k)^2}{(2t + 2k)^2}
   \right ) \right \} \, .
$$
We now use Theorem~\ref{th:smooth asym} to see that
$$
f_t (k) C(\frac{t}{k}) t^{-1/2} \sim \frac{(2t+2k)^{2t+2k}}
   {k^k (2t+k)^{2t+k}} \, . 
$$ 
The asymptotics in this case appear similar to those for 
the binomial coefficient $\binom{2t+2k}{k}$.  As opposed to the 
situation with Dyck paths, one may check that $f_t (k)$ is not 
equal to a binomial coefficient.

\setcounter{equation}{0}
\section{Discussion of other methods} 
\label{ss:compare}

\subsection{GF-sequence methods}
\label{ss:GF-sequence}

The benchmark work in the area of multivariate asymptotics is still
the 1983 article of Bender and Richmond~\cite{BeRi1983}. Their main result
is a local central limit theorem~\cite[Corollary~2]{BeRi1983} with the
exact same conclusion as Theorem~\ref{th:LCLT}.  Their hypotheses are:
\begin{enumerate}
\romenumi
\item $a_\rr \geq 0$;
\item $F$ has an algebraic singularity of order $q \notin \{ 0 , -1 , -2 , 
\ldots \}$ on the graph of a function $z_d = g(z_1 , \ldots , z_{d-1})$;
\item $F$ is analytic and bounded in a larger polydisk, if one excludes
a neighborhood of $\Imag (z_j) = 0$ for each $j$;
\item $B$ is nonsingular.
\end{enumerate}

Comparing this to the results presented herein, we find both
methodological and phenomenological differences.  One methodological
difference is that they view a $d$-variate generating function $F$ as
a sequence $\{ F_n \}$ of $(d-1)$-variate generating functions.  Their
main result on coefficients of $F$ is derived as a corollary of a
result~\cite[Theorem~2]{BeRi1983} on sequences satisfying $F_n \sim C_n g
h^n$ for some appropriately smooth $g$ and $h$. As we have remarked,
this approach is natural for some but not all applications, and leads
to some asymmetry in the hypotheses and conclusions.

A more important methodological difference is that while we always
work in the analytic category, Bender and Richmond use a blend of
analytic and smooth techniques\footnote{Analyticity is used only once,
when they rotate the quadratic form $B$.}. This manifests itself in
the hypotheses: where we require meromorphicity in a slightly enlarged
polydisk, they require that the function $g$ be in $C^3$, that the
residue $(1 - z / g)^q F$ be in $C^0$, and that $F$ be analytic away
from the real coordinate planes.  While our hypotheses are stronger in
this regard, we know of no applications where their assumptions hold
without $(1 - z/g)^q F$ being analytic.  Bender and Richmond gain
generality by allowing $q$ to be nonintegral.  This is further
exploited by Gao and Richmond, where the singularity is allowed to be
algebraico-logarithmic~\cite[Corollary~3]{GaRi1992}. On the other
hand, their methods entail estimates and therefore cannot handle
cancellation, such as occurs when $a_\rr$ have mixed signs and the
dominant singularity lies beyond the domain of convergence of $F$ (see
subsections~\ref{ss:edge} and~\ref{ss:riordan}).

Phenomenologically, there is a significant difference in generality
between our methods and those of Bender, Richmond, Gao, \emph{et al}.
Their results govern only the case where a local central limit theorem
holds -- indeed it seems they are interested mainly or only in this
case.  Other behaviors of interest, which have been analyzed in the
literature by various means, include Airy-type limits (see
Section~\ref{ss:core}), polynomial growth (see
Section~\ref{ss:integer}), and elliptic-type limits.  Central limit
behavior results from smooth points $\zz (\rr)$ with nondegenerate
quadratic approximations to $h_\rr (\zz (\rr))$, while Airy-type
limits result from degenerate quadratic approximations, polynomial
growth or corrections result when $\zz (\rr)$ is a multiple point, and
elliptic-type limits result from bad points.  By restricting our
exposition to the simplest cases, we have stayed mainly within the
smooth point case, where $\{ a_\rr \}$ obeys a LCLT and the methods of
Bender \emph{et al.} apply, in most cases equally well as the results
in Section~\ref{ss:smooth}.  But the advantage of multivariate
analytic methods is that they can in principle be applied to any case
in which $F$ is meromorphic or algebraic.  Thus, in addition to the
further generality covered in Section~\ref{ss:multiple} of this paper,
the same general method has been used to produce Airy limit
results~\cite{PeWi2002,Llad2003}, and is being applied to algebraic
generating functions and meromorphic functions with bad point
singularities, the simplest of which are quadratic cones. While this
work is not yet published, it appears that it will provide another
proof of the elliptic limit results for tiling statistics on the Aztec
Diamond~\cite{CEP1996} that is capable of generalizing to any quadratic
cone singularity.  This would prove similar behavior in two cases
where such behavior is only conjectured, namely cube
groves~\cite{PeSp2005} and quantum random walks (Chris Moore,
personal communication), as well as unifying these results with the
analyses of the coefficients of the generating function
$1/(1-x-y-z+4xyz)$ that arise in the study of super ballot
numbers~\cite{Gess1992} and Laguerre
polynomials~\cite{GRZ1983}.

\subsection{The diagonal method}
\label{ss:diagonal}

There is a third method for obtaining multivariate asymptotics that
deserves mention, namely the so-called ``diagonal method''.  This
derives a univariate generating function $f(z) = \sum_n a_{n,n} z^n$
for the main diagonal of a bivariate generating function $F(x,y) =
\sum_{rs} a_{rs} x^r y^s$.  The asymptotics of $a_{n,n}$ may then be
read off by standard univariate means from the function $f$. The
method may be adapted to compute a generating function for the
coefficients $a_{np,nq}$ along any line of rational slope. This
method, long known in various literatures, entered the combinatorics
literature in~\cite{Furs1967,HaKl1971}; our
exposition is taken from~[Section 6.3]\cite{Stan1999}.

While this elegant method produces an actual generating function,
which is more informative than the diagonal asymptotics, its scope is
quite limited.  First, while asymptotics in any rational direction may
be obtained, the complexity of the computation of $a_{np,nq}$
increases with $p$ and $q$.  Thus there is no continuity of
complexity, and no way to obtain uniform asymptotics or asymptotics in
irrational directions.  Secondly, the result is strictly bivariate.
Thirdly, even when the diagonal method may be applied, the computation
is typically very unwieldy.

As an example, consider the generating function
$$ 
F(x, y) = \sum_{m,n} a_{mn} x^m y^n = 
   \frac{1 + xy + x^2 y^2}{1 - x - y + xy - x^2 y^2}
$$
which~\cite{MuZa2002} enumerates binary words without zig-zags
(a zigzag is defined to be a subword $010$ or $101$ --- the
terminology comes from the usual correspondence of such words with
Dyck paths, where $0, 1$ respectively correspond to the steps $(1, 1),
(1, -1)$). Here $m, n$ respectively denote the number of $0$'s and
$1$'s in the word. The main diagonal enumerates zigzag-free words with
an equal number of $0$'s and $1$'s.  The solutions to $\denom = 0, x
\denom_x = y \denom_y$ are given by $x = y = 1/\phi$, $\phi = (1 +
\sqrt{5})/2$, and $x = y = (1 + \sqrt{3} i)/2$. Thus $\contrib$ is a
singleton $\{ (1/\phi , 1/\phi ) \}$ and the first order asymptotic is
readily computed to be
$$ 
a_{nn} \sim \phi^{2n} \frac{2}{\sqrt{n\pi\sqrt{5}}} \, .  
$$

The computation for any other diagonal is analogous, with the same
amount of computational effort, and the asymptotics are uniform over
any compact subset of directions keeping away from the coordinate
axes.

To obtain the same result via the diagonal method requires the
following steps. For each fixed $t$ near $0$, we compute the integral 
$$ 
D(z) := \sum_n a_{nn} z^n = \frac{1}{2\pi i} 
   \int_{\mathcal{C}_t} F(z, t/z) \, \frac{dz}{z} 
   = \frac{1}{2\pi i} \int_{\mathcal{C}_t} 
   \frac{1 + t + t^2}{-z^2 + (1 + t - t^2) z - t} 
$$
where the contour is a circle that encloses all the poles of $F(z,
t/z)/z$ satisfying $z(t) \to 0$ as $t \to 0$. Since $F(z, 0)/z$ has a
single simple pole at $z = 0$, the same is true of $F(z, t/z)/z$ for
sufficiently small $t$.  In this simple example  we can explicitly
solve for the pole $z(t)$ and compute its residue so that we obtain
the result
$$ 
\sum_n a_{nn} z^n = \sqrt{\frac{1 + z + z^2}{1 - 3z + z^2}} \, .  
$$
In other cases, after some manipulation we obtain $\sum_n a_{nn} z^n$
implicitly as the solution of an algebraic equation.  We are then
faced with the problem of extracting asymptotics, which can probably
be done using univariate techniques. However, we have left the realm
of meromorphic series and this can complicate matters.  In the example
above, the branching occurs outside the domain of convergence, and the
asymptotics are controlled by the dominant pole at $\phi^{-2}$ (the
minimal zero of the denominator). Thus one obtains the same asymptotic
as above, after some effort.

If we want to repeat this computation with $\sum_n a_{p n, qn} z^n$,
we are required to find all small poles of the function $F(z^q , t /
z^p)$.  It is unlikely that these may be found explicitly, which
complicates the task of finding which ones go to zero and computing
the residues there.

Finally, another serious problem faced by the diagonal method is that
while the diagonal of a rational series in $d = 2$ variables is always
algebraic, a fact which can itself be proved by the diagonal method,
in $d \geq 3$ variables, the diagonal of a rational series must be
D-finite but may not be algebraic~\cite{Lips1988}. Thus a
description of the diagonal generating function is more challenging. 
 For example, consider the
generating function $F(x, y, z) = (1 - x - y - z)^{-1} = \sum a_{rst}
x^r y^s z^t$, whose diagonal coefficient $a_{n,n,n}$ is the
multinomial coefficient $\binom{3n}{n,n,n}$.  This is known not to be
algebraic, since its asymptotic leading term $C \alpha^n n^{-1}$ is
not of the right form for an algebraic function.  It is completely
routine to derive this asymptotic using the methods of the present
article, but any method that relies on an exact description of the
diagonal will clearly require substantial extra work.

\section*{Acknowledgment}

Many thanks to Manuel Lladser and Alex Raichev for a number of comments on 
preliminary drafts.  

\bibliographystyle{alpha}

\bibliography{Mgf}

\end{document}